\documentclass[DIV=12]{scrartcl}
\usepackage[sb]{libertine} 
\usepackage[T1]{fontenc}
\usepackage{textcomp}
\usepackage[varqu,varl]{zi4}
\usepackage[amsthm,upint]{libertinust1math} 
\usepackage[scr=boondoxo,bb=boondox]{mathalpha} 
\usepackage{bm}
\newcommand{\symbf}[1]{\bm{#1}}
\newcommand{\symcal}[1]{\mathcal{#1}}
\newcommand{\symbfcal}[1]{\boldsymbol{\mathcal #1}}
\usepackage[article]{mymath-conf-pdflatex} \usepackage{booktabs}
\author{ Nils Margenberg \thanks{Helmut Schmidt University, Faculty of
    Mechanical and Civil Engineering, Holstenhofweg 85, 22043 Hamburg, Germany}
  \thanks{Corresponding author: \href{mailto:margenbn@hsu-hh.de}{\texttt{margenbn@hsu-hh.de}}} \and
  Franz X. Kärtner \thanks{Center for Free Electron Laser Science (CFEL),
    Deutsches Elektronen-Synchrotron (DESY) \& Department of Physics, University
    of Hamburg, Notkestraße 85, 22607 Hamburg, Germany} \thanks{The Hamburg Center for Ultrafast
    Imaging, Luruper Chaussee 149, 22761 Hamburg, Germany} \and Markus
  Bause\footnotemark[1] } \usepackage{placeins,multicol} \usepackage{tikz}
\usetikzlibrary{3d,arrows.meta,decorations.markings,decorations.pathreplacing,shapes,
  positioning, fit, quotes, angles, fit, calc, patterns,fadings,external}
\usepackage{colortbl} \date{} \title{Optimal Dirichlet Boundary Control by
  Fourier Neural Operators Applied to Nonlinear Optics}
\AtEveryBibitem{%
  \clearfield{urlyear}%
  \clearfield{urldate}%
}%
\AtEveryCitekey{%
  \clearfield{urlyear}%
  \clearfield{urldate}%
}%
\newcommand{\cg}{\symcal{F}}
\usepackage{biblatex}
\begin{document}
\maketitle
\vspace{-7ex}
\begin{abstract}
  \small
  We present an approach for solving optimal Dirichlet boundary control problems
  of nonlinear optics by using deep learning. For computing high resolution
  approximations of the solution to the nonlinear wave model, we propose higher
  order space-time finite element methods in combination with collocation
  techniques. Thereby, \(C^{l}\)-regularity in time of the global discrete is
  ensured. The resulting simulation data is used to train solution operators
  that effectively leverage the higher regularity of the training data. The
  solution operator is represented by Fourier Neural Operators and Gated
  Recurrent Units and can be used as the forward solver in the optimal Dirichlet
  boundary control problem.

  The proposed algorithm is implemented and tested on modern high-performance
  computing platforms, with a focus on efficiency and scalability. The
  effectiveness of the approach is demonstrated on the problem of generating
  Terahertz radiation in periodically poled Lithium Niobate, where the neural
  network is used as the solver in the optimal control setting to optimize the
  parametrization of the optical input pulse and maximize the yield of
  \(0.3\,\)THz-frequency radiation.

  We exploit the periodic layering of the crystal to design the neural networks.
  The networks are trained to learn the propagation through one period of the
  layers. The recursive application of the network onto itself yields an
  approximation to the full problem. Our results indicate that the proposed
  method can achieve a significant speedup in computation time compared to
  classical methods. A comparison of our results to experimental data shows the
  potential to revolutionize the way we approach optimization problems in
  nonlinear optics.

\end{abstract}

\section{Introduction}
\label{sec:org80141a0}
\subsection{Physical Problem and Machine Learning approach}
\label{sec:Physical-Problem-and-ML}
Nonlinear optical phenomena play a fundamental role in a lot of applications,
including the development of innovative optical sources. As high-intensity
lasers become more accessible and complexity increases, the simulation of
nonlinear optical phenomena gains importance in order to achieve optimal
performance and reduce the cost and time for empirical studies.

In this work we are concerned with the generation of Terahertz (THz) radiation,
which spans the frequency range of \qtyrange{0.1}{30}{\tera\hertz}. Thus, it is
positioned between the microwave and infrared electromagnetic frequency bands.
THz radiation offers great potential for a wide range of ultrafast
spectroscopic, strong field and imaging applications. However, a persistent
challenge in current research lies in the limited availability of compact THz
sources capable of delivering both high field strength and high repetition
rates. We address this limitation through the development of machine learning
techniques to elucidate and optimize THz generation in nonlinear crystals. By
leveraging these approaches, we aim to pave the way for the next generation of
compact and efficient light sources for spectroscopic applications, thereby
enabling significant advancements in the
field~\autocite{krtnerAXSISExploringFrontiers2016}.

In this paper, we develop machine learning techniques to solve an optimal
control problem that arises in the optimization of THz radiation generation in
nonlinear crystals. Specifically, the problem can be formulated as an optimal
Dirichlet boundary control problem, which requires the repeated solution of the
forward problem. While we previously developed accurate simulation methods for a
class of nonlinear dispersive wave equations in nonlinear
optics~\autocite{margenbergAccurateSimulationTHz2023}, each simulation entails a
significant computational effort. Consequently, the solution of the forward
problem with this method is impractical for the integration into the optimal
control problem, which necessitates different approaches.
\begin{figure}[h]
  \centering
  \includegraphics[width=\linewidth]{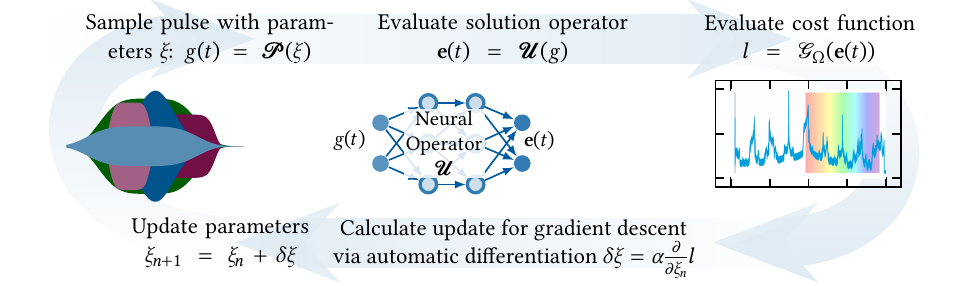}
  \caption{Sketch of the algorithm presented in this paper. We use a
    differentiable sampler, which generates the controls and a Neural Operator
    which solves the forward problem. The parametrization of the control is then
    optimized with a gradient descent method. The graph on the right shows the
    spectrum of the electric field. Here the goal is to maximize the
    \(0.3\,\)THz-frequency radiation, marked by the small bar at the
    left.}\label{fig:algsketch}
\end{figure}

The key ideas of our method are sketched in Fig.~\ref{fig:algsketch}. The first
aspect of the method we develop is the differentiability of a program that is
implemented using established Artificial Neural Network (ANN) libraries. This
paradigm is known as differentiable
programming~\autocite{rackauckasUniversalDifferentialEquations2020}.

The second main idea in our algorithm builds on the periodicity of the material
parameters. We consider a problem in nonlinear optics which involves a
periodically poled nonlinear crystal. We learn a solution operator \(\symbf{U}\)
to the forward problem in a single period of the crystal. By recursive
application of \(\symbf{U}\) onto itself we approximate the solution operator
for multiple layers. The resulting solution operator can then be integrated into
the solution of an optimal control problem.

In our work, we adopt a hybrid approach that combines classical and mature
numerical methods, specifically finite element methods, with deep learning
techniques. We use numerical methods where they have clear advantages over
machine learning approaches, while we use ANNs where numerical methods are not
feasible or efficient. Our approach to solving an optimal Dirichlet boundary
control problem exemplifies this paradigm, which is the focus of this work. In
particular, we extend our previous work
\autocite{margenbergAccurateSimulationTHz2023} on space-time finite element
methods by using higher order variational time discretizations presented in
\autocite{anselmannNumericalStudyGalerkin2020,anselmannGalerkinCollocationApproximation2020}.
The resulting finite element solution has global \(C^{1}\)-regularity in time.
We then use the resulting simulation data to train a solution operator that
effectively leverages the higher regularity of the training data. The numerical
and machine learning methods are implemented and tested on modern
high-performance computing platforms, with a focus on efficiency and
scalability.

\subsection{Related works}
\label{sec:orgbb2c17c}
Partial differential equations (PDEs) play a fundamental role in science and
engineering. They describe natural phenomena and processes in a lot of
scientific fields and provide a mathematical framework to model these phenomena.
Despite the significant advances in recent decades, challenges persist, e.\,g.\
in the context of addressing the solution of large-scale systems of nonlinear
equations.

\subsubsection*{Machine Learning for partial differential equations}
\label{sec:org6eefa6a}
Approaches for the solution of PDEs using ANNs trace back more than 20 years,
e.\,g.~\autocite{lagarisArtificialNeuralNetworks1998}. The idea is to directly
parametrize the solution to the PDE as an ANN. The network is then trained by
incorporating the differential equation, along with the boundary conditions,
into the loss function. In~\autocite{Weinan} Weinan and Yu minimize an energy
functional, resembling the variational formulation used in FEM. On the other
hand, DeepXDE~\autocite{luDeepXDEDeepLearning2021},
PINN~\autocite{raissi2019physics}, and the Deep Galerkin
Method~\autocite{sirignano2018dgm} use different approaches where the strong
residual of the PDE is minimized. This is done through collocation methods on
randomly selected points within the domain and on the boundary. Karniadakis and
Zhang proposed VPINNs \autocite{kharazmiHpVPINNsVariationalPhysicsinformed2021},
where the cost function is the variational formulation, which is optimized by
sampling test functions. A comprehensive review of PINN and related approaches
in the field of Scientific Machine Learning can be found in
\autocite{cuomoScientificMachineLearning2022}. The authors of
\autocite{mattheakisPhysicalSymmetriesEmbedded2019,chenSymplecticRecurrentNeural2020,jinSympNetsIntrinsicStructurepreserving2020,hernndezStructurepreservingNeuralNetworks2021}
develop ANNs that ensure the symplectic structure of Hamiltonian mechanics,
which improves generalization and accuracy. Based on Koopmann operator
representation, the authors
of~\autocite{ginDeepLearningModels2020,pan2020physics} train an ANN to represent
a coordinate transformation that linearizes a nonlinear PDE.

A recent approach to solving PDEs involves learning solution operators using
artificial neural networks. This technique involves approximating solution
operators using ANNs, which can potentially enable the solution of complex
problems. A significant advantage of this approach is that once the solution
operator is trained, it can be applied to other scenarios. Training ANNs is
computationally expensive, which makes PINNs and related approaches not
competitive to classical simulation methods
\autocite{grossmannCanPhysicsInformedNeural2023}. Evaluating an ANN on the other
hand is computationally cheap, making neural operators appealing: Once a
solution operator is trained, it can be generalized to other scenarios, which
only requires the evaluation of the network. Various architectures exist for
constructing these neural operators. In \autocite{lu2021learning}, Lu et al.\ construct an architecture called DeepONets by
iterating a shallow network proposed in~\autocite{chen1995universal}. This type
of network consists of a trunk network which is applied to an input function and
a branch net which is applied to an element from the domain of the operator. In
\autocite{lanthalerErrorEstimatesDeepONets2022} the authors prove an error estimate for the
DeepONet architecture. Other approaches to learn solution operators are inspired
by reduced basis methods
\autocite{bhattacharyaModelReductionNeural2021,nelsen2020random,opschoor2020deep,schwab2019deep,oleary-roseberryDerivativeinformedProjectedNeural2022,fresca2022poddlrom}.
Based on low rank decompositions the authors of~\autocite{khoo2019switchnet}
introduce an ANN with low-rank structure to approximate the inverse of
differential operators. In~\autocite{kovachkiNeuralOperatorLearning2023}, the
authors construct such a network as the tensor product of two networks, which
also carries similarities with the DeepONet
architecture~\autocite{lu2019deeponet}.

Building on~\autocite{fan2019bcr, fan2019multiscale, kashinath2020enforcing},
the Fourier Neural Operator (FNO) architecture are developed in
\autocite{liFourierNeuralOperator2020,kovachkiNeuralOperatorLearning2023}. In
\autocite{kovachkiUniversalApproximationError2021} the authors prove a universal
approximation property and error bounds. Based on FNOs, new architectures are
developed, e.\,g.\ Neural Inverse Operators
\autocite{molinaroNeuralInverseOperators2023} which are used to solve inverse
problems.

The idea of designing and interpreting ANNs using continuity becomes
increasingly popular. A notable example of this is the formulation of ResNet as
a continuous time process with respect to the depth parameter
\autocite{haber2017stable,eProposalMachineLearning2017}. See also
\autocite{antilOptimalTimeVariable2022} for an extension to adaptive
timestepping, where the timestepsize is a parameter, which can be optimized.
Similarly, works linking ANNs and dynamical systems observe that problems
arising in deep learning can be recast into optimal contol problems on
differential equations
\autocite{jiequnhanDynamicalSystemsAndOptimal2022,eMeanfieldOptimalControl2019,liuDeepLearningTheory2019,seidmanRobustDeepLearning2020,benningDeepLearningOptimal2019}.
Recent works employ deep learning techniques to address computational challenges
encountered in solving optimal control problems; The works
\autocite{weinanEmpoweringOptimalControl2022,bensoussanChapter16Machine2022} and
references therein serve as a good basis for a comprehensive survey. Existing
work is mainly concerned with stochastic control
\autocite{ruthottoMachineLearningFramework2020,carmonaConvergenceAnalysisMachine2021,carmonaConvergenceAnalysisMachine2021a}.
This is in contrast to this work, where we are concerned with Dirichlet optimal
control problems.

\subsubsection*{Space-Time Finite Element Methods}
\label{sec:orgb267e30}
We describe the numerical simulations of nonlinear optical phenomena in the
context of space-time finite element methods
\autocite{kcherVariationalSpaceTime2014}. Specifically we use time
discretization of higher order and regularity
\autocite{anselmannNumericalStudyGalerkin2020,anselmannGalerkinCollocationApproximation2020}.
Other investigations on space-time finite element methods were conducted in
\autocite{drflerSpaceTimeDiscontinuousGalerkin2016,drflerParallelAdaptiveDiscontinuous2019},
where numerical results with an adaptive algorithm are presented. Further work
relevant for electromagnetic problems is the PhD thesis
\autocite{findeisenParallelAdaptiveSpaceTime2016} and references therein.
Various alternative methods for discretizing wave equations via space-time
finite element methods exist, and they are discussed in depth in
\autocite{langerSpaceTimeMethodsApplications2019}. Notable examples include the
works of
\autocite{banjaiTrefftzPolynomialSpaceTime2017,gopalakrishnanMappedTentPitching2017},
as well as more recent developments such as those presented in
\autocite{perugiaTentPitchingTrefftzDG2020,steinbachCoerciveSpacetimeFinite2020}.
These works and their references serve as a good basis for a comprehensive
survey of recent developments in space-time discretization techniques for linear
wave equations.

The advantages of the variational time discretization include the natural
integration with the variational discretization in space and that it naturally
captures couplings and nonlinearities. These features facilitate the use of
concepts such as duality and goal oriented adaptivity
\autocite{bauseFlexibleGoalorientedAdaptivity2021}. The concepts of variational
space-time discretization also offers a unified approach to stability and error
analysis as shown in \autocite{matthiesHigherOrderVariational2011}. Furthermore,
the use of space-time FEM allow us to solve the wave equation together with the
arising ADEs in one holistic
framework~\autocite{margenbergAccurateSimulationTHz2023}. Once the formulation
is established, the methods can be extended in a generic manner. For instance,
we introduced the physical problem we are concerned with in this work
in~\autocite{margenbergAccurateSimulationTHz2023} and extend it here to the
family of Galerkin-collocation methods
\autocite{anselmannGalerkinCollocationApproximation2020}.

\section{Notation and Mathematical Problem}
\label{sec:org9aa113d}
Let \(\mathcal{D}\subset \R^d\) with \(d\in\brc{1,\,2,\,3}\) be a bounded domain
with boundary \(\partial\mathcal{D}=\Gamma_D\) and \(I=(0,\,T]\) a bounded time
interval with final time $T>0$. By \(H^m(\mathcal{D})\) we denote the Sobolev
space of \(L^2(\mathcal{D})\) functions with derivatives up to order \(m\) in
\(L^2(\mathcal{D})\). For the definition of these function spaces we refer to
\autocite{evansPartialDifferentialEquations2010}. We let
\(L \coloneq L^2(\mathcal{D})\), \(V=H^1(\mathcal{D})\) and
\(V_0=H^1_{0}(\mathcal{D})\) be the space of all \(H^1\)-functions with
vanishing trace on the Dirichlet part of the boundary \(\Gamma_D\). We denote
the \(L^{2}\)-inner product by \(\dd{\bullet}{\bullet}\). For the norms we use
\(\| \bullet \| \coloneq \| \bullet\|_{L^2(\mathcal{D})}\) and
\(\| \bullet \|_m \coloneq \| \bullet\|_{H^m(\mathcal{D})}\) for \(m\in \N\) and
\(m\geq 1\). By \(L^2(0,\,T;\,B)\), \(C([0,\,T];\,B)\) and \(C^q([0,\,T];\,B)\),
for \(q\in \N\), we denote the standard Bochner spaces of \(B\)-valued functions
for a Banach space \(B\), equipped with their natural norms. Further, for a
subinterval \(J\subseteq [0,\,T]\), we will use the notations \(L^2(J;\,B)\),
\(C^m(J;\,B)\) and \(C^0(J;\,B)\coloneq C(J;\,B)\) for the corresponding Bochner
spaces. Further, we define the function spaces, that we need below for the
variational formulation of the model equations.
\begin{definition}{Function spaces for the variational formulations}{fnspace}
  \noindent
  \begin{subequations}
    \label{eq:solnweakspace}
    \begin{align}
      \label{space} W(I)&=\brc{\vct w\in L^2(I;\,V)\suchthat \partial_t \vct w\in L^2(I;\,L)}\,,\\
      \label{space-hom} W_0(I)&=\brc{\vct w\in L^2(I;\,V_0)\suchthat \partial_t \vct w\in L^2(I;\,L)}\,,\\
      \label{nlspace} W_{\text{nl}}(I)&=\brc{\vct w\in L^2(I;\,V)\suchthat \partial_t \vct w\in L^2(I;\,L),\:\partial_t (\abs{\vct w}\vct w)\in L^2(I;\,L)}\,.
    \end{align}
  \end{subequations}
\end{definition}
\noindent
In~\eqref{nlspace} we denote by $\abs{\vct w}$ the contraction of $\vct w$, i.e.
$\abs{\vct w}=\sum_{i=1}^d\vct w_i$.
\subsubsection*{Mathematical model problem from nonlinear optics}
\label{sec:org83f7b99}
In this work we study nonlinear dispersive wave propagation, that is modeled by
the following coupled partial differential equation
(cf.~\autocite{margenbergAccurateSimulationTHz2023,abrahamConvolutionFreeMixedFiniteElement2019}).
Its physical background and application is discussed further in
Section~\ref{sec:org8d1a271}.
\begin{problem}{Nonlinear dispersive wave equation}{lor-ade}
  \begin{subequations}\label{eqwave:lor-nade2}
    \begin{align}
      \label{eqwave:lor-nade2a}
      \partial_{tt}\vct P + \Gamma_{0}\partial_{t} \vct P + \nu_{t}^2\vct P -
      (\varepsilon_{\Omega}-\varepsilon_{\omega})\nu_{t}^2 \vct E&=0 &&\text{on}\quad \mathcal{D}\times I\,,\\[.3ex]
      \label{eqwave:lor-nade2b}
      -\Delta \vct E + \varepsilon_{\omega} \partial_{tt} \vct E + (\varepsilon_{\Omega}-
      \varepsilon_{\omega})\nu_{t}^{2}\vct E- \nu_{t}^{2}\vct P -\Gamma_{0}\partial_{t}\vct P+
      \chi^{(2)}\partial_{tt}(\abs{\vct E}\vct E)&=\vct f &&\text{on}\quad \mathcal{D} \times I\,,\\[.3ex]
      \vct E(0) = \vct E_0, \quad
      \partial_t \vct E(0) = \vct E_1, \quad
      \vct P(0) = \vct P_0,\quad
      \partial_{t}\vct P(0) &= \vct P_1 &&\text{on}\quad\; \mathcal{D}\,,\\[.3ex]
      \label{eqwave:lor-nade2d} \vct E  &= g^{\vct E} &&\text{on}\quad\Gamma_D\,.
    \end{align}
  \end{subequations}
\end{problem}
By \(\vct E\) we denote the electric field, by \(\vct P\) the polarization and
\(\Gamma_{0},\,\nu_t,\,\varepsilon_{\omega},\,\varepsilon_{\Omega}\in\R_{+}\)
are material parameters. We further define
$\varepsilon_{\Delta}=\varepsilon_{\Omega}- \varepsilon_{\omega}$. The boundary
condition \(g^{\vct E}\) is a prescribed trace on \(\Gamma_D\) and \(\vct f\) is
an external force acting on the domain. To simplify the notation and enable
better numerical treatment lateron, we have already expressed
Problem~\ref{problem:lor-ade} in normalized quantities. Specifically, we have
transformed the equations and quantities using the transformation
\(\tilde{t}=c_{0}t\), where $c_0$ is the speed of light in vacuum. This
normalization is consistently applied throughout this work. Therefore, we omit
the tilde notation, as we already did in~\eqref{eqwave:lor-nade2}. For the
numerical approximation we reformulate Problem~\ref{problem:lor-ade} as a
first-order system in time; cf.~Problem~\ref{problem:lor-stm}. To this end we
introduce the auxiliary variables
\begin{subequations}\label{eq:aux}
  \begin{multicols}{2}
    \noindent
    \begin{equation}
      \vct U=\partial_{t}\vct P + \Gamma_{0} \vct P\,,
    \end{equation}
    \begin{equation}
      \vct A=\varepsilon_{\omega} \partial_{t} \vct E - \Gamma_{0} \vct P+\chi^{(2)}\partial_{t}(\abs{\vct E}\vct E)\,.
    \end{equation}
  \end{multicols}
\end{subequations}
We tacitly assume that Problem~\ref{problem:lor-ade} has a sufficiently regular,
unique solution. The proof of existence and uniqueness for the nonlinear
system~\eqref{eqwave:lor-nade2} extends beyond the scope of this work. However,
it is crucial for our subsequent mathematical arguments and formulations that
the solution to Problem~\ref{problem:lor-ade} is regular enough such that all
the mathematical arguments and formulations used below are well-defined and the
application of higher order discretization techniques becomes reasonable. This
regularity, in turn, imposes certain conditions on the data, coefficients, and
geometric properties of the domain;
cf.~\autocite{evansPartialDifferentialEquations2010}. Under the assumption of
the existence of a unique and smooth solution to~\eqref{eqwave:lor-nade2}, this
solution satisfies the following weak formulation.
\begin{problem}{Weak formulation of the nonlinear dispersive wave equation
    (\ref{eqwave:lor-nade2})}{lor-stm}
  For given data $\vct f\in L^2(I;\,L)$, boundary conditions
  $g^{\vct e}\in L^2(I;\,H^{\scriptscriptstyle 1/2}(\Gamma_D))$ and initial conditions
  $(\vct u_0,\,\vct p_0,\,\vct a_0,\,\vct e_0)\eqcolon\vct v_0\in L^3\times V$
  find
  \(\vct v \coloneq (\vct u,\,\vct p,\,\vct a,\,\vct e) \in
  \paran{W(I),\,W(I),\,W(I),\,W_{\text{nl}}(I)}\eqcolon \mat W(I)\) such that
  $\vct e\restrict{\Gamma_D}=g^{\vct e}$ and for all
  \(\vct \Phi \in \paran{W_0(I)}^4\)
  \begin{equation}
    \label{eq:variational-cont}
    \symbfcal{A}(\vct v)(\vct \Phi)=\mat F(\vct \Phi)
  \end{equation}
  is satisfied. The functional $\mat F \colon \paran{W_0(I)}^{4}\to \R$ and the
  semilinear form, which is linear in the second argument,
  $\symbfcal{A}\colon \mat W(I)\times\paran{W_0(I)}^4\to\R$ are given by
  \begin{subequations}
    \label{eq:variational-cont-def}
    \begin{align}
      \notag\symbfcal{A}(\vct v)(\vct \Phi)&\coloneq \int_0^T\dd{\partial_{t}\vct
                                             p}{\symbf \phi^0} + \Gamma_{0} \dd{\vct p}{\symbf \phi^0} - \dd{\vct u}{\vct
                                             \phi^0}  \drv t+ \dd{\vct u(0)}{\symbf \phi^0(0)}\\[3pt]
                                           &\notag\phantom{\coloneq} +\int_0^T\nu_{t}^2\dd{\vct p}{\symbf \phi^1} -
                                             \varepsilon_{\Delta}\nu_{t}^2\dd{\vct
                                             e}{\symbf \phi^1} + \dd{\partial_{t}\vct u}{\symbf \phi^1}\drv t + \dd{\vct p(0)}{\symbf \phi^1(0)}\\[3pt]
                                           &\notag\phantom{\coloneq} +\int_0^T\varepsilon_{\omega} \dd{\partial_{t} \vct e}{\symbf \phi^2} - \Gamma_{0}
                                             \dd{\vct p}{{\symbf \phi^2}}
                                             +\chi^{(2)}\dd{\partial_{t}(\abs{\vct E}\vct E)}{\symbf \phi^2} - \dd{\vct
                                             a}{\symbf \phi^2}\drv t + \dd{\vct a(0)}{\symbf \phi^2(0)}\\[3pt]
                                           &\phantom{\coloneq} +\int_0^T\dd{\symbf \nabla \vct e}{\symbf \nabla\symbf \phi^3} + \varepsilon_{\Delta}\nu_{t}^{2}\dd{\vct e}{\symbf \phi^3}
                                             -
                                             \nu_{t}^{2}\dd{\vct p}{\symbf \phi^3}+\dd{\partial_{t}\vct a}{\symbf \phi^3}\drv t+\dd{\vct e(0)}{\symbf \phi^3(0)}\,,\\
      \mat F(\symbf \phi) &\coloneq \dd{\vct u_0}{\symbf \phi^0(0)}
                          +\dd{\vct p_0}{\symbf \phi^1(0)}
                          +\dd{\vct a_0}{\symbf \phi^2(0)}
                          +\dd{\vct e_0}{\symbf \phi^3(0)}+\int_0^T\dd{\vct f}{\symbf \phi^3}\drv t\,.
    \end{align}
  \end{subequations}
\end{problem}
\noindent
We note that all integrals in~\eqref{eq:variational-cont} are well-defined in
the function space $\mat W(I)$, due to the Definition~\ref{def:fnspace}. To
obtain higher regularity of the solution, stricter assumptions $\vct f$ and
$\vct v_0$ may have to be imposed.
\begin{remark}{Weakly imposed initial conditions}{weakinitial}
  In~\eqref{eq:variational-cont-def}, the expressions $\vct w(0)$ for
  $\vct w\in\brc{\vct u,\,\vct p,\,\vct a,\,\vct e}$ are well-defined when we consider
  the continuous embedding $W\hookrightarrow C(\bar{I};\,V)$,
  cf.~\autocite[Chapter XVIII, Theorem
  1]{dautrayMathematicalAnalysisNumerical1999}. We further note that the test
  space $W_0$ is
  dense in $L^2(I;\,V_0)$, as stated in~\autocite[Chapter 2, Corollary
  2.1]{bruchhuserGoalOrientedSpaceTimeAdaptivity2022}.
\end{remark}
\noindent
Based on Remark~\ref{rem:weakinitial}, we comment on the variational
problem~\eqref{eq:variational-cont}.
\begin{itemize}\itemsep1pt \parskip0pt \parsep0pt
\item For convenience, the initial conditions of Problem~\ref{problem:lor-stm}
  are incorporated in the variational equation~\eqref{eq:variational-cont}
  through the forms~\eqref{eq:variational-cont-def}. The Sobolev embedding
  $\mat W(I)\hookrightarrow {C(\bar{I};,V)}^{4}$ ensures the well-defined
  pointwise evaluation of functions in $\mat W(I)$ within the
  forms~\eqref{eq:variational-cont-def}.
\item According to Remark~\ref{rem:weakinitial}, the test space ${W_0(I)}^4$ is
  densely embedded in the Hilbert space ${L^2(I;V_0)}^4$. This dense embedding
  is an indispensable requirement for the proper formulation of
  Problem~\ref{problem:lor-stm}.
\item The variables $(\vct u,\,\vct p,\,\vct a,\,\vct e)$ belong to the solution
  space $\mat W(I)$. Although weaker assumptions about $\vct u$, $\vct p$ and
  $\vct a$ would have been sufficient for the existence of the space-time
  integrals in~\eqref{eq:variational-cont}. However, we adopt this stronger
  assumption since we use an $H^1(\mathcal{D})$-conforming approximation for all
  variables in Section~\ref{sec:org5743c30}. This concept follows the lines
  of~\autocite{bangerthAdaptiveGalerkinFinite2010}.
\end{itemize}
Under the above-made assumptions we now define the solution operator that is
associated with Problem~\ref{problem:lor-stm} and its weak
formulation~\eqref{eq:variational-cont}.
\begin{definition}{Solution Operator}{abstract}
  Consider the nonlinear Problem~\ref{problem:lor-stm}. The solution operator
  \begin{equation}
    \label{eq:abstract-inv}
    \begin{split}
      \symbfcal{S}\colon D(\symbfcal{S})\subset L^2(I;\,H^{\scriptscriptstyle 1/2}(\Gamma_D))\times L^2(I;\,L)\times (L^{3}\times V) &\to \mat
                                                                                                                    W(I),\\
      (g^{\vct e},\,\vct f,\,\vct v_0)&\mapsto\vct v\,.
    \end{split}
  \end{equation}
  is defined by the mapping of the data $\vct f$ and the initial conditions
  $\vct v_0$ to the unique solution $\vct v$ of~\eqref{eq:variational-cont},
  such that
  \begin{equation}
    \label{eq:abstract-inv-rel}
    \symbfcal{A}(\symbfcal{S}(g^{\vct e},\,\vct f,\,\vct v_0))(\vct \Phi) = \mat
    F(\vct \Phi)\quad \forall \vct \Phi\in {(W_0(I))}^{4}\,.
  \end{equation}
\end{definition}
The domain $D(\symbfcal S)$ is supposed to be a subset of sufficiently regular
functions $\vct f,\,\vct v_0$ in
$L^2(I;\,H^{\scriptscriptstyle 1/2}(\Gamma_D))\times L^2(I;\,L)\times (L^{3}\times V)$ such
that~\eqref{eq:abstract-inv-rel} admits a unique solution with the regularity
required for the numerical approximation scheme. The goal of this work is to
approximate the operator $\symbfcal S$ by an ANN, which evaluation involves low
computational costs and thereby lets an optimal control problem subject to the
Dirichlet data of Problem~\ref{problem:lor-stm} become feasible. For the
training and validation of the ANN approximate solutions to
Problem~\ref{problem:lor-stm} with high resolution are necessary. They are
computed by space-time finite element techniques of high accuracy which are
presented in Section~\ref{sec:org5743c30}.

\subsection{Physical Background}
\label{sec:org8d1a271}
Based on~\autocite{margenbergAccurateSimulationTHz2023}, we review the
model~\eqref{eqwave:lor-nade2}, with a focus on the applications and physics of
nonlinear
optics~\autocite{newIntroductionNonlinearOptics2011,boydChapterNonlinearOptical2020}.
Nonlinear and dispersive effects arise due to the interaction of waves with
atoms or molecules in a medium. The polarization \(\vct P\) of the medium
captures these interactions at a macroscopic level. The polarization can be
developed as a power series in terms of the electric field \(\vct E\). Based on
the physical settings and materials considered in this work, it is deemed
sufficient to include only the linear and quadratic terms to accurately model
the phenomena of interest. The polarization is then given by
\begin{equation*}
  \vct P(x,\,t)=\varepsilon_{0}\paran{\chi^{(1)}\otimes \vct E(x,\,t)+\chi^{(2)}\otimes
    \vct E(x,\,t)\otimes \vct E(x,\,t)}\,,
\end{equation*}
where the electric susceptibilities
\(\chi^{(n)}\colon \R\times\mathcal{D} \to \otimes_{i=0}^{n}\C^{d}\) are
tensor-valued functions which depend on the frequency and spatial coordinate. We
assume that \(\chi^{(1)}\) and \(\chi^{(2)}\) can be simplified to scalar
functions such that
\begin{equation}
  \label{eq:1}
  \chi^{(1)} \colon \R \to \C\quad \text{and}\quad \chi^{(2)}\colon \mathcal{D}\to \R\,.
\end{equation}
We note that \(\chi^{(1)}\) doesn't depend on spatial coordinates and the
material is homogeneous w.\,r.\,t.\ to the linear susceptibility. Further, we
only consider instantaneous nonlinearities, which means that the nonlinear
susceptibilities are frequency independent. We formulate the dispersive
electromagnetic wave equation
\begin{equation}\label{eqnlwave}
  -\symbf\Delta \vct E +\partial_{tt}\varepsilon_{r}* \vct E
  +\chi^{(2)}\partial_{tt}(\abs{\vct E}{\vct E})
  =\vct f.
\end{equation}
Here \(\varepsilon_{r}\) is the relative electric permittivity for which
\(\varepsilon_{r}=n^{2}=1+\chi^{(1)}\) holds, where \(n\) is the refractive
index. A simple model introduced by Lorentz, which describes the electric
permittivity \(\varepsilon_{r}\) as a function of the frequency \(\nu\) is given
by
\begin{equation}\label{eqlorentzfreq} {\varepsilon_{r}(\nu)}=
  \varepsilon_{\omega}
  + \frac{(\varepsilon_{\Omega}-\varepsilon_{\omega})
    \nu_{t}^{2}}{\nu_{t}^{2}-\nu^{2}+\iu\Gamma_{0}\nu}\,.
\end{equation}
The physical model for~\eqref{eqlorentzfreq} is an electron bound to the nucleus
by a force governed by Hooke's law with characteristic frequency \(\nu_{t}\).
\(\Gamma_{0}\) is the damping coefficient and \(\varepsilon_{\Omega}\) and
\(\varepsilon_{\omega}\) are the low and high frequency limits of the relative
electric permittivity. In the time domain~\eqref{eqlorentzfreq} gives rise to
the convolution term \([\varepsilon_{r}(\nu) * \vct E](t)\) in~\eqref{eqnlwave}.
To avoid the computationally expensive evaluation of this convolution, we derive
an auxiliary differential equation (ADE), as given by~\eqref{eqwave:lor-nade2a}.
By substituting~\eqref{eqwave:lor-nade2a} into~\eqref{eqnlwave}, we
obtain~\eqref{eqwave:lor-nade2b}, which results in the formulation of
Problem~\ref{problem:lor-ade}.

\section{Variational Space-Time Discretization for Nonlinear Dispersive Wave
  Equations}
\label{sec:org5743c30}
In this section we present the numerical approximation scheme that we use for
highly resolved and accurate computations of solutions to the weak
form~\eqref{eq:variational-cont} of the nonlinear dispersive wave problem in
Problem~\ref{problem:lor-ade}.

The approach discretizes the continuous system~\eqref{eqwave:lor-nade2} by
enforcing differentiability in time constraints on the trial space of piecewise
polynomials in combination with variational conditions, based on the weak
formulation~\eqref{eq:variational-cont} and collocation conditions, deduced from
the strong form~\eqref{eqwave:lor-nade2}. The collocation conditions are imposed
at the end point of the subintervals of the time mesh. Due to the
differentiability in time, we will observe that the collocation conditions are
also satisfied at the initial time points of the subintervals. These schemes are
referred to as Galerkin-collocation methods, for short $\text{GCC}^s(k)$ where
$s$ denotes the differentiability with respect to the time variable and $k$ the
order of the polynomials of the trial space. Galerkin-collocation schemes have
been introduced and studied for acoustic waves
in~\autocite{anselmannNumericalStudyGalerkin2020,anselmannGalerkinCollocationApproximation2020}.
For the choice $r=k$ ($r$ being the order of approximation in space),
convergence of order $k+1$ in space and time is shown for the fully discrete
approximation of the solution and its time derivative. In our simulations
presented in Section~\ref{sec:org4ff2e73} we put $k=3$. In the numerical
investigations of Section~\ref{sec:org4ff2e73}, we will see that
Galerkin-collocation are strongly adapted to the accurate and efficient
numerical simulation of nonlinear dispersive phenomena.

The collocation conditions allow us to reduce the size of the discrete
variational test space, which leads to increased efficiency compared to standard
Galerkin-Petrov approaches, as presented
in~\autocite{kcherVariationalSpaceTime2014} for example. Galerkin-collocation
schemes lead to discrete solutions of higher order regularity in time. For
instance, by employing the \(\text{GCC}^{1}(3)\) method, the simplest scheme
from this family of time discretization techniques, we obtain solutions of
\(C^{1}\)-regularity in time, which is particularly advantageous for wave
problems. We also exploit the increased regularity in our optimal contral method
by neural networks in Section~\ref{sec:orga943742}.

For the time discretization, we split the time interval \(I\) into a sequence of
\(N\) disjoint subintervals \(I_n=(t_{n-1},\,t_n]\), \(n=1,\dots,\,N\). For a
Banach space \(B\) and \(k\in \mathbb{N}_{0}\) we define
\begin{equation}
  \label{eq:timespace}
  \mathbb{P}_{k}(I_{n};\,B)=\brc{w_{\tau_{n}}\colon I_{n}\to B\suchthat
    w_{\tau_{n}}(t)=\sum_{j=0}^kW^jt^j\;\forall t\in I_n,\:W^j\in B \;\forall j}\,.
\end{equation}
For \(r\in \mathbb{N}\) we define the finite element space that is built on the
spatial mesh as
\begin{equation}
  \label{eq:fespace}
  \symbfcal{V}_{h}=\brc{v_h\in C(\bar{\mathcal{D}})\suchthat v_h\restrict{K} \in
    \mathcal{Q}_r(K)\;\forall K \in \mathcal{T}_h}\,, \qquad\symbfcal{V}_{h,\,0}=\symbfcal{V}_{h}\cap V_0\,,
\end{equation}
where \(\mathcal{Q}_r(K)\) is the space defined by the reference mapping of
polynomials on the reference element with maximum degree \(r\) in each variable.
From now on we choose the piecewise polynomial degrees in~\eqref{eq:timespace}
and~\eqref{eq:fespace} to $k=3$ and $r=3$. The trial and test space for our
discrete problem are then defined by
\begin{equation}
  \label{eq:disc-spaces}
  \begin{alignedat}{1}
    \symbf{X}_{\tau,\,h}=\brc{w\in C^{1}(\bar{I}; \symbfcal V_h) \suchthat w\restrict{I_n} \in \mathbb{P}_{3}(I_{n};\,\symbfcal V_h) \;\forall n=1,\dots,\,N}\,,\\
    \symbf{Y}_{\tau,\,h}=\brc{w\in L^2(I; \symbfcal V_{h,\,0}) \suchthat w\restrict{I_n} \in \mathbb{P}_{0}(I_{n};\,\symbfcal V_{h,\,0}) \;\forall n=1,\dots,\,N}\,.
  \end{alignedat}
\end{equation}
We impose global $C^1$-regularity on $\symbf{X}_{\tau,\,h}$, which corresponds
to a spline-type discretization in time. We chose the global time-discrete space
of piecewise constant functions as $\symbf{Y}_{\tau,\,h}$. Thereby, we need to
fix additional degrees of freedom in order to ensure solvability. To this end,
we combine the $C^1$-regularity constraints with the strong form of the
equations at the endpoints of each subinterval $I_n$. Then, collocation
conditions are a result of the imposed global $C^1$-regularity. This is
different from~\autocite{becherVariationalTimeDiscretizations2021}, where the
collocation conditions are imposed, which then imply the $C^1$-regularity. The
different construction is due to the nonlinear character of the system. For
simplicity regarding the prescription of inhomogeneous boundary conditions we
make the following assumption.
\begin{assumption}{Inhomogeneous Dirichlet Boundary conditions}{inhom}
  We impose an implicit restriction on the set of admissible boundary conditions
  $g^{\vct e}$. We assume that there exists a function $g_{\tau,\,h}$ in
  $C^1(\bar{I};\,\symbfcal{V}_h)$ such that
  \[g_{\tau,\,h}^{\vct e}\restrict{\Gamma_D}=g^{\vct e} \quad \forall t \in
    \bar{I}\,.\] For prescribing more general boundary conditions suitable
  interpolation operators applied to the boundary values are required. For
  brevity and since this is a standard technique, it is not considered here.
\end{assumption}

We let
$(\vct e_{0,\,h},\, \vct a_{0,\,h},\, \vct p_{0,\,h},\, \vct u_{0,\,h})\eqcolon
\vct v_{0,\,h}\in \paran{\symbfcal V_h}^4$, which are appropriate finite element
approximations of the initial values \(\vct v_0\). Here, we use interpolation in
\(\symbfcal V_h\).
We introduce
\begin{equation}
\label{eq:eval-at-tn}
\partial_t^i\vct w_{n,\,h}=\partial_t^i\vct w_{\tau,\,h}(t_n)\,,
\end{equation}
and discretize Problem~\ref{problem:lor-stm} with the \(\text{GCC}^{1}(3)\)
method. From the local problems~\ref{problem:lor-gcc-local} we derive the
following global in time fully discrete formulation.
\begin{problem}{$\boldsymbol{C^{1}}$-regular in time Galerkin-collocation scheme
    for (\ref{eq:variational-cont})}{lor-gcc}
  \allowdisplaybreaks%
  For given data and boundary conditions
  $\vct f_{\tau,\,h},\:g_{\tau,\,h}^{\vct e}\in
  C^1(\bar{I};\,\symbfcal{V}_{h,\,0})$, find $\vct u_{\tau,\,h}$,
  $\vct p_{\tau,\,h}$, $\vct a_{\tau,\,h}$ and $\vct e_{\tau,\,h}$ such that
  $\vct e_{\tau,\,h}=g_{\tau,\,h}^{\vct e}\:\text{on}\:\bar{I}\times\Gamma_D$ and
  for all
  \((\symbf \phi_{1,\,h}^0,\dots,\symbf \phi_{1,\,h}^{11},\dots,\,\symbf  \phi_{N,\,h}^0,\dots,\,\symbf \phi_{N,\,h}^{11},\,\symbf  \phi_{\tau,\,h}^0,\dots,\,\symbf \phi_{\tau,\,h}^3)\eqcolon \symbf \phi_{\tau,\,h}
  \in {\symbfcal V_{h,\,0}}^{12N}\times\mat Y_{\tau,\,h}^{4}\)
  \begin{equation}
    \label{eq:abstract-discrete-form}
    \symbfcal{A}_{\tau,\,h}(\vct v_{\tau,\,h})(\vct
    \Phi_{\tau,\,h})=\mat F_{\tau,\,h}(\symbf \phi_{\tau,\,h})\,,
  \end{equation}
  is satisfied, where $\mat F \colon \mat Y_{\tau,\,h}^{4}\to \R$ and
  $\symbfcal{A}\colon \mat X_{\tau,\,h}^{4}\times \paran{{\symbfcal
      V_h}^{12N}\times\mat Y_{\tau,\,h}^{4}}\to\R$ are given by
  \begin{subequations}\label{eq:wavephy-variational-all}
    \begin{align}
      \notag\symbfcal{A}_{\tau,\,h}(\vct v_{\tau,\,h})(\symbf \phi_{\tau,\,h})&\coloneq
                                                                              \int_0^T\dd{\partial_{t}\vct p_{\tau,\,h}}{\symbf \phi_{\tau,\,h}^{0}} +
                                                                              \Gamma_{0} \dd{\vct p_{\tau,\,h}}{\symbf \phi_{\tau,\,h}^{0}} - \dd{\vct
                                                                              u_{\tau,\,h}}{\symbf \phi_{\tau,\,h}^{0}} \drv t\\
      \notag&+\int_0^T\nu_{t}^2\dd{\vct p_{\tau,\,h}}{\symbf              \phi_{\tau,\,h}^{1}} -\varepsilon_\Delta\nu_{t}^2\dd{\vct e_{\tau,\,h}}{\symbf              \phi_{\tau,\,h}^{1}} + \dd{\partial_{t}\vct u_{\tau,\,h}}{\symbf              \phi_{\tau,\,h}^{1}}\drv t\\
      \notag&+\int_0^T\varepsilon_{\omega} \dd{\partial_{t} \vct
              e_{\tau,\,h}}{\symbf \phi_{\tau,\,h}^{2}} - \Gamma_{0}\dd{\vct
              p_{\tau,\,h}}{{\symbf \phi_{\tau,\,h}^{2}}}+\chi^{(2)}\dd{\partial_{t}(\abs{\vct
              e_{\tau,\,h}}\vct e_{\tau,\,h})}{\symbf \phi_{\tau,\,h}^{2}} - \dd{\vct
              a_{\tau,\,h}}{\symbf \phi_{\tau,\,h}^{2}}\drv t\\
      \label{eq:varcon}&+\int_0^T\dd{\symbf \nabla \vct
              e_{\tau,\,h}}{\symbf \nabla\symbf \phi_{\tau,\,h}^{3}} +
              \varepsilon_{\Delta}\nu_{t}^{2}\dd{\vct e_{\tau,\,h}}{\symbf              \phi_{\tau,\,h}^{3}} - \nu_{t}^{2}\dd{\vct p_{\tau,\,h}}{\symbf              \phi_{\tau,\,h}^{3}}+\dd{\partial_{t}\vct a_{\tau,\,h}}{\symbf              \phi_{\tau,\,h}^{3}}\drv t\\
      \notag&\!\!\!\!\!\!\!\!\!\!\!\!\!\!\!\!\!\!\!\!\!\!\!\!\!\!\!\!\!\!
              +\dd{\vct u_{\tau,\,h}(0)}{\symbf \phi_{\tau,\,h}^0(0)}
              +\dd{\vct p_{\tau,\,h}(0)}{\symbf \phi_{\tau,\,h}^1(0)}
              +\dd{\vct a_{\tau,\,h}(0)}{\symbf \phi_{\tau,\,h}^2(0)}
              +\dd{\vct e_{\tau,\,h}(0)}{\symbf \phi_{\tau,\,h}^3(0)}\\
      \label{eq:initcon}&\!\!\!\!\!\!\!\!\!\!\!\!\!\!\!\!\!\!\!\!\!\!\!\!\!\!\!\!\!\!
                          +\dd{\partial_{t}\vct u_{\tau,\,h}(0)}{\symbf \phi_{\tau,\,h}^4(0)}
                          +\dd{\partial_{t}\vct p_{\tau,\,h}(0)}{\symbf \phi_{\tau,\,h}^5(0)}
                          +\dd{\partial_{t}\vct a_{\tau,\,h}(0)}{\symbf \phi_{\tau,\,h}^6(0)}
                          +\dd{\partial_{t}\vct e_{\tau,\,h}(0)}{\symbf \phi_{\tau,\,h}^7(0)}\\
      \notag&+\sum_{n=1}^{N}\Biggr(\dd{\partial_{t}\vct p_{n,\,h}}{\symbf \phi_{n,\,h}^{8}} +
              \Gamma_{0} \dd{\vct p_{n,\,h}}{\symbf \phi_{n,\,h}^{8}} - \dd{\vct
              u_{n,\,h}}{\symbf \phi_{n,\,h}^{8}}\\
      \notag&\quad+\nu_{t}^2\dd{\vct p_{n,\,h}}{\symbf              \phi_{n,\,h}^{9}} -\varepsilon_\Delta\nu_{t}^2\dd{\vct e_{n,\,h}}{\symbf              \phi_{n,\,h}^{9}} + \dd{\partial_{t}\vct u_{n,\,h}}{\symbf              \phi_{n,\,h}^{9}}\\
      \notag&\quad+\varepsilon_{\omega} \dd{\partial_{t} \vct
              e_{n,\,h}}{\symbf \phi_{n,\,h}^{10}} - \Gamma_{0}\dd{\vct
              p_{n,\,h}}{{\symbf \phi_{n,\,h}^{10}}}+\chi^{(2)}\dd{\partial_{t}(\abs{\vct
              e_{n,\,h}}\vct e_{n,\,h})}{\symbf \phi_{n,\,h}^{10}} - \dd{\vct
              a_{n,\,h}}{\symbf \phi_{n,\,h}^{10}}\\
      &\quad+\dd{\symbf \nabla \vct
              e_{n,\,h}}{\symbf \nabla\symbf \phi_{n,\,h}^{11}} +
              \varepsilon_{\Delta}\nu_{t}^{2}\dd{\vct e_{n,\,h}}{\symbf              \phi_{n,\,h}^{11}} - \nu_{t}^{2}\dd{\vct p_{n,\,h}}{\symbf              \phi_{n,\,h}^{11}}+\dd{\partial_{t}\vct a_{n,\,h}}{\symbf              \phi_{n,\,h}^{11}}
              \Biggl)\,,\label{eq:wavephy-variational-lhs}\\
      \mat F_{\tau,\,h} (\symbf \phi_{\tau,\,h})&\coloneq\int_0^T\dd{\vct
                                                      f_{\tau,\,h}}{\symbf \phi_{\tau,\,h}^{3}}\drv
                                                      t+ \sum_{n=1}^{N}\dd{\vct
                                                      f_{\tau,\,h}(t_n)}{\phi_{n,\,h}^{3}}\\
      \notag &+\dd{\vct u_{0,\,h}}{\symbf               \phi_{\tau,\,h}^0(0)}+\dd{\vct p_{0,\,h}}{\symbf \phi_{\tau,\,h}^1(0)}+\dd{\vct
               a_{0,\,h}}{\symbf \phi_{\tau,\,h}^2(0)}+\dd{\vct e_{0,\,h}}{\symbf               \phi_{\tau,\,h}^3(0)}\\
      &\!\!\!\!\!\!\!\!\!\!\!\!+\dd{\partial_t\vct u_{0,\,h}}{\symbf \phi_{\tau,\,h}^4(0)}+\dd{\partial_t\vct p_{0,\,h}}{\symbf \phi_{\tau,\,h}^5(0)}+\dd{\partial_t\vct a_{0,\,h}}{\symbf \phi_{\tau,\,h}^6(0)}+\dd{\partial_t\vct e_{0,\,h}}{\symbf \phi_{\tau,\,h}^7(0)}\,.\label{eq:wavephy-variational-rhs}
    \end{align}
  \end{subequations}
\end{problem}
In our implementation, we use a local test basis supported on the subintervals
$I_n$ in Problem~\ref{problem:lor-gcc}. This leads to a time marching scheme
with the local Problem~\ref{problem:lor-gcc-local} to be solved in each of the
time steps.
We comment on the fully discrete Problem~\ref{problem:lor-gcc}.
\begin{itemize}\itemsep1pt \parskip0pt \parsep0pt
\item In constrast to \autocite{anselmannGalerkinCollocationApproximation2020},
  collocation conditions are a result of the imposed global $C^1$-regularity. As
  already mentioned in \autocite[Remark
  3.4]{anselmannGalerkinCollocationApproximation2020} the approach of imposing
  global $C^1$-regularity is also valid. We also show
  this in detail in Appendix~\ref{sec:org85a3097}.
\item After breaking Problem~\ref{problem:lor-gcc} into local problems, we
  can put the equations of the proposed \(\text{GCC}^1(3)\) approach in their
  algebraic forms (cf.~Appendix~\ref{sec:org85a3097}) and get a nonlinear system
  of equations. The common approach of handling the nonlinear problem is a
  linearization by means of Newton's method. In every Newton step we have to
  solve a linear system of equations, of which we give a detailed description
  in the Appendix~\ref{sec:org85a3097}.
\item Hermite polynomials are ideal for wave problems, particularly those with
  high frequencies. Moreover, they offer significant advantages for the
  numerical solution of the nonlinear wave equations~\eqref{eqwave:lor-nade2} by
  reducing the computational cost of assembling matrices and residuals in
  Newton's method. These advantages result from the sparse structure of the
  nonlinear term given by Hermite polynomials as trial functions:
  \begin{align*}
    \int_{I_n}  \chi^{(2)}\dd{\partial_{t}(\abs{\vct e}\vct e)}{\zeta_1}\drv
    t &=\chi^{(2)}  {(|\vct e_0|\;|\vct e_1|\;|\vct e_2|\;|\vct e_3|)}^{\top}%
        \underbrace{\int_{I_n}{(\partial_{t}(\xi_i\xi_j))}_{i,\,j=0,\dots,\,3}\drv t}_{
        {= \diag\paran{(-1,\,0,\,1,\,0)}\in \R^{4\times 4}}}\,%
        (\vct e_0\;\vct e_1\;\vct e_2\;\vct e_3)\\
      &=\chi^{(2)} \abs{\vct e_2}\vct e_2-\chi^{(2)}\abs{\vct e_0}\vct e_0\,,
  \end{align*}
  where $\brc{\vct e_{i}}_{i=0}^3$ denote the coefficient functions of the $i$-th
  time basis function in $\symbfcal{V}_h$. This also applies to third order
  nonlinearities $\chi^{(3)}$ with two non-vanishing terms.
  \begin{align*}
    \int_{I_n}  \chi^{(3)}\dd{\partial_{t}(\abs{\vct e}^{2}\vct e)}{\zeta_1}\drv
    t &=\chi^{(3)}  {(|\vct e_0|^{2}\;|\vct e_1|^{2}\;|\vct e_2|^{2}\;|\vct e_3|^{2})}^{\top}%
        \int_{I_n}{(\partial_{t}(\xi_i\xi_j\xi_k))}_{i,\,j,\,k=0,\dots,\,3}
        (\vct e_0\;\vct e_1\;\vct e_2\;\vct e_3)\\
      &=\chi^{(3)} \abs{\vct e_2}^{2}\vct e_2-\chi^{(3)}\abs{\vct e_0}^{2}\vct e_0.
  \end{align*}
\end{itemize}
We consider the abstract space-time discrete
form~\eqref{eq:abstract-discrete-form}.
The global formulation puts the work in this section in context
to the abstract problem introduced in Definition~\ref{def:abstract} and,
together with an analogous formulation of the solution operator, becomes useful in
the next section.
\begin{definition}{Discrete Solution Operator}{abstract-discrete}
  Consider Problem~\ref{problem:lor-gcc} given in
  variational formulation. Then the solution operator $\symbfcal{S}_{\tau,\,h}$
  to~\eqref{eq:abstract-discrete-form}, which maps $g_{\tau,\,h}^{\vct e}$, $\vct f_{\tau,\,h}$ and the initial conditions $\vct v_{0,\,h}$ to the solution $\vct v_{\tau,\,h}$ is defined through
  \begin{equation}
    \label{eq:abstract-inv-disc}
    \begin{split}
    \symbfcal{S}_{\tau,\,h}\colon D(\symbfcal{S}_{\tau,\,h})\subset
      C^1(\bar{I};\,\symbfcal{V}_h)\times C^1(\bar{I};\,\symbfcal{V}_{h,\,0})\times\symbfcal{V}_h^4 &\to \mat X_{\tau,\,h}^4,\\
      (g_{\tau,\,h}^{\vct e},\,\vct f_{\tau,\,h},\,\vct v_{0,\,h})&\mapsto\vct v\,.
    \end{split}
  \end{equation}
\end{definition}
In order to ensure well-posedness we assume that \(\symbfcal{S}_{\tau,\,h}\) is
a bijection. Then, for given data $\vct f_{\tau,\,h}$ we find a unique solution
\((\vct e_{\tau,\,h},\,\vct a_{\tau,\,h},\,\vct p_{\tau,\,h},\,\vct
u_{\tau,\,h}) = \vct v_{\tau,\,h}\in \mat X_{\tau,\,h}^4\) which satisfies
\(\symbfcal{S}_{\tau,\,h}(g_{\tau,\,h}^{\vct e},\,\vct f_{\tau,\,h},\,\vct v_{0,\,h}) = \vct v_{\tau,\,h}\). We note that
\(\symbfcal{S}_{\tau,\,h}\), \(\vct v_{\tau,\,h}\), \(g_{\tau,\,h}^{\vct e}\) and \(f_{\tau,\,h}\)
approximate \(\symbfcal{S}\), \(\vct v\) and \(\vct f\)
in~\eqref{eq:abstract-inv}. In the next section we introduce two types of ANNs,
which we consider for training a discrete solution operator
\begin{equation}
\label{eq:idea}
\symbfcal{U}\approx\symbfcal{S}_{\tau,\,h}\,.
\end{equation}
$\symbfcal{U}$ is trained with accurate approximations \(\vct v_{\tau,\,h}\)
obtained by numerical solutions and is subsequently used to solve an optimal
control problem.

\section{Artificial Neural Networks}
\label{sec:ANN}
Neural networks exist in various types. In this section we briefly review the
architecture of the neural networks that we use below to learn the discrete
solution operator defined in Problem~\ref{def:abstract-discrete}. In
Section~\ref{sec:orga943742}, the neural networks are then applied to accelerate
optimization processes for Dirichlet boundary control of the pump pulse for
terahertz generation.
\subsection{Fourier Neural Operators (FNO)}
\label{sec:orgbe9b164}
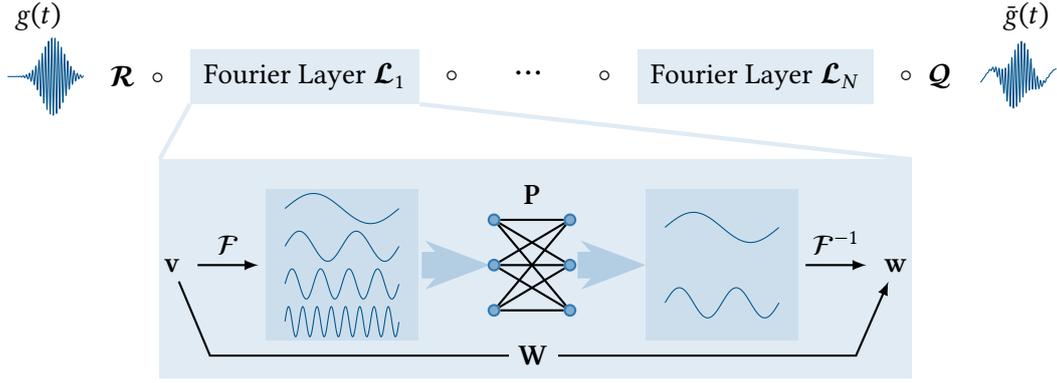
\begin{figure}
  \centering
  \begin{tikzpicture}[domain=-0.5:0,>=latex,scale=2]
    \tikzstyle{entity} = [rectangle,%
    anchor=south west,%
    every entity]%
    \tikzstyle{every entity} = [fill=helmholtzblue!12, align=left,inner
    sep=5pt] %
    \tikzstyle{every pin edge}=[<-]%
    \tikzstyle{neuron}=[circle,fill=helmholtzblue!80,minimum size=4pt,inner
    sep=0pt]%
    \tikzstyle{hidneuron}=[neuron,
    fill=helmholtzblue!50,draw=helmholtzblue!80,thick];%
    \begin{scope}[yshift=1.25cm,node distance=.2cm]
      \node[entity,anchor=west] (l1) at (0,0) {Fourier Layer $\symbfcal
        L_{1}$};%
      \node[entity,anchor=east] (ln) at (4.5,0) {Fourier Layer
        $\symbfcal L_{N}$};%
      \node[right = of l1] (c1) {\Large $\circ$};%
      \node[left = of ln] (cn) {\Large $\circ$};%
      \node[left = of l1] (q) {$\symbfcal{R}\;$ {\Large $\circ$}}; \node[right =
      of ln] (r) {{\Large $\circ$} $\;\symbfcal{Q}$}; \node (cd) at
      ($(c1)!0.5!(cn)$) {\Large $\cdots$}; \draw[color=myblue, samples=2000,
      xshift=5.7cm] plot[smooth]
      (\x,{.2*exp(-((\x+0.25)*(\x+0.25))*100)*cos(300*\x r)+.05*cos(30*\x
        r)}); %
      \node (b) at (5.5,0.4) {$\bar{g}(t)$}; %
      \draw[color=myblue, samples=2000,xshift=-0.7cm]
      plot[smooth](\x,{.2*exp(-((\x+0.15)*(\x+0.25))*100)*cos(300*\x r)}); %
      \node (a) at (-1,0.4) {$g(t)$}; %
    \end{scope}
    \fill[helmholtzblue!12] (-0.2,-0.75) node[anchor=north west]{} rectangle
    (4.75,0.7);%
    \draw[line join=bevel,ultra thick,helmholtzblue!12] (l1.south west) --
    (-0.2,0.7); \draw[line join=bevel,ultra thick,helmholtzblue!12] (l1.south
    east) -- (4.75,0.7); \node[anchor=east] (aa) at (-0,0) {$\vct v$}; %
    \node[anchor=west] (bb) at (4.5,0) {$\vct w$}; %
    \draw[->,looseness=2,thick] (aa) -- ++ (0.25,-0.6) --
    node[fill=helmholtzblue!12,pos=0.5]{$\mat W$} ++ (4.25,0) -- (bb);%
    \draw[->,thick] (0.05,0) -- ++ (.4,0) node[above, midway]{$\cg$};
    \fill[helmholtzblue!20] (0.5,-0.5) node[anchor=north west]{} rectangle
    (1.5,0.5);%
    \fill[helmholtzblue!20] (3,-0.5) node[anchor=north west]{} rectangle
    (4,0.5);%
    \begin{scope}[domain=-0.75:0]
      \draw[shift={(1.375,0.375)},color=myblue, samples=20]
      plot[smooth](\x,{.1*sin(2.666*pi*\x r)}); %
      \draw[shift={(1.375,0.125)},color=myblue, samples=40]
      plot[smooth](\x,{.1*sin(2.666*2*pi*\x r)}); %
      \draw[shift={(1.375,-0.125)},color=myblue, samples=80]
      plot[smooth](\x,{.1*sin(2.666*4*pi*\x r)}); %
      \draw[shift={(1.375,-0.375)},color=myblue, samples=160]
      plot[smooth](\x,{.1*sin(2.666*8*pi*\x r)}); %
      \draw[shift={(3.875,0.25)},color=myblue, samples=20]
      plot[smooth](\x,{.1*sin(2.666*pi*\x r)}); %
      \draw[shift={(3.875,-0.25)},color=myblue, samples=80]
      plot[smooth](\x,{.1*sin(2.666*2*pi*\x r)}); %
    \end{scope}
    \draw[->,thick] (4.05,0) -- ++ (.4,0) node[above, midway]{$\cg^{-1}$};%
    \begin{scope}[line width=10pt,helmholtzblue!30,]
      \draw[-{Stealth[scale=0.5]}] (1.525,0) -- ++ (0.45,0);%
      \draw[-{Stealth[scale=0.5]}] (2.55,0) -- ++ (0.45,0);%
    \end{scope}
    \begin{scope}[xshift=2cm,xscale=0.5,yscale=0.3,thick]
      \node[anchor=north] at (0.5,2) {$\mat P$};
      \foreach \name / \y in {-1,0,1}%
      \path node[hidneuron] (H-\name) at (0,\y cm) {};%
      \foreach \name / \y in {-1,0,1}%
      \path node[hidneuron] (H2-\name) at (1,\y cm) {};%
      \foreach \source in {-1,0,1}%
      \foreach \dest in {-1,0,1}%
      \path (H-\source) edge (H2-\dest);%
    \end{scope}
  \end{tikzpicture}
  \caption{Sketch of the FNO. First, the input is lifted to a higher dimensional
    space by $\symbfcal{R}$, which is defined in~\eqref{eq:fno-R}. Then apply
    the Fourier layers~\eqref{eq:fno-layer} and finally project onto the target
    by $\symbfcal{Q}$, as defined in~\eqref{eq:fno-Q}. In the Fourier layer
    defined by~\eqref{eq:fno-layer}, a Fourier Transform $\cg$ is applied and a
    linear transform filters out the higher modes. Then, the inverse Fourier
    transform $\cg^{-1}$ is applied. A linear transform $\mat{W}$ acts similar
    to a skip connection.}\label{fig:fno}
\end{figure}

FNO is a recently introduced type of ANN that proposes a novel method for
combining neural networks with Fourier analysis, mainly to solve differential
equations~\autocite{liFourierNeuralOperator2020}. Within the framework of Neural
Operators, a universal approximation theorem and error bounds have been
developed for the FNO in \autocite{kovachkiUniversalApproximationError2021}. The
key innovation of the FNO is a new type of layer, the Fourier layer
(cf.~Fig.~\ref{fig:fno} and~\eqref{eq:fno-layer}). In the Fourier layer the Fourier series
is used to efficiently compute the convolution of the input function with a set
of integration kernels, represented in the frequency domain.

Here, we briefly introduce FNOs for complex-valued functions $\vct v\in
L^1(\mathbb{T}^d)$ on the unit torus
\(\mathbb{T}^d\), in order to restrict ourselves to
1-periodic functions. For details we refer to~\autocite[Section
3.1]{grafakosClassicalFourierAnalysis2014}. The Fourier transform of a function
\(\vct v\colon \mathbb{T}^d \to \C^n\) is denoted by
\(\cg\colon L^2(\mathbb{T}^d;\C^n) \to \ell^2(\Z^d;\C^n)\). Similarly,
\(\cg^{-1}\colon \ell^2(\Z^d;\C^n)\to L^2(\mathbb{T}^d;\C^n)\) denotes the
Fourier inversion. More precisely, for a function
\(v \in L^2(\mathbb{T}^d;\C)\) the Fourier transform is defined by (cf.~\autocite[Definition
  3.1.1]{grafakosClassicalFourierAnalysis2014})
\begin{align}
  \label{eq:fourier-t}(\cg v) (l) &= \int_{\mathbb{T}^d}v(x) \exp(-2\pi \iu \dd{l}{x})\drv x\,,
                \qquad l \in \Z^d\,. \\
  \intertext{For a
  function \(w \in \ell^2(\Z^d;\C)\) is given by Fourier inversion (cf.~\autocite[Proposition 3.2.5]{grafakosClassicalFourierAnalysis2014})
  by}
  \label{eq:fourier-it}(\cg^{-1} w) (x) &= \sum_{l \in \Z^d} w (l) \exp(2\pi \iu \dd{l}{x})\,, \qquad x \in \mathbb{T}^d\,.
\end{align}
For vector-valued functions, the formulas~\eqref{eq:fourier-t}
and~\eqref{eq:fourier-it} are applied componentwise. We note that, for an
integrable function $\vct u$ on $\R^n$ with Fourier transform
$\widehat{\vct u}$, the Fourier series and Fourier inversion can be seen as the
restriction of the classical Fourier transform to $\Z^n$. Together with the
Poisson summation formula~\autocite[Theorem
3.2.8]{grafakosClassicalFourierAnalysis2014}, we know that the Fourier expansion
equals the periodization of the function $\vct v$ on $\R^n$. This gives us a
perspective on the extension of the Fourier series, and therefore FNOs, to
non-periodic functions.
\begin{definition}{Fourier Neural Operator (FNO)}{fno}
  An FNO
  $\symbfcal{N}\colon L^2(\mathbb{T}^d,\R^{d_i})\to
  L^2(\mathbb{T}^d,\,\R^{d_o})$ is a mapping consisting of a concatenation of
  functions such that
  \begin{equation}
    \label{eq:fno-al}
    \symbfcal{N}(\vct v)=\symbfcal{Q}\circ \symbfcal{L}_{N}\circ \cdots \circ \symbfcal{L}_1\circ
    \symbfcal{R}(\vct v)\,,
  \end{equation}
  with a lifting operator $\symbfcal{R}$ and a projection operator $\symbfcal
  Q$, represented by matrices $\mat R\in \R^{n\times d_i}$ and $\mat
  Q\in\R^{d_{o}\times n}$, respectively
  \begin{multicols}{2}
    \noindent
    \begin{equation}\label{eq:fno-R}
      \begin{split}
        \symbfcal{R}\colon L^2(\mathbb{T}^d;\,\R^{d_i})&\to L^2(\mathbb{T}^d;\,\R^n)\,,\\
        \vct v&\mapsto \symbfcal{R} \vct v\,,\\
        (\symbfcal{R} v)(x)&=\mat R (v(x)),\;\mat R\in \R^{n\times d_i}\,.
      \end{split}
    \end{equation}
    \begin{equation}\label{eq:fno-Q}
      \begin{split}
        \symbfcal{Q}\colon L^2(\mathbb{T}^d;\,\R^n)&\to L^2(\mathbb{T}^d;\,\R^{d_o})\,,\\
        \vct v&\mapsto \symbfcal{Q} \vct v\,,\\
        (\symbfcal{Q} v)(x)&=\mat Q (v(x)),\;\mat Q\in \R^{d_o\times n}\,.
      \end{split}
    \end{equation}
  \end{multicols}
  A Fourier layer $\symbfcal{L}_k$ is given by
  \begin{equation}
    \label{eq:fno-layer}
    \symbfcal{L}_{k}(\vct v) =\vct{\sigma} \big(\mat{W}_k \vct v + \vct b_k + \underbrace{\cg^{-1} \paran{P_n \cg(\vct v)}}_{\displaystyle \symbfcal{K}_n\vct v}\big)\,,
  \end{equation}
  where $\mat{W}_k\in \R^{n\times n}$ is a weight matrix and $\vct b_k\in \R^n$
  a bias vector and $P_n\colon \Z^d\to\C^{n\times n}$,
  $l \mapsto P_n(l)\in \C^{n\times n}$ are the weights of the modes $l\in \Z^d$
  and $\vct \sigma \colon \R^n\to\R^n$ is an activation function, for instance the
  $\tanh$ function is applied.
\end{definition}
Concerning Definition~\ref{def:fno} we note the following.
\begin{itemize}\itemsep1pt \parskip0pt \parsep0pt
\item Let us consider \((\cg \vct v)(l) \in \C^{n}\). In order to ensure that
  \(\symbfcal{K}_n\vct v\) in~\eqref{eq:fno-layer} is real-valued for real-valued \(\vct v\) conjugate
  symmetry in the parametrization is enforced by
  \begin{equation}
    \label{eq:conj-symm}
    P_n(-l)_{j,\,k} = P_n^*(l)_{j,\,k}\,, \;
    j=1,\dots,\,m, \; k=1,\dots,\,n\, \quad \forall l \in Z_{l_{\text{max}}};
  \end{equation}
\item In~\eqref{eq:fno-R} and~\eqref{eq:fno-Q}, $\symbfcal{R}$ and
  $\symbfcal{Q}$ are both locally acting operators. They are represented by the
  matrices $\mat R$ and $\mat Q$, that are to be trained.
\item We restrict the domain of the FNO to $\mathbb{T}^d$, in order to consider only 1-periodic
  functions. The Poisson summation formula lets us lift this restriction.
  Similarly, in~\autocite[Lemma 41]{kovachkiUniversalApproximationError2021} the authors
  show that FNOs can be generalized to domains with
  Lipschitz boundary.
\item The activation function $\vct \sigma \colon \R^n\to\R^n$ with
  $\symbf\sigma \vct v= {(\sigma(v_{1}),\dots,\,\sigma(v_{n}))}^{\top}\in \R^n$
  is a componentwise applied scalar- and real-valued, non-polynomial function
  $\sigma\in C^{\infty}(\R)$, which is globally Lipschitz-continuous.
\end{itemize}
We sketch the FNO in Fig.~\ref{fig:fno}. The key feature of FNO architectures
are the convolution-based integral kernels $\symbfcal{K}_n$, that are non-local.
This enables learning operators with a global character, such as operators
arising in the simulation of PDEs. Another major factor in the efficiency is
that in the discrete case we are able to use the Fast Fourier Transform (FFT) to
compute \(\symbfcal{K}_n\vct v\) in~\eqref{eq:fno-layer},
if the computational mesh is uniform. This is sketched in the following.

\subsubsection*{The Discrete Setting}
\label{sec:org74ff5eb}
Let the \(D_{J}\subset \mathbb{T}^d\) be a set of \(J \in \mathbb{N}\) uniformly distributed points with resolution \(s_1 \times \cdots \times
s_d = J\) in the domain
 \(\mathbb{T}^d\), \(v \in \C^{J \times n}\) and \(\cg (v) \in \C^{J \times n}\).
The multiplication by the weight tensor \(P \in \C^{J \times m \times n}\) is defined
by the operation
\begin{equation}
\label{eq:fft_mult}
\bigl( P \cdot (\cg v) \bigr)_{k,\,l} = \sum_{j=1}^{n} P_{k,\,l,\,j}  (\cg v)_{k,\,j}\,, \qquad k=1,\dots,\,J, \quad l=1,\dots,\,m.
\end{equation}
The Fourier transform \(\cg\) can be replaced by the Fast Fourier Transform (FFT). For \(v \in \C^{J
\times n}\), \(k = (k_1, \ldots, k_{d}) \in \mathbb{Z}_{s_1} \times \cdots \times
\mathbb{Z}_{s_d}\), and \(x=(x_1, \ldots, x_{d}) \in \mathbb{T}^d\), the FFT \(\widehat{\cg}\) and its
inverse \(\widehat{\cg}^{-1}\) are defined as
\begin{align*}
    (\widehat{\cg} v)_l(k) = \sum_{x_1=0}^{s_1-1} \cdots \sum_{x_{d}=0}^{s_d-1} v_l(x_1, \ldots, x_{d}) \exp\paran{- 2\iu \pi \sum_{j=1}^{d} \frac{x_j k_j}{s_j} }\,,\;\text{for}\;l=1,\dots,\,n, \\
    (\widehat{\cg}^{-1} v)_l(x) = \sum_{k_1=0}^{s_1-1} \cdots \sum_{k_{d}=0}^{s_d-1} v_l(k_1, \ldots, k_{d}) \exp\paran{2\iu \pi \sum_{j=1}^{d} \frac{x_j k_j}{s_j} }\,,\;\text{for}\;l=1,\dots,\,n.
\end{align*}
The parameters $\mat W$, $\vct b$, $\mat P$ of the Fourier layers in Definition~\ref{def:fno} are learned in Fourier space, where they
can be expressed in terms of the Fourier coefficients of the input functions.
When the network is used to evaluate functions in physical space, it simply
amounts to projecting onto the basis functions
\(\exp\paran{2\pi i \langle x, k \rangle}\), which are well-defined for all
\(x \in \C^d\). This allows the network to evaluate functions at any desired
resolution, without being tied to a specific discretization scheme. The
implementation of the FNO using the FFT restricts the geometry and
discretization to uniform mesh discretizations of $\mathbb{T}^d$. In practice
FNOs can be extended to other domains by padding the input with zeros. The loss is
computed only on the original domain during training. The Fourier neural
operator extends the output smoothly to the padded domain, as discussed
in~\autocite{kovachkiNeuralOperatorLearning2023}.

\subsection{Recurrent neural networks with memory}
\label{sec:orgd48a483}
Recurrent neural networks (RNNs) are an extension to Feed-forward neural
networks that use an activation variable \(\vct a_n^k\in\mathbb{R}^{p_k}\) to propagate
information over discrete time steps, making them suitable for time series and
sequential data. An extension of this model uses network nodes with memory.
These neural networks are effective in modeling long-term dependencies and can
overcome the vanishing gradient problem that recursive neural networks face
\autocite{hochreiterGradientFlowRecurrent2001}. In our work, we employ Gated
Recurrent Units (GRUs)~\autocite{choLearningPhraseRepresentations2014}.
\begin{definition}{Gated Recurrent Unit (GRU)}{gru}
A Gated Recurrent Neural Network $\mathcal{N}\colon X^N\to Y^N$ maps a
sequence of elements of a finite dimensional inner product space $X$ to
a sequence of elements of a finite dimensional inner product space $Y$.
It consists of a concatenation multiple GRUs, i.e.
\begin{equation}
  \label{eq:grnn}
\mathcal{N} = \symbfcal{G_1}\circ\cdots\circ\symbfcal{G}_L
\end{equation}
A GRU $\symbfcal{G_k}\colon \R^{N}\times\R^{q}\to\R^{N\times q}\times\R^{q}$, $k\in\{1,\dots,\,L\}$ is defined by the
following equations:
\begin{subequations}
  \label{eq:gru}
\begin{align}
  \label{eq:gru-z}\vct z_n^k &= \vct \sigma_z \big( \mat W_k^{(z)} \, \vct h_n^{k-1} + \mat U_k^{(z)} \, \vct h_{n-1}^k + \vct b_k^{(z)} \big)\,,
  \\[3pt]
  \label{eq:gru-r}\vct r_n^k &= \vct \sigma_r \big( \mat W_k^{(r)} \, \vct h_n^{k-1} + \mat U_k^{(r)} \, \vct h_{n-1}^k + \vct b_k^{(r)} \big)\,,
  \\[3pt]
  \label{eq:gru-h}\vct h_n^k &= \vct z_n^k \odot \vct h_{n-1}^k + (1 - \vct z_n^k) \odot \vct \sigma_h \big( \mat W_k^{(h)} \, \vct h_n^{k-1} + \mat U_k^{(h)} (\vct r_n^k \odot \vct h_{n-1}^k) + \vct b_k^{(h)} \big)\,,
\end{align}
\end{subequations}
where $n\in \{1,\dots,\,N\}$, $\odot$ denotes the element-wise product and $\mat U_k^{( \cdot )}$ $\mat W_k^{( \cdot )}$
and $\vct b_k^{( \cdot )}$ are weight matrices and bias vectors determined by
training. By $\vct h_n^0$ we denote the
input to the GRU.
\end{definition}
In~\eqref{eq:gru}, the update gate vector \(\vct z_n^k \in \R^q\) defined
in~\eqref{eq:gru-z} determines the contribution of the previous hidden output
\(\vct h_{n-1}^k\) to the current output \(\vct h_n^k\) (cf.~\eqref{eq:gru-z}),
while the reset gate vector \(\vct r_n^k \in \R^q\) defined by
(cf.~\eqref{eq:gru-r}) controls the nonlinearity of the cell. Together, they
control the memory of a GRU cell, determining to what extent information from
the past is carried over to the present output.

\section{Optimal control with Neural Operators}
\label{sec:orga943742}
Optimal control problems (OCP) are important in several branches of science and
engineering. Finding efficient solutions to these problems remains a challenging
task. Neural operators can represent the dynamics of complex systems
efficiently. Their combination with OCPs has the potential to yield novel
solutions by replacing the oftentimes costly solution of the forward problem. In
this section, we investigate the use of neural operators for solving OCPs, with
a focus on Dirichlet boundary conditions as constraints. We apply this technique
to the problem of THz generation in a periodically poled crystal
(cf.~Section~\ref{sec:org8d1a271}) and propose a
novel approach to optimize the input pulse with the goal to maximize the
efficiency of optical to THz generation. This can be formulated as an optimal
boundary control problem, where we seek the Dirichlet boundary conditions that
yield the maximum optical to THz conversion.
\subsection{Optimal Dirichlet boundary control}
First, we state a general optimal Dirichlet boundary control problem, which
serves as the foundation for our proposed method.
\begin{definition}{Function spaces for optimal control problem}{fn-ocp}
  Define $S=I\times \Gamma_D$ and $Q=I\times \mathcal{D}$ and the set
  of admissible controls
  \begin{equation}
    \label{eq:adm-controls}
    U_{\text{ad}}=\brc{u\in L^2(S)\suchthat u_a\leq u\leq u_b\:\text{a.\,e.\ in}~S,\:u_{a},\:u_b\in L^2(S)}\,.
  \end{equation}
\end{definition}
Let \(\mathcal{J}\) be a Gateaux differentiable
functional. For the state \(\vct y\in W(I)\) and the control
\(u\in U_{\text{ad}}\) we consider the following optimization problem.
\begin{problem}{Optimal Dirichlet boundary control}{opcon}
  For \(\mathcal{J}\colon \mat W(I)\times L^2(S)\to \R\), \((\vct y,
  u)\mapsto\mathcal{J}(\vct y,\,u)\) solve
  \begin{subequations}
    \begin{align}
      \label{eq:goal}&\min_{u\in U_{\text{ad}}} &\mathcal{J}(\vct y,\,u) &= \mathcal{G}(\vct y) + \frac{\alpha}{2}\norm{u}^2\,,\;\alpha>0\,,\\
      &\text{subject to }& \vct y&=\symbfcal{S}(u,\,\vct f,\,\vct y_0)\,.
    \end{align}
  \end{subequations}
\end{problem}
The operator \(\symbfcal{S}\) is the abstract solution operator of the PDE
introduced for our application in~\eqref{eq:abstract-inv}, by which the optimization problem is
constrained. The control \(u\) enters through the Dirichlet boundary condition.
The functional \(\mathcal{G}\colon \mat W(I) \to \R\) is left to
be defined for the application.

We now derive an OCP similar to Problem~\ref{problem:opcon} for a
setting where Problem~\ref{problem:lor-ade} provides the initial and boundary
conditions and the PDE for the optimization problem. In practice, we define the
functional \(\mathcal{G}\) in~\eqref{eq:goal} such that radiation at
frequency \(f_{\Omega}\) is optimized. Since \(\max \mathcal{J}= - \min
(-\mathcal J)\) is satisfied, we restrict ourselves to the description of minimization
problems.
\begin{definition}{Cost function for optimizing generation of THz
    radiation}{cost-thz}
  Let $y_{c}\colon W_{\text{nl}}(I) \to L^{2}(I;\,\R)$ and $\psi \colon \R \to \R$ be
  given by
  \begin{equation}
    \label{eq:cost-fn-dc}
    y_{c}(\vct y)= \int_{B_{\varepsilon}(c)} \vct y(x,\, t)\drv x\,,
  \end{equation}
  \begin{equation}
    \label{eq:cost-fn-psi}
    \psi(\nu) = \symbf 1_{(f_{\Omega}-r,\,f_{\Omega}-r)}\exp\paran{\frac{r^2}{(\nu-r-f_{\Omega})(\nu+r+f_{\Omega})}}\,,\:r>0\,,
  \end{equation}
  where the ball $B_{\varepsilon}(c)$ around the control point $c\in \mathcal{D}$ is chosen such that
  $B_{\varepsilon}(c)\cap \Gamma_{D}=\emptyset$ is satisfied. We define the
  cost functional $\mathcal{G}$ as
  \begin{equation}
    \label{eq:cost-fn}
    \mathcal{G}_{\Omega}(\vct y)= \int_{f_{\Omega}-r}^{f_{\Omega}+r}\cg (y_{c}(\vct y,\,t))(\nu)^{2}
    \psi(\nu)\drv \nu\,.
  \end{equation}
\end{definition}
We note that \(y_{c}\in L^2(I;\,\R)\) and therefore its Fourier transform exists.
The parameter \(r\) is chosen such that \(\psi\) is sufficiently close to the
indicator function at \(f_{\Omega}\). In the discrete case, we specify this more
precisely. For the state \(\vct E\in W_{\text{nl}}(I)\) and the control \(g^{\vct e}\in
U_{\text{ad}}\) we study the following optimization problem.
\begin{problem}{Optimal Dirichlet boundary control for THz generation}{opconTHz}
  For the solution operator~\eqref{eq:abstract-inv} to
  Problem~\ref{problem:lor-stm}, solve the optimization problem
  \begin{subequations}
    \begin{align}
      &\max_{u\in U_{\text{ad}}} && &\pushleft{\mathcal{J}(\vct E,\,g^{\vct e}) = \mathcal{G}_{\Omega}(\vct E) + \frac{\alpha}{2}\norm{g^{\vct e}}^2\,,}\\
      \notag
      &\text{subject to }&& & \pushleft{(\vct u,\,\vct p,\,\vct a,\,\vct
                              e)^{\top}=\symbfcal{S}(g^{\vct e},\,\symbf{f}\,,\vct v_0)}\,.
    \end{align}
  \end{subequations}
\end{problem}
We provide realistic parameters for this problem in Section~\ref{sec:org4ff2e73}.
In the formulation of Problem~\ref{problem:opconTHz}, we can
replace the solution operator $\symbfcal{S}$, the data $\symbf{f}$, and $\vct v_0$
with their discrete counterparts as defined in
Definition~\ref{def:abstract-discrete} in a straightforward manner. We can
evaluate the cost function in the discrete setting using an FFT. While different
methods exist for the solution of OCPs similar
to~\ref{problem:opcon}, to the best of our knowledge the nonlinear wave equation
of Problem~\ref{problem:opconTHz} has not yet been investigated in theory or
practice. Within this work we concentrate on the algorithmic and practical
aspects of solving Problem~\ref{problem:opcon}. For an overview over optimal
control theory and solution methods we refer to
\autocite{manzoniOptimalControlPartial2021,hinzeOptimizationPDEConstraints2009a}
and references therein and, more specifically for hyperbolic problems, to
\autocite{gugatOptimalBoundaryControl2015}.
\subsection{Optimal Control for THz generation with Neural operators}
Even in one space dimension solving Problem~\ref{problem:opconTHz} by using the
variational space-time methods we presented so far is infeasible for scenarios
of practical interest due to the substantial computational burden imposed by the
solution of the forward problem;
cf.~\autocite{margenbergAccurateSimulationTHz2023}. In order to focus the
presentation on the essential ideas, we restrict ourselves to the
one-dimensional case for the remainder of this section.
In~\autocite{margenbergAccurateSimulationTHz2023} we observed that this is a
reasonable restriction to make from a practical point of view.

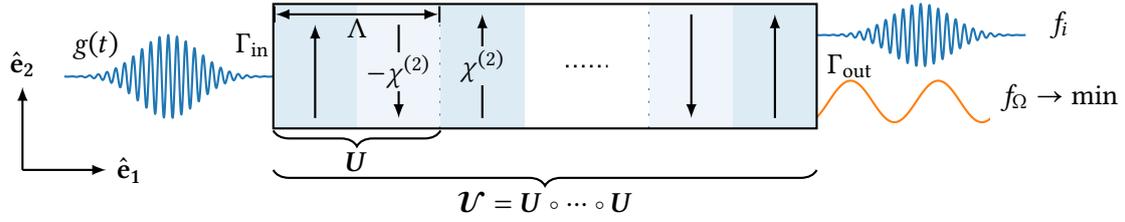
\begin{figure}
  \centering
  \begin{tikzpicture}[xscale=5.5,yscale=5.5,domain=-0.5:0,>=latex,thick]
    \draw[color=mpblue, samples=2000]
    plot[smooth](\x,{.125+.1*exp(-((\x+0.25)*(\x+0.25))*100)*cos(300*\x r)}); %
    \node (a) at (-0.425,0.2) {$g(t)$}; %
    \draw[->] (-0.6,-0.1) -- ++ (0.2,0.0) node[right] {$\vct{\hat e_1}$}; %
    \draw[->] (-0.6,-0.1) -- ++ (0.0,0.2) node[above] {$\vct{\hat e_2}$}; %
    \draw[color=mpblue, samples=2000, xshift=1.8cm] plot[smooth]
    (\x,{.225+.075*exp(-((\x+0.25)*(\x+0.25))*100)*cos(300*\x r)}); %
    \node (b) at (1.88,0.25) {$f_i$}; %
    \draw[color=mporange, samples=2000, xshift=1.8cm,yshift=-0.01cm]
    plot[smooth] (\x,{.075+.05*cos(30*\x r)}); %
    \node[fill=white] (c) at (1.88,0.075) {$f_{\Omega}\rightarrow\min$}; %
    \fill [color=mpblue!15] (0,0) rectangle (0.2,0.3); %
    \fill [color=mpblue!7] (0.2,0) rectangle (0.4,0.3); %
    \draw [color=mpblue,loosely dotted] (0.4,0) -- (0.4,0.3); %
    \fill [color=mpblue!15] (0.4,0) rectangle (0.6,0.3); %
    \draw [color=mpblue,loosely dotted] (0.9,0) -- (0.9,0.3); %
    \fill [color=mpblue!7] (0.9,0) rectangle (1.1,0.3); %
    \fill [color=mpblue!15] (1.1,0) rectangle (1.3,0.3); %
    \draw[decorate,decoration={brace,amplitude=6pt,mirror}] (0,-0.005) --
    node[fill=white,pos=0.5,inner sep=0pt,yshift=-12pt]{$\symbf{U}$} ++
    (0.4,0.); %
    \draw[decorate, decoration={brace,amplitude=6pt,mirror}] (0,-0.1) --
    node[fill=white,pos=0.5,inner
    sep=0pt,yshift=-12pt]{$\symbfcal{U}=\symbf{U}\circ\cdots \circ\symbf{U}$} ++
    (1.3,0.); %
    \draw[|<->|, decorate, decoration={amplitude=3pt}] (0,0.275) --
    node[pos=0.5,inner sep=0pt,yshift=-5pt]{$\Lambda$} ++ (0.4,0.); %
    \path[draw=black] (0,0) -- ++(1.3,0) -- node[pos=0.5,
    right]{$\Gamma_{\text{out}}$} ++(0,0.3) -- ++(-1.3,0) -- node[pos=0.33,
    left,inner sep=1pt]{$\Gamma_{\text{in}}$}cycle; %
    \draw[->] (.1,.025)-- ++ (0,0.225); %
    \draw[<-] (.3,.025) -- node[fill=mpblue!7,inner
    sep=1pt,pos=0.5]{$-\chi^{(2)}$}++ (0,0.225); %
    \draw[->] (.5,.025) -- node[fill=mpblue!15,inner
    sep=1pt,pos=0.5]{$\chi^{(2)}$} ++ (0,0.25); %
    \draw[<-] (1.0,.025) -- ++ (0,0.25); %
    \draw[->] (1.2,.025) -- ++ (0,0.25); %
    \draw[dotted,line width=1pt] (0.7,.15) --++(.1,0); %
  \end{tikzpicture}
  \caption{Sketch of using ANNs for the simulation of wave propagation in
    periodically poled nonlinear materials. The solution operator $\symbf U$,
    solves the wave equation~\eqref{eqwave:lor-nade2a} over one period.
    $\symbf U$ is learned from simulation data generated by space-time FEM. The
    training data for the $\symbf U$ is collected at the interfaces of two
    periods, which are separated by a distance of
    $\Lambda$.}\label{fig:periodsop}
\end{figure}

We propose an algorithm
that relies on ANNs to accelerate the solution of the PDE, allowing for a more
efficient optimization of the control parameters. In Fig.~\ref{fig:periodsop} we
sketch its key idea. The goal is to train neural operators \(\symbfcal{U}\)
based on accurate numerical simulations which generalize well to
\(U_{\text{ad}}\). Then they are used as the forward solver in the optimal
Dirichlet boundary control problem.
The first cornerstone of the method is to consider only controls that are
feasible in practice. By this we can implement a differentiable sampler of
Dirichlet boundary data in a deep learning library of our choice and concatenate
it with the solution operator.
\begin{definition}{Admissible controls for
    Problem~\ref{problem:opconTHz}}{diri-ocp}
  The Dirichlet data, i.\,e.~the control in Problem~\ref{problem:opconTHz}, is of
  the form
  \begin{subequations}
    \label{eq:sampling}
    \begin{align}
      \label{eq:samples}
      g(t)&=\exp\paran{-{\big(2\log
            2{\big(\tfrac{t}{\tau}\big)}^{2}\big)}^{p}}
            \sum_{i=1}^n a_i\cos\paran{\varphi_i + 2 \pi\paran{\frac{1}{2} \zeta_i t^2 + f_it}}\,.\\
      \intertext{For some fixed $n$ the set of
      parameters $\Xi$ and a sampler $\symbfcal{P}$ which maps these parameters
      to the pulses are given by}
      \notag\Xi&=\big\{(\tau,\,p,\,a_0,\,\varphi_0,\,\zeta_0,\,f_0,\dots,\,a_n,\,\varphi_n,\,\zeta_n,\,f_n)\in \R_{+}^{2+4n}\suchthat
                 \tau\leq \tau_{\max},\,p\leq p_{\max},\\
      \label{eq:xin}
          &\qquad a_{i}\leq
            a_{\max},\,\varphi_{i}\leq \varphi_{\max},\,\zeta_{i} \leq
            \zeta_{\max},\,f_{i}\leq f_{\max}
            \;\forall\,i=1,\dots,\,n\big\}\,,\\[.25ex]
      \label{eq:pulse-sampler}
      \symbfcal{P}\colon \Xi
          &\to L^2(\Gamma_{\text{in}}\times I),\;
            (\tau,\,p,\,a_0,\,\varphi_0,\,\zeta_0,\,f_0,\dots,\,a_n,\,\varphi_n,\,\zeta_n,\,f_n)\mapsto
            g(t)\,.
    \end{align}
  \end{subequations}
\end{definition}
In~\eqref{eq:sampling}, the parameter \(\tau\) is the full width half maximum, \(p\) the order
of the supergaussian, \(a_i\) the amplitude, \(\phi_i\) the phaseshift,
\(\zeta_i\) the quadratic chirprate and \(f_i\) the center frequency. The upper
bounds in \eqref{eq:xin} are given through the limitations of the experimental
setup. Note that the image of \(\symbfcal{P}[\Xi]\subset L^2(S)\) is the set of
admissible controls \(U_{\text{ad}}\) in Problem~\ref{problem:opconTHz}.

The second idea of the method is the differentiability of a program written with
an established ANN library: For most operations on the datastructures of these
libraries a method for calculating its gradient is already implemented. The last
idea in our algorithm builds on the periodicity of the material parameters (cf. Fig.~\ref{fig:periodsop}). The
example in this work is based on periodically poled crystals, where
\(\chi^{(2)}\) govern the nonlinear processes. The periodicity of \(\chi^{(2)}\)
can be used to learn a solution operator \(\symbf{U}\) to the forward problem in
only one period of the \(\chi^{(2)}\) parameter. We formalize this concept in the discrete setting: Consider the
discrete solution operator given in~\eqref{eq:abstract-inv-disc}. The goal is to
approximate \(\symbfcal{S}_{\tau,\,h}(g(t))\) by
\begin{equation}
\label{eq:rec-solop}
\symbf{U}\circ\cdots\circ\symbf{U}(g(t))\eqcolon\symbfcal{U}\approx\symbfcal{S}_{\tau,\,h}(g(t))\,.
\end{equation}
On the other hand, for the efficient solution of Problem~\ref{problem:opconTHz} we
don't need the full space-time solution. We only need the solution at some
collocation points
\begin{equation}
\label{eq:collocation-points}
J_{\mathcal{D}}
=\bigcup_{i=1}^m \brc{x_{\mathcal{D},\,i}}\,,\quad
x_{\mathcal{D},\,i}\coloneq\Lambda i\in \mathcal{D}\,,
\end{equation}
where \(m\) is the number of periods in the crystal.
For the
generation of training data we evaluate and save the solution \(\vct E(x,\,t)\)
at the points in \(J_{\mathcal{D}}\). Motivated by the fact that \(J_{\mathcal{D},\,i+1}\)
is the set \(J_{\mathcal{D},\,i}\), shifted in positive \(x_1\)-direction, we construct an operator
\(\symbf{U}\) which maps the time trajectory of the electric field \(\vct E\) at
period \(i\) to the time trajectory of \(\vct E\) at period \(i+1\). In order to
construct a suitable solution operator we define the space
\begin{equation}
\label{eq:no-space}
V_{\tau}(I)=\brc{\vct w\in L^2(I;\,\R)\suchthat \vct w(t)\restrict{I_n}\in \mathbb{P}_3(I_n,\,\R)}\,.
\end{equation}
\begin{definition}{Neural Operator for Problem~\ref{problem:opconTHz}}{no}
  Let $i\in\brc{1,\dots,\,m}$, $x_{a}\in J_{\mathcal{D},\,i}$,
  $x_{b}=x_a+\vct e_1\Lambda\in J_{\mathcal{D},\,i+1}$ and
  $p(t;\,x_{a})\in V_{\tau}(I)$. The neural operator
  $\symbf{U}=\mat I\circ\symbfcal{N}\circ \mat{T}$ is constructed
  such that
  \begin{equation}
    \label{eq:nop}
    \begin{split}
      \symbf{U}\colon V_{\tau}(I)&\to V_{\tau}(I)\,,\\
      \vct p(t;\,x_{a})&\mapsto \vct{\hat{p}} (t;\,x_b)\,,
    \end{split}
  \end{equation}
  where $\symbfcal{N}$ is one of the networks introduced in Section~\ref{sec:ANN}, $\mat T$ evaluates $\vct p\in V_{\tau}(I)$ at the endpoints of the subinterval and $\mat I$ is the Hermite-type interpolator~\eqref{eq:hermite_interpolation}
  applied on each subinterval $I_n$,
  \begin{multicols}{2}
    \noindent
    \begin{equation}
      \label{eq:no-T}
      \begin{split}
        \mat T\colon V_{\tau}&\to \R^{N+1}\,,\\
        \vct p(t;\,x_a)&\mapsto (\vct p(t_0;\,x_a),\dots,\,\vct p(t_N;\,x_a))\,,
      \end{split}
    \end{equation}
    \begin{equation}
      \label{eq:no-That}
      \begin{split}
        \mat{I}\colon \R^{N+1}&\to V_{\tau}\,,\\
        \vct u&\mapsto  \vct{\hat p}(t;\,x_b)\,.
      \end{split}
    \end{equation}
  \end{multicols}
\end{definition}
For the evaluation of \(\mat I\) we use automatic differentiation in order to
obtain \(\partial_t\symbfcal{N}\). Further, we note that the coefficients of the
polynomials and the values of \(\vct{\hat p}\) at the time endpoints of \(I_n\) coincide,
which makes the evaluation computationally cheap. With these preparations,
the computation of \(\vct{p}(t;x_{a})=\symbf{U}(g(t))\), \(x_{a}\in
J_{\mathcal{D},\,0}=\Gamma_{\text{in}}\) is well-defined and
\(\vct{p}(t;x_{a}+\Lambda)\) is the
time trajectory of the electric field \(\vct E\) at period \(1\).
We can iterate this to obtain the time trajectory at period \(i\),
\begin{equation}
\label{}
\vct{p}_{i-1}(t;\,x_{a}+\Lambda (i-1))=\symbf{U}^{i-1}(g(t))=\vct{p}_{i}(t;\,x_{a}+\Lambda i)\,.
\end{equation}
With the solution operator defined, we can formulate the algorithm for the
solution of the optimal Dirichlet boundary control problem. Through the
differentiability of \(\symbfcal{P}\circ \symbfcal{U}\) we can calculate the
gradients of the parameters \(\xi\in\Xi\) with respect to the cost function in
Problem~\ref{problem:opconTHz}. Then we can use the well-known gradient descent
algorithm or Newton's method for the solution of Problem~\ref{problem:opconTHz}.
In Algorithm~\ref{alg:opcondl} we describe the steps using a simple gradient
descent method, which can also be tracked in Fig.~\ref{fig:algsketch}. The
extension to Newton's method is straightforward. In
Appendix~\ref{sec:orgb4fe698} we give an abstract formulation of how a solution
operator for a full space-time approximation \(\symbfcal{U}\) can be obtained.
\begin{algo}{Optimal control based on deep learning}{opcondl}
  \begin{algorithm}[H]
    \NoCaptionOfAlgo %
    \LinesNumbered %
    \SetKwInOut{Input}{Input} %
    \SetKwProg{Function}{function}{}{end} %
    \SetKwRepeat{Do}{do}{while} %
    \Function{optimize-pulse-parameters($\symbfcal{P}$, $\symbfcal{U}$) :
        $\xi_{\text{opt}}$}{ %
      \label{alg:init}
      Initialize pulse parameters $\xi_0 \in \Xi$, Set $n=0\,$\;
      \Do{$\abs{\delta\xi}>\varepsilon\,$}{ %
        $g(t)=\symbfcal{P}(\xi)$ \tcp*{Sample pulse}%
        $\vct E(x,\,t)=\symbfcal{U}(g)$ \tcp*{Evaluate solution operator}%
        $\displaystyle\delta\xi=\alpha\frac{\partial}{\partial \xi_n}\mathcal{G}_{\Omega}(\vct E(x,\,t))$ \tcp*{Calculate update for gradient descent}%
        $\xi_{n+1}=\xi_n + \delta\xi$ \tcp*{Update pulse parameters}%
        $n= n+1\,$\; %
      } %
      Return $\xi_{\text{opt}}=\xi_n\,$\;
    } %
  \end{algorithm}
  \vspace*{-8pt}
\end{algo}
\FloatBarrier

\section{Numerical Experiments}
\label{sec:org4ff2e73}
We present numerical studies of the proposed neural operators for solving optimal
control problems. First, we investigate and validate their ability to efficiently represent
the dynamics of a simple test case. Then we extend the test to our method of
optimal control via neural operators by adding a set of constraints and solving
the resulting optimal control probem. Finally, we demonstrate the feasibility of
the proposed approach by applying our methods to the physical problem of THz generation
and compare the results to experimental data. In this section we only consider
settings in one space dimension, since otherwise numerical simulations are too time-consuming. In
\autocite{margenbergAccurateSimulationTHz2023} we also restricted ourselves to 1D
without notable limitations.

\subsection{Implementation aspects}
\label{sec:org0dba394}
We implemented our numerical simulations using
\texttt{deal.II}~\autocite{arndtDealIILibrary2021}, a finite element toolbox that offers
efficient and scalable parallelization with MPI. To solve the nonlinear systems
of equations, we employ a Newton-Krylov method. For the linear systems of
equations that arise for each Newton iteration, we use the generalized minimal
residual method (GMRES) with the algebraic multigrid solver MueLu~\autocite{MueLu}.
MueLu serves as a preconditioner with a single sweep for every GMRES
iteration. We implemented the ANNs and the optimal control method proposed here with the
\texttt{C++} interface of
\texttt{PyTorch}~\autocite{paszkePyTorchImperativeStyle2019}, \texttt{libtorch}.
PyTorch also supports parallelization with MPI, which is used throughout this
work.
\subsection{Domain truncation}
\label{sec:org6202cdc}
In numerical simuations, wave propagation and other physical processes have to
be truncated to bounded regions. To this end, we extend \(\mathcal{D}\)
by a Perfectly Matched Layer (PML) on the right-hand side
\(\mathcal{D}_{F}=\mathcal{D}\cup\mathcal{D}_{\text{PML}}\).
We only consider the 1D case where \(\mathcal{D}=[0,\,L]\subset \R\)
is a bounded and closed interval. The PML can be written as
\(\mathcal{D}_{\text{PML}}=(L,\,L_{\mathcal{D}_{F}}]\) with
\(L_{\text{PML}}\coloneq L_{\mathcal{D}_F}-L\).
Inside the PML-region we have the problem
\begin{subequations}\label{eqpml:lor-ade-cfs-pml}
  \begin{align}
    \partial_{tt}\vct P + \Gamma_{0}\partial_{t}\vct P
    + \nu_{t}^2\vct P - \kappa_{x}\varepsilon_\Delta\nu_{t}^2 \vct E&=0&&\text{on}\quad \mathcal{D}_{\text{PML}} \times I\,,\\[.3ex]
    \partial_{t}\vct R+\alpha_{x}\vct R -\varepsilon_{\omega}\sigma_{x}\vct E&=0&&\text{on}\quad \mathcal{D}_{\text{PML}} \times I\,,\\[.3ex]
    \partial_{t} \vct Q + \tilde{\alpha}_{x} \vct Q - \tilde{\sigma}_{x} \partial_{x}
    \vct E &= 0&&\text{on}\quad \mathcal{D}_{\text{PML}} \times I\,,\\[.3ex]
    -\symbf \nabla\cdot \kappa_{x}^{-1}\symbf \nabla \vct E + \partial_{x} \vct Q +
    \kappa_{x}\varepsilon_{\omega} \partial_{tt} \vct E +
    \kappa_{x}(\varepsilon_{\Omega}- \varepsilon_{\omega})\nu_{t}^{2}\vct E\quad&\notag\\
    - \nu_{t}^{2}P -\Gamma_{0}\partial_{t}\vct P +\partial_{t}
    \paran{\varepsilon_{\omega}\sigma_{x}\vct E - \alpha_{x}\vct R}&=0&&\text{on}\quad \mathcal{D}_{\text{PML}} \times I\,,\\[.3ex]
    \vct E(0) = 0 \quad \partial_t \vct E(0)  & = 0 &&\text{on}\quad\; \mathcal{D}_{\text{PML}}\,,\\[.3ex]
    \vct E  &= 0 &&\text{on}\quad\Gamma_D\cap\bar{\mathcal{D}_{\text{PML}}} \times I\,.
  \end{align}
\end{subequations}
A more in-depth presentation with further discussion and references for PML can be found
in~\autocite{margenbergAccurateSimulationTHz2023}.

\subsection{Numerical convergence test of the space-time finite element method}
\label{sec:org02cbd0d}
Here we verify the numerical methods we developed for the forward problem. To
this end we prescribe a function as the solution to the equations in Problem
\ref{problem:lor-ade}. We use the residual of this function as a source
term, which in turn makes the prescribed function the solution. We use the
Galerkin--collocation method proposed in Problem~\ref{problem:lor-stm} for
the time discretization and the finite element space \(\symbfcal V_{h}\) defined in
\eqref{eq:fespace} for the spatial discretization. Consequently, we expect
fourth-order convergence.
\begin{table}[htbp]
\caption{\label{tab:org6f7a060}Calculated errors for the electric field \(\vct E\) and the auxiliary variable \(A\) for \(\text{GCC}^1(3)\).}
\centering
\small
\begin{tabular}{llrlrlrlr}
\toprule
\(k\) & \(L^\infty-L^2 (\vct E)\) & EOC & \(L^\infty-L^2 (A)\) & EOC & \(L^2-L^2 (\vct E)\) & EOC & \(L^2-L^2 (A)\) & EOC\\[0pt]
\midrule
\(\tau_0\) & \num{7.1334e-01} & - & \num{1.3573e+04} & - & \num{2.2943e-02} & - & \num{4.0377e+02} & -\\[0pt]
\(\tau_0\times 2^{-1}\) & \num{3.9379e-02} & 4.18 & \num{8.9722e+02} & 3.92 & \num{1.3091e-03} & 4.13 & \num{3.0466e+01} & 3.73\\[0pt]
\(\tau_0\times 2^{-2}\) & \num{2.5447e-03} & 3.95 & \num{6.0263e+01} & 3.90 & \num{8.3673e-05} & 3.97 & \num{2.0390e+00} & 3.90\\[0pt]
\(\tau_0\times 2^{-3}\) & \num{1.6050e-04} & 3.99 & \num{3.8358e+00} & 3.97 & \num{5.2351e-06} & 4.00 & \num{1.2894e-01} & 3.98\\[0pt]
\(\tau_0\times 2^{-4}\) & \num{1.0054e-05} & 4.00 & \num{2.4085e-01} & 3.99 & \num{3.2623e-07} & 4.00 & \num{8.0604e-03} & 4.00\\[0pt]
\(\tau_0\times 2^{-5}\) & \num{6.2863e-07} & 4.00 & \num{1.5070e-02} & 4.00 & \num{2.0340e-08} & 4.00 & \num{5.0310e-04} & 4.00\\[0pt]
\(\tau_0\times 2^{-6}\) & \num{3.9293e-08} & 4.00 & \num{9.4217e-04} & 4.00 & \num{1.2695e-09} & 4.00 & \num{3.1415e-05} & 4.00\\[0pt]
\bottomrule
\end{tabular}
\caption{\label{tab:org90460fc}Calculated errors for the auxiliary variables \(P\) and \(U\) for \(\text{GCC}^1(3)\).}
\centering
\small
\begin{tabular}{llrlrlrlr}
\toprule
\(k\) & \(L^\infty-L^2 (U)\) & EOC & \(L^\infty-L^2 (P)\) & EOC & \(L^2-L^2 (U)\) & EOC & \(L^2-L^2 (P)\) & EOC\\[0pt]
\midrule
\(\tau_0\) & \num{1.2918e+00} & - & \num{1.6725e-03} & - & \num{4.0207e-02} & - & \num{5.2300e-05} & -\\[0pt]
\(\tau_0\times 2^{-1}\) & \num{6.1522e-02} & 4.39 & \num{7.6061e-05} & 4.46 & \num{1.9033e-03} & 4.40 & \num{2.3583e-06} & 4.47\\[0pt]
\(\tau_0\times 2^{-2}\) & \num{3.7217e-03} & 4.05 & \num{4.4730e-06} & 4.09 & \num{1.1401e-04} & 4.06 & \num{1.3805e-07} & 4.09\\[0pt]
\(\tau_0\times 2^{-3}\) & \num{2.3102e-04} & 4.01 & \num{2.7585e-07} & 4.02 & \num{7.0556e-06} & 4.01 & \num{8.4946e-09} & 4.02\\[0pt]
\(\tau_0\times 2^{-4}\) & \num{1.4416e-05} & 4.00 & \num{1.7193e-08} & 4.00 & \num{4.3980e-07} & 4.00 & \num{5.2888e-10} & 4.01\\[0pt]
\(\tau_0\times 2^{-5}\) & \num{9.0060e-07} & 4.00 & \num{1.0741e-09} & 4.00 & \num{2.7465e-08} & 4.00 & \num{3.3023e-11} & 4.00\\[0pt]
\(\tau_0\times 2^{-6}\) & \num{5.6281e-08} & 4.00 & \num{6.7133e-11} & 4.00 & \num{1.7160e-09} & 4.00 & \num{2.0634e-12} & 4.00\\[0pt]
\bottomrule
\end{tabular}
\end{table}
We choose a 1D test case in the domain \(\mathcal{D}=[0,\,\num{0.001955}]\) over
the time interval \(I=[0,\,\num{1e-13}]\). As the electric field we choose
\begin{equation}
  \label{eq:testefield1}
  \vct E(x,\,t)=\sin \left(2 \pi \omega_{2} \left(x-n_{2} t\right)\right)+
  \sin \left(2 \pi \omega_{1} \left(x-n_{1} t\right)\right)\,.
\end{equation}
To compute the error in the physical domain and exclude error contributions from
within the PML region, we introduce a weighting function \(l\colon \mathcal{D}
\to \R\) that is equal to one in the physical domain and zero in the PML region:
\begin{equation*}
l(x)=
  \begin{cases} 0,\;x\in\mathcal{D}_{\text{PML}}\,,\\
    1,\;x\in\mathcal{D}\,.
  \end{cases}
\end{equation*}
Furthermore we multiply \(l\) by the source term to restrict it to the physical
domain. Thereby the solution inside \(\mathcal{D}\) is given
by~\eqref{eq:testefield1}. Then it propagates into \(\mathcal{D}_{\text{PML}}\)
where it is attenuated to the point of vanishing.
We study the errors \(e_{\vct Z}=\vct Z(x,\,t)-\vct Z_{\tau,\,h}(x,\,t)\) for \(\vct Z\in\{\vct E,\,\vct A,\,\vct P,\,\vct U\}\) in the norms
\begin{equation}
  \label{eq:norms}
  \norm{e_{\vct Z}}_{L^{\infty}(L^2)}=\max_{t\in I}\paran{ \int_{\mathcal{D}} \lvert e_{\vct Z}\rvert^{2}\drv x}^{\frac{1}{2}}\quad\text{and}\quad
  \norm{e_{\vct Z}}_{L^2(L^2)}=\paran{\int_{I}\int_{\mathcal{D}} \lvert
    e_{\vct Z}\rvert^{2}\drv x\drv t}^{\frac{1}{2}}\,.
\end{equation}
We abbreviate the error quantities \(\norm{e_{\vct Z}}_{L^{\infty}(L^2)}\) and
\(\norm{e_{\vct Z}}_{L^{2}(L^2)}\) by \(L^{\infty}\text{-}L^{2}(\vct Z)\) and
\(L^{2}\text{-}L^{2}(\vct Z)\) for
\(z\in \brc{\vct E,\,\vct A,\,\vct P,\,\vct U}\). The errors are calculated by
simultaneous refinement in space and time. In Table~\ref{tab:org6f7a060} we
observe the fourth order convergence in the variables \(\vct E\) and \(\vct A\).
For the auxiliary variables \(\vct U\) and \(\vct P\) we observe the same
convergence rates in Table~\ref{tab:org90460fc}, which highlights the advantage
of modelling these auxiliary variables with differential equations; cf.
\autocite{margenbergAccurateSimulationTHz2023}.

\subsection{Test of the solution operator and optimal control methodology}
\label{sec:org9eb2948}
In order to evaluate our algorithm, we construct an artificial test case similar
to a numerical convergence test studied before. This test aims to provide empirical evidence
of the algorithm's capability to solve complex problems. By this
rigorous evaluation, we hope to gain insights into its strengths and weaknesses
and identify optimal parametrizations that may be employed in the practical
test case.

\subsubsection{Training and testing of the solution operator}
\label{sec:org529e4c6}
In a first step we construct a test for the solution operator, where we consider
plane waves in vacuum. We generate the training data from plane waves at
frequencies
\(f\in\{\SI{291.56}{\tera\hertz},\allowbreak\,\SI{290.56}{\tera\hertz},\allowbreak\,2\cdot\SI{291.56}{\tera\hertz},\allowbreak\,2\cdot\SI{290.56}{\tera\hertz}\}\).
These frequencies are in the range of what we encounter in practice for the THz
generation. For the test we choose 2 plane waves with frequencies \(f_{1}\), \(f_2\) drawn from a continuous uniform distribution with support
\((\SI{290.56}{\tera\hertz},\,2\cdot\SI{291.56}{\tera\hertz})\). Then we add two
more plane waves by chosing the frequencies \(f_{3}=f_{1}+\SI{1}{\tera\hertz}\)
and \(f_{4}=f_{2}+\SI{1}{\tera\hertz}\). We choose a 1D test case in the domain
\(\mathcal{D}=[0,\,\num{3.e-14}c_{0}]\) over the time interval
\(I=[0,\,\num{1e-14}]\). The spatial domain has 3 periods with
\(\Lambda=\num{1.e-14}c_{0}\). Therefore, the solution operator is always applied
3 times onto itself. This is done throughout this subsection.

\begin{figure}[htbp]
\centering
\begin{minipage}{0.45\textwidth}
\begin{center}
\begin{tabular}{crrrr}
\toprule
\(l\) $\backslash$ \(w\) & 4 & 8 & 16 & 32\\[0pt]
\midrule
1 & 1317 & 2665 & 6993 & 22177\\[0pt]
2 & 1593 & 3761 & 11361 & 39617\\[0pt]
4 & 2145 & 5953 & 20097 & 74497\\[0pt]
8 & 3249 & 10337 & 37569 & 144257\\[0pt]
\bottomrule
\end{tabular}
\end{center}
\end{minipage}
\begin{minipage}{0.45\textwidth}
\begin{center}
\begin{tabular}{crrrr}
\toprule
\(l\) $\backslash$ \(w\) & 4 & 8 & 16 & 32\\[0pt]
\midrule
1 & 88 & 272 & 928 & 3392\\[0pt]
2 & 208 & 704 & 2560 & 9728\\[0pt]
4 & 448 & 1568 & 5824 & 22400\\[0pt]
8 & 928 & 3296 & 12352 & 47744\\[0pt]
\bottomrule
\end{tabular}
\end{center}
\end{minipage}
\captionof{table}{\label{orgf123153}The number of parameters of the FNO (left)
  and GRU (right) for different parametrizations: The number of layers \(l\) and
  the width of the layers \(w\).}
\vspace*{10pt}
  \includegraphics[width=\linewidth]{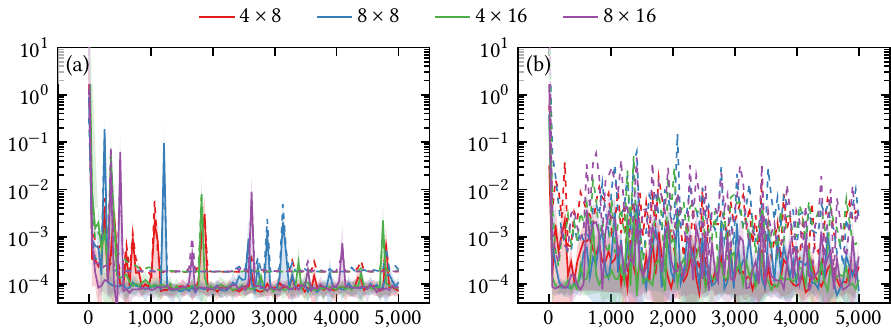}
  \captionof{figure}{The average training (solid) and validation (dashed) loss
    ($l_{C^{1}}^{2}$ Eq.~\eqref{eq:lc1}) for FNO (a) and GRU (b) architectures,
    smoothened with gnuplots smooth acsplines method. In order to avoid clutter
    and focus on the relevant results, we exclude larger network
    parametrizations that do not provide any additional accuracy benefits from
    the figure. The notation for the legend entries is $l\times w$, where $l$ is
    the number of layers and $w$ the layer width.}\label{fig:so-losses}
  \vspace*{10pt}
  \includegraphics[width=\linewidth]{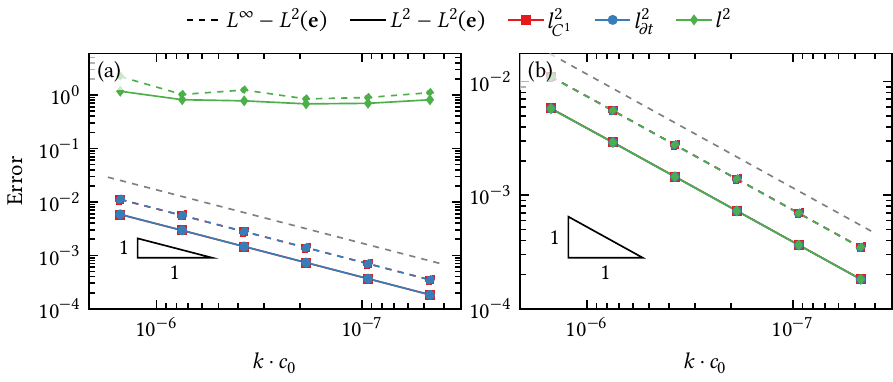}
  \captionof{figure}{Convergence of the $L^{\infty}(L^2)$ (dashed) and $L^{2}(L^2)$
    (solid) loss for FNO (a), GRU (b). The dashed grey line corresponds to
    linear convergence $O(n)$.}\label{fig:so-convergence}
\end{figure}

We test FNOs and GRUs. Each of these models are trained and evaluated with
varying numbers of layers and layer widths to evaluate their performances. In
Table~\ref{orgf123153} we collect all configurations used here. For some of
them, Fig.~\ref{fig:so-losses} shows the loss curves. The legend entries are
named according to the rows and columns in Table~\ref{orgf123153}. It is evident
that although both models achieve the same level of accuracy, the GRU exhibits
significant instability and oscillation in its loss. We attempted to address
this issue by using an annealing learning rate during training and conducted
extensive tuning of the hyperparameters, but the instability persisted. We note
that one batch of training data used for Fig.~\ref{fig:so-losses} already
contains \(\num{10000}\) timesteps, so the issue could be related to the long
term stability of the GRU. However, the prediction of 10000 timesteps is low
compared to our practical example in the following section. The FNO on the other
hand converges fast, especially for networks with 8 layers compared to the ones
with 4 layers (cf.~\ref{fig:so-losses}\,(a)). Due to the simplicity of the
problem setting, the best models exhibit a similar loss across all
architectures, eventhough the number of trainable parameters varies by multiple
orders; cf.~Table~\ref{orgf123153}. Although the GRUs are smaller in these
scenarios, the average training time is eight times longer for the same number
of layers $l$ and width $w$.

In order to show the advantage of the higher regularity time discretization, we
use differnet loss functions during training. We consider the three loss functions
\begin{subequations}
  \label{eq:losses}
  \begin{align}
    \label{eq:lc1}
    l_{C^{1}}^{2}(\vct E,\vct{\hat E})&=\paran{\frac{1}{\abs{I}\abs{\mathcal{D}}}\int_{I}\int_{\mathcal{D}} \lvert
                                        \vct{\hat E}-\vct E\rvert^{2}\drv x\drv t}\,,\\
    \label{eq:lpt}
    l_{\partial t}^{2}(\vct E,\vct{\hat E})&=\frac{1}{N}\sum_{i=1}^{N} \abs{\vct{\hat E}_{i}- \vct E_i}^2+\frac{1}{N}\sum_{i=1}^{N} \abs{\partial_t \vct{\hat E}_{i}-\partial_t \vct E_i}^2\,,\\
    \label{eq:l2}
    l^{2}(\vct E,\vct{\hat E})&=\frac{1}{N}\sum_{i=1}^{N} \abs{ \vct{\hat E}_{i}- \vct E_i}^2\,.
  \end{align}
\end{subequations}
Here \(l_{C^{1}}\) is motivated by the higher order time discretization and to
evaluate the integrals we integrate the Hermite-type polynomials on the
subintervals analytically. The second loss makes use of the data provided by the
higher order time discretization but only considers the error in the collocation
points (the subinterval endpoints). Since the losses themselves are difficult to
compare, we study the errors \(e_{\vct E}=\vct{\hat E}(x,\,t)-\vct E(x,\,t)\) in
the norms given in~\eqref{eq:norms} with the abbreviations
\(L^{\infty}\text{-}L^{2}(\vct E)\) and \(L^2\text{-}L^2(\vct E)\). In
Fig.~\ref{fig:so-convergence}, we evaluate the networks on successively refined
time meshes in line with a numerical convergence test. For each refinement we
use a new ANN which is trained as mentioned at the beginnig of this section. The
errors are then evaluated during the testing of the ANNs, which we also
described above.

For the GRU the test results are stable, despite the high oscillations observed
during the training time. Furthermore, all three loss functions lead to similar
results. The GRUs do not benefit from the two loss functions which include the
time derivative.

The FNO on the other hand profits from the added information and is otherwise
stuck at high errors. Even at high time resolution, the network encounters
difficulty in distinguishing frequencies that are just \(\SI{1}{\tera\hertz}\)
apart. Including the time derivative via the loss functions~\eqref{eq:lc1}
or~\eqref{eq:lpt} is important for effectively training the network. However,
the difference between them is negligible. Interestingly we are able to observe
linear convergence for both networks. 

\subsubsection{Computational efficiency of the solution operator}
\label{sec:org3158f28}
For the training of the ANNs, we implemented a distributed training algorithm
similar to the one in~\autocite{lianCanDecentralizedAlgorithms2017}. In contrast
to~\autocite{lianCanDecentralizedAlgorithms2017}, we sychronize the network
parameters by averaging them over all processes. For the available resources of
5 nodes, each with 2 GPUs, we are unable to determine any significant
performance gains from considering only the neighboring MPI processes. Our
distributed implementation is not equivalent to the sequential implementation
due to the synchronization of the network parameters and not the gradients,
which reduces the computational overhead. In the tests we run in this work,
there is no disadvantage to this approach, yielding the same accuracies up to
machine precision. In \autocite{lianCanDecentralizedAlgorithms2017} the authors
show that their decentralized algorithm, which is related to our approach, leads
to the same convergence rate as the vanilla SGD.

\begin{figure}[htbp]
  \centering
  \includegraphics[width=\linewidth]{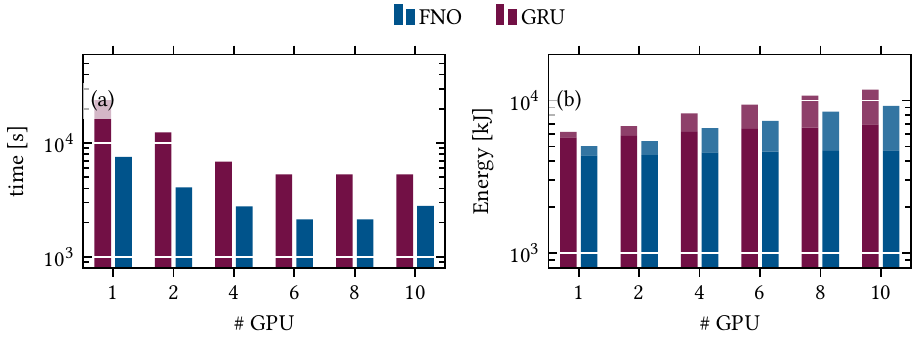}
  \caption{Runtime and energy consumption of GRUs and FNOs tested on up to 5
    nodes of an HPC cluster with 2 NVidia A100 each. The tests were conducted
    with an FNO model of size ($8\times 16$) and a GRU of size ($8\times 16$).
    The datasets are the largest presented in the solution operator
    tests.}\label{fig:rt-energy}
  \vspace*{10pt}
  \includegraphics[width=\linewidth]{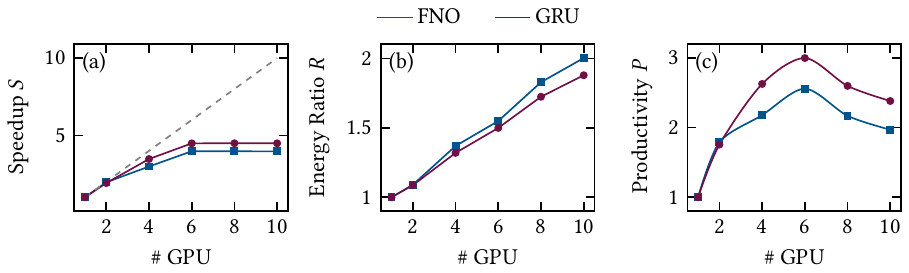}
  \caption{The speedup $S$, energy ratio $R$ and productivity $P$ calculated
    from the data of Fig.~\ref{fig:rt-energy} for our scaling
    tests.}\label{fig:srp}
  \vspace*{10pt}
\let\mc\multicolumn%
\begin{minipage}{0.45\textwidth}
\begin{center}
\begin{tabular}{crrrr}
\toprule
  \(l\) $\backslash$ \(w\) & \mc{2}{c}{8} & \mc{2}{c}{16}\\[0pt]
                           & SO & OCP & SO & OCP\\[0pt]  
\cmidrule(r{.66em}){1-1}\cmidrule(lr{.66em}){2-3}\cmidrule(lr{.66em}){4-5}
4 & \num{11609} & \num{130.37} & \num{11678} & \num{176.40}\\[0pt]
8 & \num{11647} & \num{223.73} & \num{11646} & \num{199.06}\\[0pt]
\bottomrule
\end{tabular}
\end{center}
\end{minipage}
\begin{minipage}{0.45\textwidth}
\begin{center}
\begin{tabular}{crrrr}
\toprule
  \(l\) $\backslash$ \(w\) & \mc{2}{c}{8} & \mc{2}{c}{16}\\[0pt]
                           & SO & OCP & SO & OCP\\[0pt]
\cmidrule(r{.66em}){1-1}\cmidrule(lr{.66em}){2-3}\cmidrule(lr{.66em}){4-5}
4 & \num{23796} & \num{153.18} & \num{23756} & \num{145.81} \\[0pt]
8 & \num{23652} & \num{151.23} & \num{23764} & \num{169.07} \\[0pt]
\bottomrule
\end{tabular}
\end{center}
\end{minipage}
\captionof{table}{\label{tab:ocp-wall-time}The wall time for the different stages of our
  proposed method in seconds: the solution operator training (SO) and the solution of the
  OCP. The table on the left contains results for the FNO and the one on the
  right contains results for the GRU for different parametrizations: The number
  of layers \(l\) and the width of the layers \(w\).}
\end{figure}

Fig.~\ref{fig:rt-energy} shows the strong scaling of the algorithm and the
corresponding energy consumption for the two network architectures under
consideration. The tests are run on an HPC cluster with 5 GPU nodes, each with 2
Nvidia A100 GPUs and 2 Intel Xeon Platinum 8360Y CPUs. The scaling tests are
performed with 1 GPU as a baseline and with 1 to 5 nodes, always using both GPUs
on the node. The number of MPI processes are equal to the number of GPUs such
that one MPI process uses one GPU and CPU. We further note that, as shown in
Table~\ref{tab:ocp-wall-time}, the training and evaluation times are equal
across one architecture, although the number of parameters differ significantly
(cf.~Table~\ref{orgf123153}). However the implementation in PyTorch is optimized
for larger networks and the ones we use are too small, to make a difference in
computation times.

Fig.~\ref{fig:rt-energy}
and~\ref{fig:srp}\,(a) illustrate the near-optimal scaling performance for up to
4 GPUs. For 6 GPUs, the impact of synchronization costs becomes noticeable, as
shown in Fig.~\ref{fig:srp}\,(a) and saturates afterwards. A further comparison
of our implementation with the asynchronous implementation available only in
PyTorch's Python interface would require significant effort since we have
exclusively used PyTorch's \texttt{C++} interface. Such a comparison for the
assessment of our implementation is beyond the scope of this work.
Here, we concentrate on evaluating the strong scaling test by means of the
speedup \(S\) and energy ratio \(R\),
\begin{subequations}
  \begin{multicols}{2}
    \noindent
    \begin{equation}
      S=\frac{t_{\text{wall}}(1)}{t_{\text{wall}}(n)}\,,
    \end{equation}
    \begin{equation}
      R=\frac{E(n)}{E(1)}\,,
    \end{equation}
  \end{multicols}
\end{subequations}
\noindent
where \(n\) is the number of GPUs (which coincides with the number of MPI
processes in this study), \(E(n)\) is the energy consumed by the CPU, memory and
GPU and \(t_{\text{wall}}\) is the wall clock time. The energy consumption of
the CPU and memory is almost constant, with the increase in energy consumption
primarily attributable to additional GPUs. Furthermore, the costs for CPU and
memory are high and the energy consumption of the GPUs only became larger when
using 10 GPUs. Overall our implementation exhibits great performance for this
small artificial problem. This is confirmed by the productivity
metric \(P\) in Fig.~\ref{fig:srp} (c), which is defined as the ratio of \(S\)
and \(R\), as per \autocite{anselmannEnergyefficientGMRESMultigridSolver2023}.
The optimal productivity lies at 6 GPUs and the results are promising for
scaling to larger problems.

\subsubsection{Optimal control through Deep Neural Networks}
\label{sec:org9e82767}
\begin{figure}[htbp]
  \centering {\small
    \begin{tabular}{lllllllll}
      \toprule
      Parametrization & \(\tau\,[\si{\second}]\) & \(P_1\) & \(a_{1}\) & \(f_1\,[\si{\tera\hertz}]\) & \(P_3\) & \(a_3\) & \(f_3\,[\si{\tera\hertz}]\) & \(\mathcal{J}\)\\[0pt]
      \midrule
      Initial values & \num{1.0e-14} & 6 & 1 & \num{291.56} & 6 & 1 & \(2\cdot 291.56\) & 0\\[0pt]
      GRU, \(4\times 8\) & \num{9.9e-13} & 5.99 & 0.02654 & \num{291.56} & 5.99 & 0 & \(2\cdot 291.56\) & \num[round-mode=places,round-precision=3,exponent-mode=scientific,print-zero-exponent=true,tight-spacing=true]{0.00017036065350346685}\\[0pt]
      GRU, \(8\times 8\) & \num{9.9e-13} & 5.99 & 0.02654 & \num{291.56} & 5.99 & 0 & \(2\cdot 291.56\) & \num[round-mode=places,round-precision=3,exponent-mode=scientific,print-zero-exponent=true,tight-spacing=true]{0.0001741007695254013}\\[0pt]
      GRU, \(4\times 16\) & \num{9.9e-13} & 5.99 & 0.02654 & \num{291.56} & 5.99 & 0 & \(2\cdot 291.56\) & \num[round-mode=places,round-precision=3,exponent-mode=scientific,print-zero-exponent=true,tight-spacing=true]{0.0001693222032205083}\\[0pt]
      GRU, \(8\times 16\) & \num{9.9e-13} & 5.99 & 0.02654 & \num{291.56} & 5.99 & 0 & \(2\cdot 291.56\) & \num[round-mode=places,round-precision=3,exponent-mode=scientific,print-zero-exponent=true,tight-spacing=true]{0.0001740310819512365}\\[0pt]
      FNO, \(4\times 8\) & \num{9.9e-13} & 5.99 & 2.8339 & \num{291.56} & 5.99 & 0 & \(2\cdot 291.56\) & \num[round-mode=places,round-precision=3,exponent-mode=scientific,print-zero-exponent=true,tight-spacing=true]{5.317621534240165e-6}\\[0pt]
      FNO, \(8\times 8\) & \num{9.9e-13} & 5.99 & 2.8338 & \num{291.56} & 5.99 & 0 & \(2\cdot 291.56\) & \num[round-mode=places,round-precision=3,exponent-mode=scientific,print-zero-exponent=true,tight-spacing=true]{5.275819862406618e-6}\\[0pt]
      FNO, \(4\times 16\) & \num{9.9e-13} & 5.99 & 2.8339 & \num{291.56} & 5.99 & 0 & \(2\cdot 291.56\) & \num[round-mode=places,round-precision=3,exponent-mode=scientific,print-zero-exponent=true,tight-spacing=true]{5.100505460091095e-6}\\[0pt]
      FNO, \(8\times 16\) & \num{9.9e-13} & 5.99 & 2.8337 & \num{291.56} & 5.99 & 0 & \(2\cdot 291.56\) & \num[round-mode=places,round-precision=3,exponent-mode=scientific,print-zero-exponent=true,tight-spacing=true]{5.1050882414502536e-6}\\[0pt]
      \bottomrule
    \end{tabular}
  }
  \captionof{table}{\label{tab:org688b584}The results for the test of the optimal control
    problem for different parametrizations of GRUs and FNOs after the training
    has converged. We show the parameters of 2 of the 4 pulses. All the
    parameters are trainable and therefore possibly affected by the updates
    during the optimization.}  
  \vspace*{10pt}
  \includegraphics[width=\linewidth]{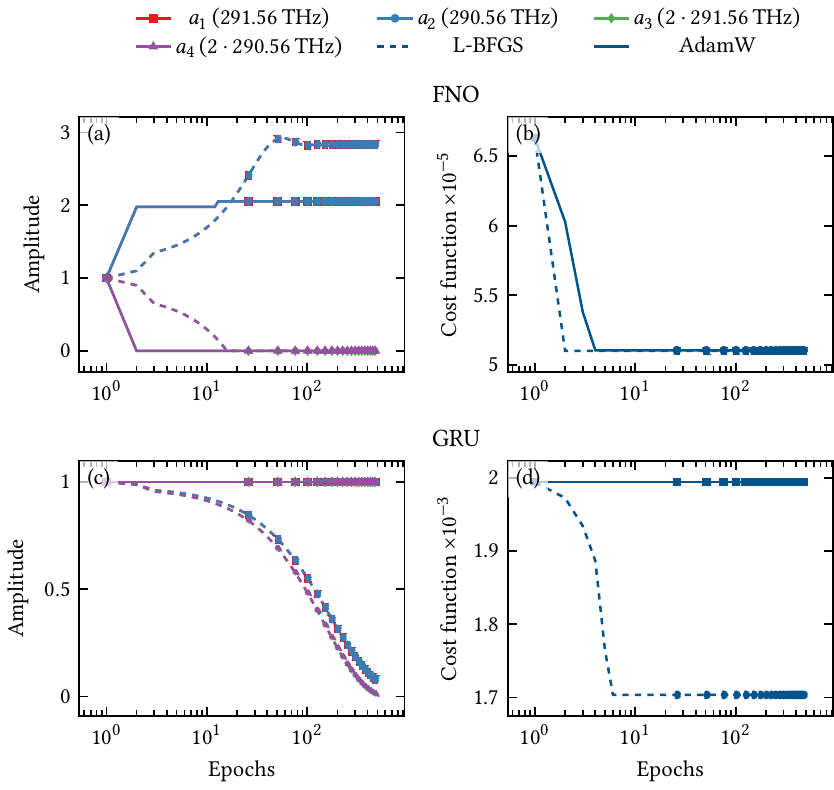}
  \captionof{figure}{The evolution of the amplitudes of the pulses (a, c) and the cost
    function (b, d) for FNOs (the first row of the subplots) and GRUs (the
    second row of the subplots). Solid lines correspond to the optimization with
    L-BFGS and dashed lines to the AdamW optimizer.}\label{fig:ocp-losses}
\end{figure}

In order to test the methodology we propose in Section~\ref{sec:orga943742}, we
construct a simple OCP based on the solution operator we obtained in the last
section. Initially, we sample 4 super Gaussian pulses parametrized as
in~\eqref{eq:samples} and define \(\Xi\) from~\eqref{eq:xin} accordingly
(\(n=4\)). We choose a supergaussian pulse of order \(p=6\), a full-width half
maximum of \(\tau=\SI{1.e-14}{\second}\) and the frequencies
\(f_1=\SI{291.56}{\tera\hertz},\,f_2=\SI{290.56}{\tera\hertz},\,f_3=2\cdot\SI{291.56}{\tera\hertz},\,f_4=2\cdot\SI{290.56}{\tera\hertz}\).
The time domain is chosen as \(I=[\num{-8e-15},\,\num{8e-15}]\) and the spatial
domain \(D=[\num{-8e-15}c_{0},\,\num{8e-15}c_{0}]\). In
Table~\ref{tab:org688b584} we list the initial parameters $\xi_0\in \Xi$ of the
pulse (cf.~Algorithm~\ref{alg:opcondl}, line~\ref{alg:init}) for all tests
performed in this section, except for
\(\varphi_1=\varphi_2=\varphi_3=\varphi_4=\zeta_1=\zeta_2=\zeta_3=\zeta_4=0\),
since they showed low sensitivity.

The setting for the OCP is described in
Problem~\ref{problem:opconTHz}. We choose the cost function such that the amplitudes
of the two high frequencies \(f_{3}\), \(f_4\) are minimized. To this end, we
put \(f_{\Omega}=\frac{1}{2}(f_3+f_4)\) and \(r=\SI{10}{\tera\hertz}\)
in~\eqref{eq:cost-fn-psi}. The solution to this minimization problem is trivial,
as setting the amplitude of the two high frequencies to zero would be sufficient
to solve the problem. For the ANNs exploiting the linearity is not
straightforward: The solution operator has to generalize from plane waves to
Gaussian pulses and is a nonlinear operator by construction. Therefore,
linearity has to be learned from the training data, and we cannot assume that we
always achieve that.

We compare two different optimization methods: the AdamW optimizer
\autocite{kingmaAdamMethodStochastic2015}, a modified version of stochastic
gradient descent, and the L-BFGS~\autocite{liuLimitedMemoryBFGS1989}, a
quasi-Newton method. In Fig.~\ref{fig:ocp-losses} we plot the development of the
amplitudes and cost function over the epochs. The first row of plots contains
the results for the FNOs. The L-BFGS method converges in 1 step and only changes
slightly afterwards. However, the two optimization routines don't lead to the
same result. The L-BFGS gets stuck in the first local minimum it finds, the AdamW
optimizer does not due to the added momentum.
We can rule out the penalty parameter as the reason. We set it to
\(\beta=\num{6.e-14}\) in all of our tests. Slight differences of other
parameters introduced during the optimization and differences in the output of
the ANN are further contributing factors. In Table~\ref{tab:org688b584} we
compare the final pulse parametrizations and the cost function of the OCP for
different networks using the AdamW optimizer. All parameters of the pulse are
subject to optimization and can be affected by updates during optimization. For
the FNO, the amplitude $a_1$ is correctly used as the main quantity to control the
optimization and other parameters show low sensitivity. The results for the
FNO are independent of the parametrization and reach nearly the same values with
respect to the parameters.

The second row of plots in Fig.~\ref{fig:ocp-losses} contains the results for
the GRUs in solving the OCP. The L-BFGS optimizer stagnates
and fails to improve the cost function, while the AdamW optimizer shows some
improvements. However, the solution operator resulting from the GRUs can't
distinguish between different frequencies. In Table~\ref{tab:org688b584} we
observe that the GRUs nearly remove all the lower frequencies. In the
previous section the accuracy and convergence behavior of the GRU and FNO were
almost the same. From the unsatisfactory results for the solution of the OCP
with GRUs we conclude that the GRU architecture is not
well-suited for this type of problem. Furthermore, the applicability of the ANNs to
Problem~\ref{problem:opconTHz} is predicted by the ability to approximate the
solution operator $\symbfcal{S}$, which is difficult to evaluate experimentally.
Nevertheless, FNOs are promising tools for this task based on their performance
in our experiments. Overall the
FNOs solve the test problem and are good candidates for the deployment to the
realistic problem from nonlinear optics.

\subsubsection{Computational efficiency of the optimal control algorithm}
For the evaluation of the ANN in our optimal control algorithm we only use one
GPU on a single node, since the distributed algorithm does not pay off in that
case due to the fast evaluation of the ANN. In Table~\ref{tab:ocp-wall-time} we
show the wall times for the training of the solution operator and the wall time
for the solution of the OCP. The wall time for the OCP is significantly lower
than the training of the solution operator. The low cost of solving the OCP we
is a considerable advantage, since we expect to reuse the trained solution operator
multiple times in the optimal control setting. Overall, the approach has great
potential for the solution of OCPs, since classical numerical solutions, even in
this artificial settings, exhibit high computational cost.
In the next section, we use highly accurate simulation data from a realistic
physical setting to train the solution operator and apply it to Problem
ref:problem:opconTHz, the OCP of maximizing THz generation.

\subsection{THz Generation in a Periodically Poled Nonlinear Crystal}\label{sec:org81ca8fd}
The main goal in this section is to show the potential of our proposed method by
applying it to a case where experimental results
are available~\autocite{olgunHighlyEfficientGeneration2022}, in order to verify that
the algorithm solves an OCP in realistic settings. That potentially leads to the
improvement of the experimental setup especially for higher intensities where
simplified models fail. As in~\autocite{margenbergAccurateSimulationTHz2023}, we
use 2 super Gaussian pulses parametrized as in~\eqref{eq:samples} and define
\(\Xi\) from~\eqref{eq:xin} accordingly (\(n=2\)). We choose a supergaussian
pulse of order \(p=6\), a full-width half maximum of
\(\tau=\SI{250}{\pico\second}\) and the frequencies
\(f_{1}=\SI{291.56}{\tera\hertz}\), \(f_{2}=\SI{291.26}{\tera\hertz}\). The
pulses are separated in center frequency by the THz frequency
\(f_{\Omega}=\SI{0.3}{\tera\hertz}\). In this section we choose a pulse with
average fluence of \(\SI{200}{\milli\joule\per\square\centi\metre}\). The
average fluence $F$ is defined as the mean of the optical intensity
$I_{\vct{E}}$ over time. The detailed definitions are given in
Appendix~\ref{sec:qoi}. In the simplified 1D case, they can be expressed as
\begin{subequations}
  \begin{multicols}{2}
    \noindent
    \begin{equation}
      \label{eq:intensity-1d}
      I_{\vct E}=\frac{1}{2}\varepsilon_0c_0 \abs{\sqrt{\varepsilon_r}*\vct E}^2\,,
    \end{equation}
    \begin{equation}
      \label{eq:fluence-1d}
      F=\frac{1}{T}\int_0^T I_{\vct{E}}(t)\drv t\,,
    \end{equation}
  \end{multicols}
\end{subequations}
\noindent
The pulse is applied at the left-hand side of the crystal by a Dirichlet boundary condition
on \(\Gamma_{\text{in}}\) (cf. Fig.~\ref{fig:periodsop}), propagates through the
domain and enters the PML where it is attenuated. The problem setting is already
sketched in Fig.~\ref{fig:periodsop}. The computational effort for these
simulations is high: The simulations presented here took 15 days on an HPC
cluster using 5 nodes, each with 2 Intel Xeon Platinum 8360Y CPUs. In this
study, we limited our investigations and numerical simulations to one spatial
dimension. This was necessitated by simulation times and the added complexity of
using PMLs in 2D and 3D. In the settings investigated here, the simplification
of reducing the simulations to one spatial dimension and neglecting the impacts
of the remaining spatial directions is not expected to significantly perturb the
results. The simulation results presented here are based on a timestep size of
\(k=\num{5.e-17}\) and average cell-size of \(\num{5.175e-8}\), which leads to
\(\num{5.e9}\) number of timesteps and \(\num{1703940}\) degrees of freedom in
space.

\subsubsection{Training and Evaluation of the Solution Operator}\label{sec:orga81bfd4}
\begin{figure}[htb]
  \centering
  \includegraphics[width=\linewidth]{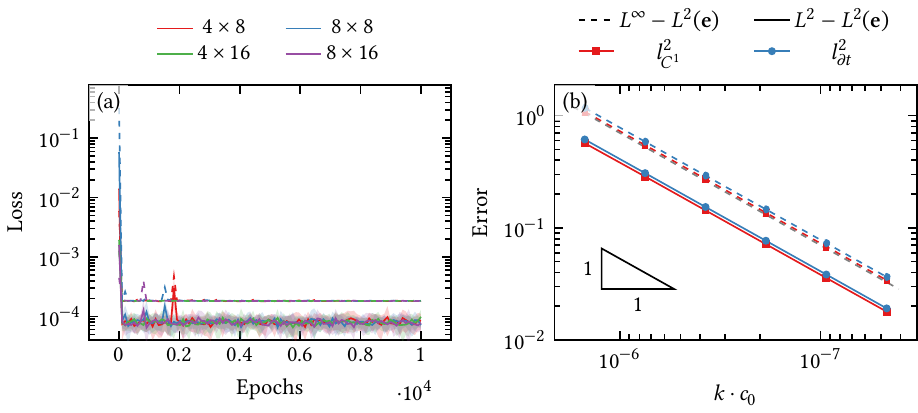}
  \caption{The average training (solid) and validation (dashed) loss for the FNO
    architecture (a). The right subplot (b) shows the convergence of the
    $L^{\infty}(L^2)$ (dashed) and $L^{2}(L^2)$ (solid) error for
    FNO.}\label{fig:so-losses-thz}
\end{figure}
\begin{figure}[htb]
  \centering
  \includegraphics[width=\linewidth]{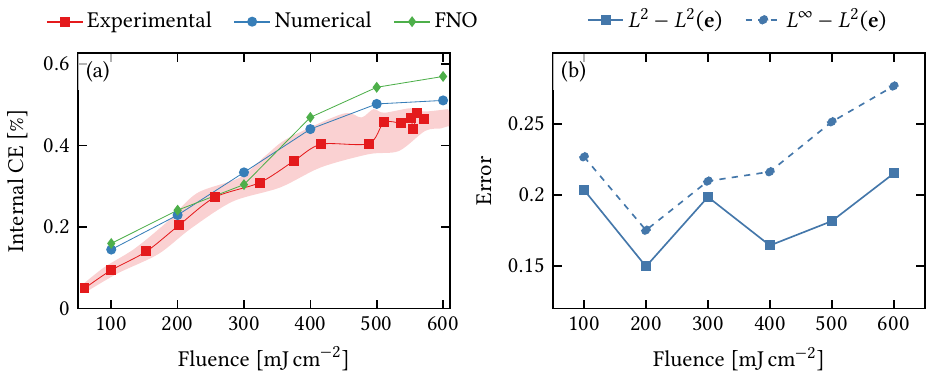}
  \caption{The left subplot shows the optical to THz conversion efficiency of
    the FNOs, simulation data and the experimental data
    from~\autocite{olgunHighlyEfficientGeneration2022}. For the experimental
    data we calculate an envelope from the standard deviation, that reflects the
    measurement inaccuracies. The right subplot shows the errors between FNOs
    and numerical solution. }\label{fig:so-convergence-thz}
\end{figure}
As in the case of artificial data in Section~\ref{sec:org9e82767}, the simulation data is used to train a
solution operator. We test only FNOs, since GRUs did not show satisfactory
results in the artificial test case. Fig.~\ref{fig:so-losses-thz}\,(a) shows the
losses of some parametrizations given in Table~\ref{orgf123153}. For the
training we split the data set obtained in the setting described above into a
training and validation set. We simulate 25 periods of the crystal, where the
first 15 periods are used for training and the last 10 periods are used for the
validation.

In Fig.~\ref{fig:so-losses-thz}\,(a) we observe fast convergence of the FNOs,
even with larger architectures, and the best models exhibit a similar loss
across all architectures. We use the two loss functions~\eqref{eq:lc1}
and~\eqref{eq:lpt}, since the added information of the time derivative proved to
be essential for good performance in the settings we investigate in this work.
In Fig.~\ref{fig:so-losses-thz}\,(b) we plot the errors
\(e=\vct{E}_{\text{GCC}^{1}(3)}-\vct{E}_{\text{FNO}}\). For different
timestep-sizes we use different FNOs, trained on simulation data obtained with
the same step size. We observe linear convergence as before in the artificial
test case.

We test the FNO on pulses \(g(t)\in \symbfcal{P}[\Xi]\) with different average
fluence \qtylist{100;200;300;400;500;600}{\milli\joule\per\square\centi\metre}.
Other parameters remain unchanged compared to the training scenario. We test the
different parametrizations of the FNOs (cf. Table~\ref{orgf123153}). We evaluate
the accuracy of the FNOs based on the internal conversion efficiency (CE) and
the errors~\eqref{eq:norms}. In Fig.~\ref{fig:so-convergence-thz}\,(a), we
compare the CE obtained from the FNO simulations with numerical simulations
obtained from the space-time finite element method presented here and
experimental results from~\autocite{olgunHighlyEfficientGeneration2022}. The
numerical simulations and FNOs are in good agreement with the experimental data and are
mostly close or within the standard deviation of the experimental data. The FNOs
perform very well, but lose some accuracy as the fluence increases. This is
expected, since we only trained it with data from simulations with a fluence of
\(\SI{200}{\milli\joule\per\square\centi\metre}\). Nevertheless, they seem to
learn the physical processes governing the THz generation accurately.

We analyze the accuracy of the FNO further in
Fig.~\ref{fig:so-convergence-thz}\,(b), where we plot the errors
\(e=\vct{E}_{\text{GCC}^{1}(3)}-\vct{E}_{\text{FNO}}\), evaluated in the
norms~\eqref{eq:norms} for different values of the average fluence. Although the
error grows with increasing fluence, for that the network was not trained,
the results are promising. The FNO is able to provide a good generalization
to pulses with higher fluence.

\subsubsection{Optimal control through the Solution Operator}
\label{sec:orgfe6e1b6}
As the final task, we consider the Problem~\ref{problem:opconTHz} in the
realistic setting and compare the results to experimental results. In this case
we want to maximize the radiation at the frequency
\(f_{\text{THz}}=\SI{0.3}{\tera\hertz}\). Therefore, we set
\(f_{\Omega}=f_{\text{THz}}\) and \(r=\SI{0.25}{\tera\hertz}\)
in~\eqref{eq:cost-fn-psi}. We test the method on the FNOs we trained in
Section~\ref{sec:orga81bfd4}. Fig.~\ref{fig:so-losses-thz} shows the losses on a
subset of the parametrizations given in Table~\ref{orgf123153}. The initial
parameters correspond to the case we considered in Section~\ref{sec:orga81bfd4}
for an average fluence of \(\SI{100}{\milli\joule\per\square\centi\metre}\). As
observed in~\ref{fig:so-convergence-thz}\,(a) the internal CE grows with
increasing fluence. In order to maximize the \(\SI{0.3}{\tera\hertz}\)-frequency
radiation the simplest improvement is an increase of the amplitude of the pulse.
We test here if this gets picked up by the FNO and it successfully optimizes
the internal CE. The internal CE is closely linked to the cost function of the
OCP, which is proportional to the intensity at \(\SI{0.3}{\tera\hertz}\).

\begin{figure}[htbp]
\centering
{\small
\begin{tabular}{llllllllll}
\toprule
Parametrization & \(\tau\,[\si{\pico\second}]\) & \(P_1\) & \(a_{1}\) & \(f_1\,[\si{\tera\hertz}]\) & \(P_2\) & \(a_2\) & \(f_2\,[\si{\tera\hertz}]\) & \(\mathcal{J}\) & \(F\) [\si{\milli\joule\per\square\centi\metre}]\\[0pt]
\midrule
Initial values & \num{250} & 6 & 1 & \num{291.26} & 6 & 1 & \(291.56\) & 0 & 100\\[0pt]
FNO, \(4\times 8\) & \num{257.268} & 6.57 & 3.7140 & \num{291.26} & 6.56 & 3.7080 & \(291.56\) & \num[round-mode=places,round-precision=3,exponent-mode=scientific,print-zero-exponent=true,tight-spacing=true]{2931.9749827289698} & 402\\[0pt]
FNO, \(8\times 8\) & \num{256.315} & 6.56 & 3.7257 & \num{291.26} & 6.56 & 3.6919 & \(291.56\) & \num[round-mode=places,round-precision=3,exponent-mode=scientific,print-zero-exponent=true,tight-spacing=true]{3266.13871331754825} & 402\\[0pt]
FNO, \(4\times 16\) & \num{306.852} & 6.59 & 3.7284 & \num{291.26} & 6.63 & 3.7332 & \(291.56\) & \num[round-mode=places,round-precision=3,exponent-mode=scientific,print-zero-exponent=true,tight-spacing=true]{2551.69235096332093} & 403\\[0pt]
FNO, \(8\times 16\) & \num{306.815} & 6.70 & 3.8010 & \num{291.26} & 6.71 & 3.7541 & \(291.56\) & \num[round-mode=places,round-precision=3,exponent-mode=scientific,print-zero-exponent=true,tight-spacing=true]{3300.54136064045696} & 406\\[0pt]
\bottomrule
\end{tabular}
}
\captionof{table}{\label{tab:org2948f27}The results for the test of the OCP for different parametrizations of FNOs after the training has converged. We show the relevant parameters of the pulses. The amplitude is normalized w.\,r.\,t.\ the initial value. The fluence \(F\) is displayed in the last column.}
\vspace*{10pt}
  \includegraphics[width=\linewidth]{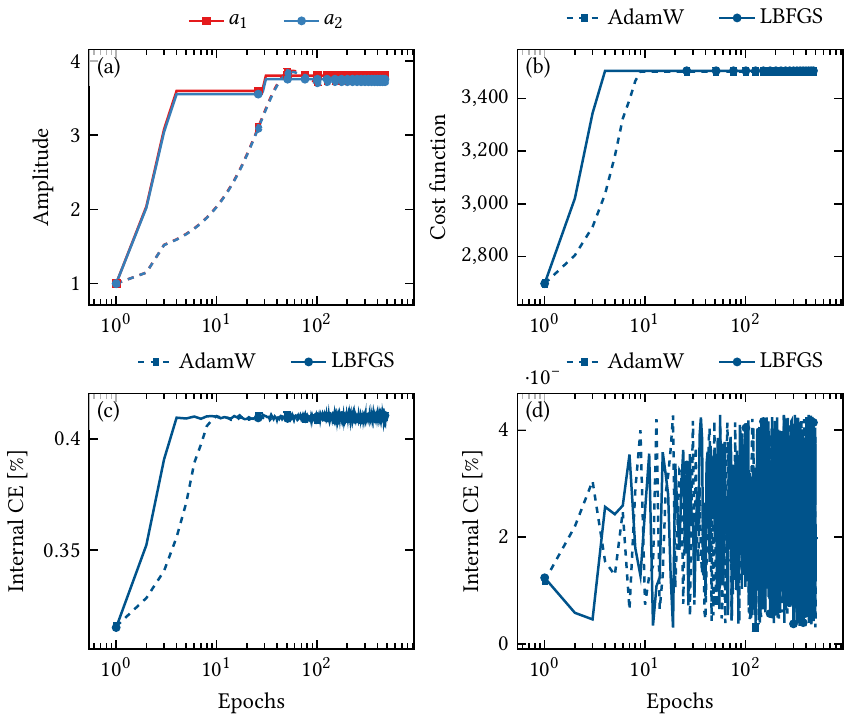}
  \captionof{figure}{The evolution of the amplitudes of the pulse (a). The Subplot (b)
    shows the evolution of the cost function, (c) the internal optical-to-THz CE
    and (d) internal CE for the second harmonic.}\label{fig:ocp-losses-thz}
\end{figure}
Comparing AdamW optimizer and L-BFGS for optimization methods, we plot the
amplitudes and cost function development over epochs in
Fig.~\ref{fig:ocp-losses-thz}\,(a), (b). Although the convergence is slower than
the artificial test case, the trajectories are overall similar. Again, L-BFGS
converges significantly faster than AdamW. Both reach similar optima, further
improvement may be limited by the regularization term. Lowering it beyond the
value we used before, lead to instabilities and implausible results. The reason
for the slower convergence of the L-BFGS method is the requirement of using a
lower learning-rate. Any attempt to use a higher learning-rate result in
stagnation.

The second row of plots in Fig.~\ref{fig:ocp-losses-thz} contains the internal
optical to THz CE and the internal CE of the second harmonic generation. The
optical to THz CE in Fig.~\ref{fig:ocp-losses-thz}\,(c) grows proportionally to the
cost function, which confirms that the FNO approximates at least a part of the
solution operator and physical model. In order to improve the internal CE, the
amplitude grows significantly, which confirms our expectation, that the
amplitude should be the main tuning parameter. We also observe this in Table~\ref{tab:org2948f27}.

In good agreement to our previous observations
in~\autocite{margenbergAccurateSimulationTHz2023}, the CE of the second harmonic
generation in Fig.~\ref{fig:ocp-losses-thz}\,(d) oscillates strongly. The reason
for the oscillations is the
phase-mismatch, which leads to oscillating negative and positive interference. This
leads to varying conversion efficiencies over the layers, depending on how close
we are to a phase-match.
In Table~\ref{tab:org2948f27}, we compare the final pulse parametrizations and
the cost function of the OCP for different networks using
the AdamW optimizer. Overall, FNOs perform well for the optimization of optical
to THz generation.

Performing the numerical simulations for obtaining the training data is the main
contributor to the computational costs. The subsequent training of the solution
operator takes 1 day and the final solution of the optimal control algorithm
takes 2 hours at most. Considering that a single numerical solution takes 15
days, the proposed approach offers great potential for optimal control problems
involving complex physics, in particular nonlinear optics. These problems are
still computationally challenging and oftentimes remain infeasible through classical
methods.

\FloatBarrier
\section*{Conclusion}
\label{sec:orgb7efbd1}
In this paper we developed methods to solve an optimal control problem arising
in nonlinear optics. To this end, we extended the Galerkin-collocation time discretization
to a nonlinear dispersive wave eqation. We observed that the method is
particularly well suited for problems arising in nonlinear optics. We confirmed
the results found in~\autocite{anselmannGalerkinCollocationApproximation2020}
by convergence tests. Although the implementation of the method is
parallelized and able to run on HPC platforms, the solution time for using them
within an optimal control loop is still too high.

We devised an algorithm which uses the simulation data with discrete solutions
of higher regularity in time to train an ANN, which is used for the forward solve.
The algorithm is applicable to a general optimal Dirichlet boundary control
problem and can be extended to other optimal control problems. Our method allows
for efficient solution of the optimal control problem, since we only require the
solution at some collocation points, and don't need the full space-time
solution. We compared GRUs and FNOs and tested their implementation on HPC
platforms and verified it by a strong scaling test. We also evaluated the energy
efficiency.
We were able to observe first order convergence,
which we only reached for the FNO with the added higher regularity. A thorough
investigation of this phenomenon is subject to future work. The GRU architecture
was not able to solve the optimal control problem satisfactorily despite its
good accuracy during the initial tests. FNOs were successful in solving the
optimal control problem. They clearly had an advantage over GRUs, since they are
designed for solving PDEs.

The optical to THz conversion efficiency achieved by the FNOs was found to be in
good agreement with experimental data. Moreover, the FNOs were successful in
optimizing the efficiency of this conversion process in an optimal boundary
control setting. Solving the whole optimal control problem with the trained
solution operator is 360 times faster than a single forward solve with the
numerical methods. Through its computational efficiency, the FNO has the
potential to enable breakthroughs in the development of high-field THz pulses,
by efficiently solving the optimization problem of maximizing the optical to THz
conversion efficiency. The rapid convergence of the training and fast evaluation
of the FNOs make them a cost-effective solution for this purpose and potentially
other complex physical problems.



\subsection*{Acknowledgement}
\label{sec:org588caf9}
NM acknowledges support by the Helmholtz-Gesellschaft grant number HIDSS-0002
DASHH. Computational resources (HPC-cluster HSUper) were provided by the project
hpc.bw, funded by dtec.bw — Digitalization and Technology Research Center of the
Bundeswehr. FXK acknowledges support througth ERC Synergy Grant (609920).

\printbibliography

\appendix
\section{Derivation of the fully discrete system}
\label{sec:org85a3097}
Here, we elaborate on the derivation of the fully discrete problems carried out
in Section~\ref{sec:org83f7b99}. We discretize Problem~\ref{problem:lor-stm} and
particularly describe the steps necessary to obtain the fully discrete, global
in time Problem~\ref{problem:lor-gcc} with the
equation~\eqref{eq:abstract-discrete-form} from~\eqref{eqwave:lor-nade2}. We
derive the local fully discrete problem, which, as discussed in
Section~\ref{sec:org5743c30}, leads to the global fully discrete
problem~\ref{problem:lor-gcc}. Finally, we describe the solution of the local,
fully discrete problems by a Newton linearization in combination with the
solvers for the arising linear systems of equations.

Following~\autocite{anselmannNumericalStudyGalerkin2020}, we define \(\{\symbf \phi_j\}_{j=1}^J
\subset \symbfcal V_h\) as a (global) nodal Lagrangian basis of \(\symbfcal V_h\) and
the Hermite-type basis of
\(\mathbb{P}_{3}(\hat I;\, \R)\), where $\hat I \coloneq [0,\,1]$:
\begin{equation*}
    \hat\xi_{0}(t)=1-3t^{2}+2t^{3}\quad \hat\xi_{1}(t)=t-2t^{2}+t^{3}\quad
    \hat\xi_{2}(t)=3t^{2}-2t^{3}  \quad  \hat\xi_{3}(t)=-t^{2}+t^{3}\,.\\
\end{equation*}
With the affine transformation
\begin{equation*}
  \mat{T}_{n}
  \colon \begin{cases}%
    \hat I &\to I_{n}\\
    \hat t &\mapsto t_{n-1} + (t_n-t_{n-1}) \hat t
  \end{cases}%
\end{equation*}
the basis \({\{\xi_{i}\}}_{i=0}^{3}\) on \(I_{n}\) is given by the composition of
\(\hat\xi_{l}\circ \mat{T}_{n}^{-1}\eqcolon \xi_{l}\) for \(l=0,\dots,\,3\).
Functions \(\vct w_{\tau,\,h} \in \mathbb{P}_{3}(I_n;\,\symbfcal{V}_{h})\) are thus represented as
\begin{equation}
  \label{gcc-rep-basis-fun}
  \vct w_{\tau,\,h}(x,\,t) =\sum_{i=0}^3 \vct w_{n,\,i}(x)\xi_{i}(t) = \sum_{i=0}^{3} \sum_{j=1}^{J} \vct
  w_{n,\,i,\,j}\symbf \phi_j(x)\xi_{i}(t),\qquad \text{for }(x,\,t)\in \Omega\times \bar{I_{n}}.
\end{equation}
We adopt the representation~\eqref{gcc-rep-basis-fun} for the variables
\(\vct U,\, \vct P,\, \vct E,\, \vct A\) and choose test functions from
\(\mathbb{P}_{0}(I_n;\,\symbfcal{V}_{h})\). A test basis of
\(\mathbb P_0 (I_n,\,\symbfcal{V}_{h})\) is then given by
\begin{equation}
  \label{eq:testbasis}
  \symbfcal B = \brc{\symbf \phi_i \vct 1_{I_n}}_{i=1}^{J}\,.
\end{equation}
Let \(\symbf A_h^0:
V_0 \to \symbfcal V_{h,\,0}\) be the discrete operator that is defined by
\begin{equation}
\label{eq:Ah}
\langle \symbf A_h \vct e_h ,\, \symbf \phi_h \rangle = \langle \symbf \nabla \vct e_h,\, \symbf \nabla \symbf \phi_h\rangle  \quad
\forall\;  \symbf \phi_h\in \symbfcal V_{h,\,0}\,.
\end{equation}
We define
$V_{g_h}\coloneq \brc{\vct v \in V\suchthat \vct v = g_h\; \text{on}\; \Gamma_D
}$ and $\mat A_h\colon V_{g_h}\to \symbfcal{V}_h$, $\vct w\mapsto \mat A_h\vct w$.
By the definition of $V_{g_h}$, $\vct w$ admits the representation $\vct w= \vct
w^0+g_h$ and we define $\mat A_h$ by
\begin{equation}
  \label{eq:Ah-g}
  \mat A_h\vct w=\mat A_h^0w^0+g_h\,. 
\end{equation}
For $\vct w \in \{\vct u,\,\vct p,\,\vct a,\,\vct e\}$, we
denote the right and left-hand limit by
\begin{equation}
  \label{eq:limits}
  \partial_t^i\vct w_{n,\,h}^{-}=\lim_{t\nearrow t_n}\partial_t^i\vct w_{\tau,\,h}(t),\quad\partial_t^i\vct
  w_{n,\,h}^{+}=\lim_{t\searrow t_n}\partial_t^i\vct w_{\tau,\,h}(t),\quad
  \text{for}\; i\in \brc{0,\,1}.
\end{equation}
Recall the fully discrete, global formulation of the $\text{GCC}^{1}(3)$ method
Problem~\ref{problem:lor-gcc}. Now consider the local problem on the interval
\(I_n\) where the trajectories \(\vct e_{\tau,\,h}(t)\),
\(\vct a_{\tau,\,h}(t)\), \(\vct p_{\tau,\,h}(t)\), and \(\vct u_{\tau,\,h}(t)\)
have already been computed for all \(t \in [0,\,t_{n-1}]\) with initial
conditions \(\vct e_{\tau,\,h}(0)=\vct e_{0,\,h}\),
\(\vct a_{\tau,\,h}(0)=\vct a_{0,\,h}\), \(\vct p_{\tau,\,h}(0)=\vct p_{0,\,h}\)
and \(\vct u_{\tau,\,h}(0)=\vct u_{0,\,h}\). Then we solve the following local problem:
\begin{problem}{Local, fully discrete, GCC$^{1}$(3) method for (\ref{eq:variational-cont})}{lor-gcc-local}
  Given
  \((\vct e_{\tau,\,h}(t_{n-1}),\,\vct a_{\tau,\,h}(t_{n-1}),\,\vct
  p_{\tau,\,h}(t_{n-1}),\,\vct u_{\tau,\,h}(t_{n-1}))\in \symbfcal{V}_h^{4}\),
  find
  \((\vct e_{\tau,\,h},\,\vct a_{\tau,\,h},\,\vct p_{\tau,\,h},\,\vct
  u_{\tau,\,h}) \in \paran{\mathbb{P}_{3}(I_{n};\,\symbfcal{V}_{h})}^{4}\) such
  that $e_{\tau,\,h}=g_{\tau,\,h}^{\vct e}\:\text{on}\:\bar{I}_n\times \Gamma_D$
  and
  \begin{subequations}\label{eq:wavephy-collocation}
    \begin{align}
      \label{eq:collocation1}\vct w_{n-1,\,h}^{+}=\vct w_{n-1,\,h}^{-}\; \forall
      \vct w\in \{\vct e,\,\vct a,\,\vct p,\,\vct u\}&\,,\\
      \label{eq:collocation2}\partial_{t} \vct w_{n-1,\,h}^{+}=\partial_{t} \vct w_{n-1,\,h}^{-}\; \forall
      \vct w\in \{\vct e,\,\vct a,\,\vct p,\,\vct u\}&\,,\\
      \label{eq:collocation3}-\vct u_{n,\,h}^{-}+\partial_t \vct p_{n,\,h}^{-}+\Gamma_{0} \vct p_{n,\,h}^{-}&=0\,,\\
      \nu_{t}^{2} \vct p_{n,\,h}^{-} -
      \nu_{t}^{2}\varepsilon_\Delta\vct e_{n,\,h}^{-} + \partial_t\vct u_{n,\,h}^{-}&=0\,,\\
      \label{eq:collocation4}-\Gamma_{0}\vct p_{n,\,h}^{-} +
      \chi^{(2)}\partial_t\paran{\abs{\vct e_{n,\,h}^{-}}\vct e_{n,\,h}^{-}}+
      \varepsilon_{\omega}\partial_{t}\vct e_{n,\,h}^{-} -\vct a_{n,\,h}^{-}&=0\,,\\
      \label{eq:collocation5}\nu_{t}^{2}\varepsilon_\Delta\vct e_{n,\,h}^{-} + \symbf{A}_h\vct  e_{n,\,h}^{-}
      -\nu_{t}^{2}\vct p_{n,\,h}^{-} + \partial_{t}\vct a_{n,\,h}^{-}&=\vct f_{n,\,h}^{-}\,,
    \end{align}
  \end{subequations}
  and for all
  \((\symbf \phi_{\tau,\,h}^{0},\,\symbf \phi_{\tau,\,h}^{1},\,\symbf \phi_{\tau,\,h}^{2},\,\symbf \phi_{\tau,\,h}^{3}) \in
  \paran{\mathbb{P}_{0}(I_{n};\,\symbfcal{V}_{h,\,0})}^{4}\),
  \begin{subequations}\label{eq:wavephy-variational}
    \begin{align}
      \int_{t_{n-1}}^{t_{n}}\dd{\partial_{t}\vct p_{\tau,\,h}^{n}}{\symbf \phi_{\tau,\,h}^{0}} + \Gamma_{0} \dd{\vct p_{\tau,\,h}^{n}}{\symbf \phi_{\tau,\,h}^{0}} - \dd{\vct u_{\tau,\,h}^{n}}{\symbf \phi_{\tau,\,h}^{0}} \drv t&=0\,,\label{eq:wavephy-variational-a}\\
      \int_{t_{n-1}}^{t_{n}}\nu_{t}^2\dd{\vct p_{\tau,\,h}^{n}}{\symbf \phi_{\tau,\,h}^{1}} -
      \varepsilon_\Delta\nu_{t}^2\dd{\vct e_{\tau,\,h}^{n}}{\symbf \phi_{\tau,\,h}^{1}} + \dd{\partial_{t}\vct u_{\tau,\,h}^{n}}{\symbf \phi_{\tau,\,h}^{1}}\drv t&=0\,,\label{eq:wavephy-variational-b}\\
      \notag\int_{t_{n-1}}^{t_{n}}\varepsilon_{\omega} \dd{\partial_{t} \vct e_{\tau,\,h}^{n}}{\symbf \phi_{\tau,\,h}^{2}} - \Gamma_{0}
      \dd{\vct p_{\tau,\,h}^{n}}{{\symbf \phi_{\tau,\,h}^{2}}}\qquad\qquad&\\
      +\chi^{(2)}\dd{\partial_{t}(\abs{\vct e_{\tau,\,h}^{n}}\vct e_{\tau,\,h}^{n})}{\symbf \phi_{\tau,\,h}^{2}} - \dd{\vct a_{\tau,\,h}^{n}}{\symbf \phi_{\tau,\,h}^{2}}\drv t&=0\,,\label{eq:wavephy-variational-c}\\
      \notag\int_{t_{n-1}}^{t_{n}}\dd{\symbf \nabla \vct e_{\tau,\,h}^{n}}{\symbf \nabla\symbf \phi_{\tau,\,h}^{3}} + (\varepsilon_{\Omega}-
      \varepsilon_{\omega})\nu_{t}^{2}\dd{\vct e_{\tau,\,h}^{n}}{\symbf \phi_{\tau,\,h}^{3}}\qquad\qquad&\\
      -
      \nu_{t}^{2}\dd{\vct p_{\tau,\,h}^{n}}{\symbf \phi_{\tau,\,h}^{3}}+\dd{\partial_{t}\vct a_{\tau,\,h}^{n}}{\symbf \phi_{\tau,\,h}^{3}}\drv t&=\int_{t_{n-1}}^{t_{n}}\dd{\vct f_{\tau,\,h}}{\symbf \phi_{\tau,\,h}^{3}}\drv t\,.\label{eq:wavephy-variational-d}
    \end{align}
  \end{subequations}
\end{problem}
We comment on the local fully discrete problem~\ref{problem:lor-gcc-local}:
\begin{itemize}\itemsep1pt \parskip0pt \parsep0pt
\item Note that we evaluate the time integrals on the right-hand side
of~\eqref{eq:wavephy-variational-d} and the boundary conditions
$g_{\tau,\,h}^{\vct e}\in C^1(\bar{I};\,\symbfcal{V}_h)$ (cf.~Assumption~\ref{ass:inhom}) using the
Hermite-type interpolation operator \(I_{\tau}\restrict{I_n}\), on \(I_n\),
defined by
\begin{equation}
\label{eq:hermite_interpolation}
I_{\tau}\restrict{I_n}g(t)
=
\hat \xi_0  (0)
 g\restrict{I_n}(t_{n-1})
+
\tau_n \hat \xi_1  (0)
\partial_t g\restrict{I_n}(t_{n-1})
+
\hat \xi_2  (1)
 g\restrict{I_n}(t_{n-1})
+
\tau_n \hat \xi_3  (1)
\partial_t g\restrict{I_n}(t_{n-1})\,.
\end{equation}
\item The collocation conditions~\eqref{eq:collocation2} need to be defined at
  the initial time $t_0$. From the initial conditions we can get the collocation conditions~\eqref{eq:collocation2} at the initial
  timepoints by setting
  \begin{subequations}\label{eq:initial-collocation}
    \begin{align}
      \partial_t \vct u_{\tau,\,h}(t_0^{-})&=\nu_t^2\varepsilon_{\Delta}\vct
                                         e_{0,\,h}-\nu_t^2\vct p_{0,\,h}\,,\\
      \partial_t \vct p_{\tau,\,h}(t_0^{-})&=\vct u_{0,\,h}-\Gamma_0\vct p_{0,\,h}\,,\\
      \partial_t \vct a_{\tau,\,h}(t_0^{-})&=\nu_t^2\vct p_{0,\,h}-\mat A_h\vct
                                         e_{0,\,h}-\nu_t^2\varepsilon_{\Delta}\vct
                                         e_{0,\,h}\,,\\
      \partial_t \vct e_{\tau,\,h}(t_0^{-})&=\varepsilon_{\Delta}^{-1}(\vct a_{0,\,h} -
                                         \chi^{(2)}\partial_t(\abs{\vct
                                         e_{0,\,h}}\vct
                                         e_{0,\,h})+\Gamma_0\vct p_{0,\,h})\,.
    \end{align}
  \end{subequations}
\item Consider a time interval $I_l$, $l=2,\dots,\,N$. We previously solved the
  Problem~\ref{problem:lor-gcc-local} on $I_{l-1}$. At $t_{l-1}$ the
  collocation conditions~\eqref{eq:collocation3}--\eqref{eq:collocation5} are
  fulfilled. For~\eqref{eq:collocation4}, we see that
  \begin{align*}
    &-\Gamma_{0}\vct p_{l,\,h}^{+} +
      \chi^{(2)}\partial_t\paran{\abs{\vct e_{l,\,h}^{+}}\vct e_{l,\,h}^{+}}+
      \varepsilon_{\omega}\partial_{t}\vct e_{l,\,h}^{+} -\vct a_{l,\,h}^{+} -\paran{-\Gamma_{0}\vct p_{l,\,h}^{-} +
      \chi^{(2)}\partial_t\paran{\abs{\vct e_{l,\,h}^{-}}\vct e_{l,\,h}^{-}}+
      \varepsilon_{\omega}\partial_{t}\vct e_{l,\,h}^{-} -\vct a_{l,\,h}^{-}}\\
    &\qquad=\chi^{(2)}\abs{\vct e_{l,\,h}^{+}}\paran{\partial_t\vct
      e_{l,\,h}^{+}}-\chi^{(2)}\abs{\vct e_{l,\,h}^{-}}\paran{\partial_t\vct
      e_{l,\,h}^{-}}
      +\chi^{(2)}\vct
      e_{l,\,h}^{+}\paran{\partial_t\abs{\vct e_{l,\,h}^{+}}}-\chi^{(2)}\vct
      e_{l,\,h}^{-}\paran{\partial_t\abs{\vct e_{l,\,h}^{-}}}=0\,,
  \end{align*}
  by using~\eqref{eq:collocation1} and~\eqref{eq:collocation2} componentwise.
  The remaining conditions~\eqref{eq:collocation3}--\eqref{eq:collocation5}
  follow immediately. Therefore, upon solving Problem~\ref{problem:lor-gcc-local} on
  $I_{l}$ the equations
  \begin{subequations}\label{eq:collocation-other}
    \begin{align}
      \label{eq:collocation-o3}-\vct u_{l,\,h}+\partial_t \vct p_{l,\,h}+\Gamma_{0} \vct p_{l,\,h}&=0\,,\\
      \nu_{t}^{2} \vct p_{l,\,h} -
      \nu_{t}^{2}\varepsilon_\Delta\vct e_{l,\,h} + \partial_t\vct u_{l,\,h}&=0\,,\\
      \label{eq:collocation-o4}-\Gamma_{0}\vct p_{l,\,h} +
      \chi^{(2)}\partial_t\paran{\abs{\vct e_{l,\,h}}\vct e_{l,\,h}}+
      \varepsilon_{\omega}\partial_{t}\vct e_{l,\,h} -\vct a_{l,\,h}&=0\,,\\
      \label{eq:collocation-o5}\nu_{t}^{2}\varepsilon_\Delta\vct e_{l,\,h} + \symbf{A}_h\vct  e_{l,\,h}
      -\nu_{t}^{2}\vct p_{l,\,h} + \partial_{t}\vct a_{l,\,h}&=\vct f_{l,\,h}\,,
    \end{align}
  \end{subequations}
  hold. This justifies the notion of a collocation method and shows that, from
  the initial timepoint on, global $C^1$-regularity is achieved by enforcing it
  from time step to time step.
\end{itemize}
We put the equations of the proposed \(\text{GCC}^1(3)\) approach in
their algebraic forms. In the variational equations
\eqref{eq:wavephy-variational}, we use the representation
\eqref{gcc-rep-basis-fun} for each component of \((\vct e_{\tau,\,h},\,\vct
a_{\tau,\,h},\,\vct p_{\tau,\,h},\,\vct u_{\tau,\,h}) \in (\mathbb
P_3(I_n;\symbfcal{V}_h))^4\) and choose the piecewise constant test functions.
We interpolate the right-hand sides in~\eqref{eq:wavephy-variational} by
applying the Hermite interpolation and evaluate the arising time integrals
analytically. The collocation conditions~\eqref{eq:wavephy-collocation} can be
recovered in their algebraic forms by using the fact that the Hermite type
polynomials and their first derivatives vanish at the locations \(x=0\) and \(x=1\),
with the exceptions \(\xi_0(0)=1\), \(\partial_t\xi_1(0)=1\), \(\xi_2(1)=1\),
\(\partial_t \xi_3(1)=1\).

Given the local Problem
\ref{problem:lor-gcc-local} on the interval \(I_n\)
and~\eqref{gcc-rep-basis-fun}, we
introduce the abbreviations
\(\vct w_{h,\,i}=\vct w_{n,\,i}(x)\in \symbfcal{V}_h\) and
\(\vct w_i=\paran{\vct w_{n,\,i,\,0},\dots,\,\vct w_{n,\,i,\,J}}^{\top}\in \R^J\) for \(\vct w \in \brc{\vct
Q,\, \vct R,\, \vct U,\, \vct P,\, \vct E,\, \vct A}\).
\begin{subequations}
Further we define
\begin{equation}
\label{eq:soln-vecs}
\vct v_{h,\,r} = (\vct e_{0,\,h}~\vct e_{1,\,h}~\vct a_{0,\,h}~\vct a_{1,\,h})^{\top} \text{ and } \vct v_{h,\,l}=(\vct e_{2,\,h}~\vct
e_{3,\,h}~\vct a_{2,\,h}~\vct a_{3,\,h})^{\top}\,.
\end{equation}
\end{subequations}
Then we condense the system of equations such that
we solve for \(\vct e_{2,\,h}\), \(\vct e_{3,\,h}\), \(\vct a_{2,\,h}\), \(\vct a_{3,\,h}\).
Solving for the unknowns \(\vct p_{2,\,h}\), \(\vct p_{3,\,h}\), \(\vct u_{2,\,h}\), \(\vct
u_{3,\,h}\) reduces to simple vector identities in their algebraic form. We write the nonlinear system of
equations in variational form for each subinterval \(I_{n}\) as
\begin{equation}
\label{gcc-lin-block-sys}
\symbfcal{A}_{h,\,n}(\vct{v}_{h,\,l})(\symbf{\Phi})=\mat{F}_{h,\,n}(\symbf{\Phi};\,\vct{v}_{h,\,r})\quad
\forall \symbf{\Phi}\in \symbfcal{V}_{h}^{4}\,,
\end{equation}
where
\(\symbfcal{A}_{h,\,n}\colon
\symbfcal{V}_{h}^{4}\times\symbfcal{V}_{h}^{4}\to\R\) is a semilinear form and
\(\mat{F}_{h,\,n}(\symbf{\Phi};\,\vct{v}_r)\) the right-hand side. Then \(\symbfcal{A}_{h,\,n}\) and the functional \(\mat{F}_{h,\,n}\)
in~\eqref{gcc-lin-block-sys} are defined through
\begin{subequations}
  \begin{equation}
    \label{gcc-system-var-S}
    \symbfcal{A}_{h,\,n}(\vct{v}_{h,\,l})(\symbf{\phi})=\symbfcal{A}_{h,\,n}^{1}(\vct{v}_{h,\,l})(\symbf{\Phi})+\symbfcal{A}_{h,\,n}^{2}(\vct{v}_{h,\,l})(\symbf{\Phi})+\symbfcal{A}_{h,\,n}^{3}(\vct{v}_{h,\,l})(\symbf{\Phi})+\symbfcal{A}_{h,\,n}^{4}(\vct{v}_{h,\,l})(\symbf{\Phi})\,,
  \end{equation}
  with the components $\symbfcal{A}_{h,\,n}^{i}$, $i=1,\dots,\,3$. The
  components represent the block structure of the system of equations in algebraic
  form. The components given as
  \begin{align}
    \notag\symbfcal{A}_{h,\,n}^{1}(\vct{v}_{h,\,l})(\symbf{\phi})&=-\frac{\Gamma_0}{{\nu_{t}^{2}}}
                                                             \dd{\symbf{\symbf \nabla} \vct e_{h,\,2}}{\symbf{\symbf \nabla}\symbf{\phi}}-
                                                             \varepsilon_{\Delta} \Gamma_0\dd{ \vct e_{h,\,2}}{\symbf{\phi}}%
                                                             -\dd{\vct a_{h,\,2}}{\symbf{\phi}}\\
                                                           &\quad + \chi^{(2)} \dd{\abs{\vct e_{h,\,3}} \vct e_{h,\,3}}{\symbf{\phi}}
                                                             +\frac{{\varepsilon_{\omega}}}{k} \dd{\vct e_{h,\,3}}{\symbf{\phi}}%
                                                             -\frac{\Gamma_0}{k {\nu_{t}^{2}}}\dd{\vct a_{h,\,3}}{\symbf{\phi}}\,,\\
    \notag\symbfcal{A}_{h,\,n}^{2}(\vct{v}_{h,\,l})(\symbf{\phi})&=\varepsilon_{\Delta}\dd{\vct
                                                             e_{h,\,2}}{\symbf{\phi}}-\frac{ {{k}^{2}} {\nu_{t}^{2}}-12}{12
                                                             {\nu_{t}^{2}}}\dd{\symbf{\symbf \nabla} \vct
                                                             e_{h,\,2}}{\symbf{\nabla}\symbf{\phi}}%
                                                             +\frac{ k \Gamma_0+6}{k{\nu_{t}^{2}}}\dd{\vct a_{h,\,2}}{\symbf{\phi}}\\
                                                           &\quad -\frac{ k \Gamma_0+6 }{12{\nu_{t}^{2}}}%
                                                             \dd{\symbf{\nabla} \vct
                                                             e_{h,\,3}}{\symbf{\nabla}\symbf{\phi}}-\frac{\varepsilon_{\Delta} \left( k
                                                             \Gamma_0+6\right) }{12}\dd{\vct e_{h,\,3}}{\symbf{\phi}} -\frac{ {{k}^{2}}
                                                             {\nu_{t}^{2}}+6 k \Gamma_0+24}{12 k {\nu_{t}^{2}}}\dd{ \vct
                                                             a_{h,\,3}}{\symbf{\phi}}\,,\\
    \notag\symbfcal{A}_{h,\,n}^{3}(\vct{v}_{h,\,l})(\symbf{\phi})&=\frac{ k {\nu_{t}^{2}}+2
                                                             \Gamma_0 }{2 {\nu_{t}^{2}}}\dd{\symbf{\nabla} \vct
                                                             e_{h,\,2}}{\symbf{\nabla}\symbf{\phi}}+ \varepsilon_{\Delta} \Gamma_0\dd{ \vct
                                                             e_{h,\,2}}{\symbf{\phi}}%
                                                             +\frac{{{k}^{2}}{\nu_{t}^{2}}-12}{{{k}^{2}}{\nu_{t}^{2}}} \dd{ \vct
                                                             a_{h,\,2}}{\symbf{\phi}}\\
                                                           &\quad +\frac{\varepsilon_{\Delta}}{k}\dd{\vct e_{h,\,3}}{\symbf{\phi}}%
                                                             -\frac{ {{k}^{2}} {\nu_{t}^{2}}-12 }{12 k {\nu_{t}^{2}}}\dd{\symbf{\nabla}
                                                             \vct e_{h,\,3}}{\symbf{\nabla}\symbf{\phi}}+ \frac{k \Gamma_0+6 }{{{k}^{2}}
                                                             {\nu_{t}^{2}}}\dd{\vct a_{h,\,3}}{\symbf{\phi}}\,,\\
    \notag\symbfcal{A}_{h,\,n}^{4}(\vct{v}_{h,\,l})(\symbf{\phi})&=\frac{2{\varepsilon_{\omega}}-{k
                                                             \varepsilon_{\Delta} \Gamma_0}}{2}\dd{ \vct e_{h,\,2}}{\symbf{\phi}}-\frac{ k
                                                             \Gamma_0}{2 {\nu_{t}^{2}}}\dd{\symbf{\nabla} \vct
                                                             e_{h,\,2}}{\symbf{\nabla}\symbf{\phi}} -\frac{ k {\nu_{t}^{2}}+2 \Gamma_0 }{2
                                                             {\nu_{t}^{2}}}\dd{\vct a_{h,\,2}}{\symbf{\phi}}\\
                                                           &\quad +\chi^{(2)} \dd{\abs{\vct e_2}\vct e_2}{\symbf{\phi}}
                                                             +\frac{k \Gamma_0}{12 {\nu_{t}^{2}}}\dd{\symbf{\nabla} \vct
                                                             e_{h,\,3}}{\symbf{\nabla}\symbf{\phi}}%
                                                             +\frac{k \varepsilon_{\Delta} \Gamma_0}{12}\dd{\vct e_{h,\,3}}{\symbf{\phi}} +
                                                             \frac{k}{12}\dd{\vct a_{h,\,3}}{\symbf{\phi}}\,,
  \end{align}
  and with an analogous splitting of $\mat{F}_{h,\,n}$
  \begin{align}
    \label{gcc-rhs}
    \mat{F}_{h,\,n}^1(\symbf{\Phi};\,\vct{v}_{h,\,r})&= 0\,,\\
    \notag \mat{F}_{h,\,n}^2(\symbf{\Phi};\,\vct{v}_{h,\,r})&=-\frac{k \Gamma_0 +6}{12 {\nu_{t}^{2}}}(\dd{\symbf{\nabla}\vct
                                                        e_{h,\,1}}{\symbf{\nabla}\symbf{\phi}}+\dd{\symbf{\nabla}\vct
                                                        e_{h,\,0}}{\symbf{\nabla}\symbf{\phi}}) - \frac{
                                                        \varepsilon_{\Delta} \left( k \Gamma_0+6\right) }{12}(\dd{\vct
                                                        e_{h,\,1}}{\symbf{\phi}}+\dd{\vct e_{h,\,0}}{\symbf{\phi}})
                                                        \\
                                               &\phantom{=}
                                                 + \frac{k \Gamma_0+6}{k
                                                 {\nu_{t}^{2}}}\dd{\vct a_{h,\,0}}{\symbf{\phi}}
                                                 + \frac{k}{12}\dd{\vct u_{h,\,1}}{\symbf{\phi}}
                                                 +\frac{k}{2}\dd{\vct u_{h,\,0}}{\symbf{\phi}}
                                                 +\frac{1}{2} \dd{\vct p_{h,\,1}}{\symbf{\phi}}
                                                 +\dd{4 \vct p_{h,\,0}}{\symbf{\phi}}\,,\\
    \notag\mat{F}_{h,\,n}^3(\symbf{\Phi};\,\vct{v}_{h,\,r})&= \frac{\varepsilon_{\Delta}}{k}\dd{\vct
                                                       e_{h,\,1}}{\symbf{\phi}}-\frac{{k^{2}} {\nu_{t}^{2}}-12 }{12 k
                                                       {\nu_{t}^{2}}}\dd{\symbf{\nabla}\vct
                                                       e_{h,\,1}}{\symbf{\nabla}\symbf{\phi}} + \frac{6
                                                       \varepsilon_{\Delta}}{k}\dd{\vct e_{h,\,0}}{\symbf{\phi}}
                                                       -\frac{{k^{2}} {\nu_{t}^{2}}-12 }{2 k
                                                       {\nu_{t}^{2}}}\dd{\symbf{\nabla}\vct
                                                       e_{h,\,0}}{\symbf{\nabla}\symbf{\phi}}\\
                                               &\phantom{=}
                                                 + \frac{ {k^{2}}
                                                 {\nu_{t}^{2}}-12 }{{k^{2}} {\nu_{t}^{2}}} \dd{\vct a_{h,\,0}}{\symbf{\phi}}
                                                 +\dd{\vct u_{h,\,0}}{\symbf{\phi}}
                                                 -\frac{1}{k}\dd{\vct p_{h,\,1}}{\symbf{\phi}}
                                                 -\frac{6}{k}\dd{\vct p_{h,\,0}}{\symbf{\phi}}\,,\\
    \notag\mat{F}_{h,\,n}^4(\symbf{\Phi};\,\vct{v}_{h,\,r})&=                \frac{k \varepsilon_{\Delta} \Gamma_0+2
                                                       {\varepsilon_{\omega}}}{2}\dd{\vct
                                                       e_{h,\,0}}{\symbf{\phi}}+\frac{k \Gamma_0}{2
                                                       {\nu_{t}^{2}}}\dd{\symbf{\nabla}\vct
                                                       e_{h,\,0}}{\symbf{\nabla}\symbf{\phi}} + \frac{k \Gamma_0}{12
                                                       {\nu_{t}^{2}}}\dd{\symbf{\nabla}\vct
                                                       e_{h,\,1}}{\symbf{\nabla}\symbf{\phi}}+\frac{k
                                                       \varepsilon_{\Delta} \Gamma_0}{12}\dd{\vct e_{h,\,1}}{\symbf{\phi}}\\
                                               &\phantom{=}+\chi^{(2)}\dd{\abs{\vct e_0}\vct e_0}{\symbf{\phi}}
                                                 +\frac{k {\nu_{t}^{2}}-2 \Gamma_0 }{2 {\nu_{t}^{2}}} \dd{\vct
                                                 a_{h,\,0}}{\symbf{\phi}}
                                                 +\frac{k}{12} \dd{\vct a_{h,\,1}}{\symbf{\phi}}\,.
  \end{align}
\end{subequations}
As a result of the condensation we further get update equations for the
variables \(\vct p_{2},\:\vct u_{2},\:\vct p_{3},\:\vct u_{3}\)
\begin{subequations}
  \label{gcc-updates}
  \begin{align}
    \vct u_{2}&=-\frac{6\Gamma_0}{k\nu_{t}^{2}}\vct u_{0} - \frac{6}{k}\vct p_{0} +
                \frac{\Gamma_0}{k\nu_{t}^{2}} \vct u_{1} - \frac{1}{k}\vct p_{1} +
                \frac{\varepsilon_{\Delta}}{k}\left(6\vct e_{0}+\vct e_{1} -
                \Gamma_0 k \vct e_{2} - \vct e_{3}\right)\,, \\
    \vct p_{2}&= \frac{\Gamma_0}{\nu_{t}^{2}}\vct u_{2} + \frac{6}{k\nu_{t}^{2}}\vct u_{0} -
                \frac{6(k\Gamma_0 -2)}{k^{2}\nu_{t}^{2}}\vct p_{0} + \frac{1}{k\nu_{t}^{2}}\vct u_{1} -
                \frac{\Gamma_0}{k\nu_{t}^{2}} \vct p_{1} - \varepsilon_{\Delta} \vct e_{2}\,, \\
    \vct u_{3}&=k\nu_{t}^{2}\vct p_{2} - k\varepsilon_{\Delta}\nu_{t}^{2}\vct e_{2}\,,\\
    \vct p_{3}&=k\Gamma_0 \vct p_{2}-k \vct u_{2}\,.
  \end{align}
\end{subequations}

The common approach of handling the nonlinear problem is a linearization by means of
Newton's method. Let
\begin{equation}
  \label{gcc-disc-var}
  \vct{x}\in \symbfcal{V}_{h} \colon
  \symbfcal{A}_{h,\,n}(\vct{x})(\symbf{\phi})=\mat{F}(\symbf{\phi})\qquad\forall
  \symbf{\phi}\in\symbfcal{V}_{h}
\end{equation}
be the variational equation related to~\eqref{gcc-lin-block-sys}. Recall that
\(\symbfcal{A}_{h,\,n}(\bullet)(\bullet)\) of~\eqref{gcc-system-var-S} is a
semi-linear form which is linear in the second argument. We assume that it is
sufficiently differentiable by means of the Gateaux derivative
\(\symbfcal{A}_{h,\,n}' (\vct{x})(\delta \vct{x},\,\symbf{\phi})\coloneq
\frac{\drv}{\drv s} \symbfcal{A}_{h,\,n}(\vct{x}+\varepsilon \delta
\vct{x})(\symbf{\phi})\restrict{\varepsilon=0}\). \(\symbfcal{A}_{h,\,n}'\)
denotes the derivative of \(\symbfcal{A}_{h,\,n}\) at
\(\vct{x}\in \symbfcal{V}_{h}\) in direction
\(\delta\vct{x}\in \symbfcal{V}_{h}\). The Newton iteration for
solving~\eqref{gcc-disc-var} with an initial guess
\(\vct{x}_{0}\in \symbfcal{V}_{h}\) iterates for \(m=0,\dots\)
\begin{equation}
  \label{gcc-newton}
  \begin{aligned}
    \delta\vct{x}_{m}\colon
    \symbfcal{A}_{h,\,n}'(\vct{x}_{m-1})(\delta\vct{x}_{m},\,\symbf{\phi}) &=
                                                                             \symbf{F}(\symbf{\phi})-\symbfcal{A}_{h,\,n}(\vct{x}_{m-1})(\symbf{\phi})
                                                                             \qquad\forall \symbf{\phi}\in \symbfcal{V}_{h}\,,\\
    \vct{x}_{m}&\coloneq \vct{x}_{m-1} + \delta\vct{x}_{m}.
  \end{aligned}
\end{equation}
Next we apply the Newton scheme to the system~\eqref{gcc-lin-block-sys}. The Gateaux derivative
\(\symbfcal{A}_{h,\,n}' (\vct{x})(\delta\vct{x},\,\symbf{\phi})\) is
\begin{equation}
  \label{gcc-gateaux}
  \begin{split}
    {\symbfcal{A}_{h,\,n}^1}' (\vct{x})(\delta\vct{x},\,\symbf{\phi})&= -\frac{\Gamma_0}{{\nu_{t}^{2}}} \dd{\symbf{\nabla}\delta
                                                                   \vct e_{h,\,2}}{\symbf{\nabla}\symbf{\phi}}- \varepsilon_{\Delta} \Gamma_0\dd{\delta
                                                                   \vct e_{h,\,2}}{\symbf{\phi}}%
                                                                   -\dd{\delta \vct a_{h,\,2}}{\symbf{\phi}}%
    \\
                                                                 &\quad+ \chi^{(2)} (\dd{\abs{\vct e_{h,\,3}}\delta \vct e_{h,\,3}}{\symbf{\phi}} +
                                                                   \dd{\abs{\delta \vct e_{h,\,3}} \vct e_{h,\,3}}{\symbf{\phi}})
                                                                   +\frac{{\varepsilon_{\omega}}}{k} \dd{\delta \vct e_{h,\,3}}{\symbf{\phi}}%
                                                                   -\frac{\Gamma_0}{k {\nu_{t}^{2}}}\dd{\delta \vct a_{h,\,3}}{\symbf{\phi}}\,,
    \\
    \notag{\symbfcal{A}_{h,\,n}^2}' (\vct{x})(\delta\vct{x},\,\symbf{\phi})&= \varepsilon_{\Delta}\dd{\delta \vct e_{h,\,2}}{\symbf{\phi}}-\frac{ {{k}^{2}}
                                                                         {\nu_{t}^{2}}-12}{12 {\nu_{t}^{2}}}\dd{\symbf{\nabla}\delta
                                                                         \vct e_{h,\,2}}{\symbf{\nabla}\symbf{\phi}}%
                                                                         +\frac{ k \Gamma_0+6}{k{\nu_{t}^{2}}}\dd{\delta  \vct a_{h,\,2}}{\symbf{\phi}}
    \\
                                                                 &\quad-\frac{ k \Gamma_0+6 }{12{\nu_{t}^{2}}}%
                                                                   \dd{\symbf{\nabla}\delta
                                                                   \vct e_{h,\,3}}{\symbf{\nabla}\symbf{\phi}}-\frac{\varepsilon_{\Delta} \left( k
                                                                   \Gamma_0+6\right) }{12}\dd{\delta \vct e_{h,\,3}}{\symbf{\phi}} -\frac{
                                                                   {{k}^{2}} {\nu_{t}^{2}}+6 k \Gamma_0+24}{12 k {\nu_{t}^{2}}}\dd{\delta
                                                                   \vct a_{h,\,3}}{\symbf{\phi}}\,,\\
    \notag{\symbfcal{A}_{h,\,n}^3}' (\vct{x})(\delta\vct{x},\,\symbf{\phi})&= \frac{ k {\nu_{t}^{2}}+2 \Gamma_0 }{2
                                                                         {\nu_{t}^{2}}}\dd{\symbf{\nabla}\delta
                                                                         \vct e_{h,\,2}}{\symbf{\nabla}\symbf{\phi}}+ \varepsilon_{\Delta} \Gamma_0\dd{\delta
                                                                         \vct e_{h,\,2}}{\symbf{\phi}}%
                                                                         +\frac{{{k}^{2}}{\nu_{t}^{2}}-12}{{{k}^{2}}{\nu_{t}^{2}}} \dd{\delta
                                                                         \vct a_{h,\,2}}{\symbf{\phi}}%
    \\
                                                                 &\quad+ \frac{\varepsilon_{\Delta}}{k}\dd{\delta \vct e_{h,\,3}}{\symbf{\phi}}%
                                                                   -\frac{ {{k}^{2}} {\nu_{t}^{2}}-12 }{12 k
                                                                   {\nu_{t}^{2}}}\dd{\symbf{\nabla}\delta
                                                                   \vct e_{h,\,3}}{\symbf{\nabla}\symbf{\phi}}+ \frac{k \Gamma_0+6 }{{{k}^{2}}
                                                                   {\nu_{t}^{2}}}\dd{\delta \vct a_{h,\,3}}{\symbf{\phi}}\,,
    \\
    \notag{\symbfcal{A}_{h,\,n}^4}' (\vct{x})(\delta\vct{x},\,\symbf{\phi})&= \frac{2{\varepsilon_{\omega}}-{k \varepsilon_{\Delta} \Gamma_0}}{2}\dd{\delta
                                                                         \vct e_{h,\,2}}{\symbf{\phi}}-\frac{ k \Gamma_0}{2
                                                                         {\nu_{t}^{2}}}\dd{\symbf{\nabla}\delta
                                                                         \vct e_{h,\,2}}{\symbf{\nabla}\symbf{\phi}} -\frac{ k {\nu_{t}^{2}}+2 \Gamma_0
                                                                         }{2
                                                                         {\nu_{t}^{2}}}\dd{\delta \vct a_{h,\,2}}{\symbf{\phi}}+\frac{k}{12}\dd{\delta \vct a_{h,\,3}}{\symbf{\phi}} \\
                                                                 &\quad+ \chi^{(2)} (\dd{\abs{ \delta \vct e_{h,\,2}} \vct e_{h,\,2}+\abs{\symbf{e_{h,\,2}}} \delta \vct e_{h,\,2}}{\symbf{\phi}})
                                                                   +\frac{k \Gamma_0}{12 {\nu_{t}^{2}}}(\dd{\symbf{\nabla}\delta
                                                                   \vct e_{h,\,3}}{\symbf{\nabla}\symbf{\phi}}+\dd{\delta \vct e_{h,\,3}}{\symbf{\phi}})\,.%
  \end{split}
\end{equation}
In every Newton step we have to solve a linear system of equations, for which we
use the GMRES method with an algebraic multigrid solver, which serves as a
preconditioner with a single sweep for every GMRES iteration. This accelerates
the convergence of the GMRES iterations.

\section{Extension of the neural operator to a full space-time approximation}
\label{sec:orgb4fe698}
 We sketch how the solution operator \(\symbf{U}\) can be extended in
order to obtain an approximation to \(\symbfcal{S}_h\). The idea is to
construct an interpolation operator \(\mat I_h\) which interpolates the solutions
\(\symbfcal{U}^{i}(\vct v)\) and \(\symbfcal{U}^{i+1}(\vct
v)\) in space. We need a finite-dimensional subspace \(W \subset V\) of dimension
\(M = \dim(W)\). A standard example would be \(W=V_h\), a classical finite element
space, neural networks are feasible as well. First note that from Eq.
\eqref{eq:collocation-points} we see that \(J_{\mathcal{D},\,i+1}\) is the set
\(J_{\mathcal{D},\,i}\), translated in positive \(x_1\)-direction. For \(i\in
\brc{1,\dots,\,m}\) we choose \(x_{a}\in J_{\mathcal{D},\,i}\) and let
\(\symbf{U}^{i-1}(g(t))=\vct p(t,\,x_{a})\). Then
\(\symbf{U}(\vct p(t,\,x_{a}))=\vct{\hat{p}} (t,\,x_b)\) where
\(x_{b}\in J_{\mathcal{D},\,i+1}\). We choose a basis
\(\brc{\phi_{1},\dots,\,\phi_J}\) of \(W\) with support points \(J_{\mathcal{D}}\).
Let \(\phi_j^a\) denote basis functions with support points in
\(J_{\mathcal{D},\,i+1}\) and \(\phi_j^b\) those with support points in
\(J_{\mathcal{D},\,i}\). On each subinterval \(I_n\) we can then interpolate
\begin{equation}
\label{eq:interp}
\symbf{U}\restrict{I_n}(x,\,t)
=\sum_{k=0}^{3}\xi_k(t)
\paran{\sum_{j=0}^{\abs{J_{\mathcal{D},\,i}}} \vct p_k\phi_j^a(x)+\vct{\hat{p}}_k\phi_j^b(x)}\,.
\end{equation}
Finally the solution operator \(\symbfcal{S}_h(g(t))\) can be approximated by
recursive application of \(\symbf{U}\) to itself,
i.\,e.\
\(\symbf{U}\circ\cdots\circ\symbf{U}(g(t))\eqcolon\symbfcal{U}\approx\symbfcal{S}_h(g(t))\).
\section{Physical Quantities and Quantites of Interest}
\label{sec:qoi}
We give some background on the quantities of interest in the simulations carried
out in this work. To this end we first introduce the Poynting vector $\vct S$
and optical power $P$ (cf.~\autocite[Chapter 6, Section 6]{jacksonClassicalElectrodynamics1999}).
\begin{multicols}{2}
  \noindent
  \begin{equation}
    \label{eq:poynting}
    \vct S = \frac{1}{\mu_0}\vct E \times \vct B\,,
  \end{equation}
  \begin{equation}
    \label{eq:op-power}
    P=\int_A \vct s\cdot n \drv s
  \end{equation}
\end{multicols}
\noindent
The flow of energy in an electromagnetic field, with the electric field $\vct E$
and the magnetic field $\vct B$ is described by the Poynting
vector~\eqref{eq:poynting}. The optical power $P$ in~\eqref{eq:op-power} is the
flux of the Poynting vector through a surface $A$. Then, the intensity is the
magnitude of the Poynting vector and the average fluence is the mean of the
intensity over time.
\begin{multicols}{2}
  \noindent
  \begin{equation}
    \label{eq:intense}
    I_{\vct E}=\abs{\vct S}=\frac{P}{\abs{A}}.
  \end{equation}
  \begin{equation}
    \label{eq:op-fluence}
    F=\frac{1}{T}\int_0^T  I_{\vct E}\drv t
  \end{equation}
\end{multicols}
\noindent
In the simplified 1D case, the optical intensity $I_{\vct{E}}$ and fluence
$F$ can be described by equations~\eqref{eq:intensity-1d}
and~\eqref{eq:fluence-1d} respectively. We note that the intensity is
proportional to the power and the fluence is proportional to the energy.
\end{document}